\documentclass[a4paper,11pt]{article}
\usepackage[english]{babel}

\usepackage{amsmath,amsfonts,amssymb, amsthm, amsbsy,amscd, mathrsfs}
\usepackage{enumerate, enumitem, xfrac, graphicx, graphics,epsfig, hyperref, fontenc, colortbl, bigints, comment}
\usepackage[utf8]{inputenc}
 \usepackage{AVDoc4}
 \usepackage{Muller2}

\numberwithin{equation}{section}

\newtheorem{prop}{Proposition}[subsection]
\newtheorem{theo}[prop]{Theorem} 
\newtheorem{defnt}[prop]{Definition} 
\newtheorem{lem}[prop]{Lemma}
\newtheorem{rem}[prop]{Remark}
\newtheorem{cor}[prop]{Corollary}
 
\begin{document}
 	\title{Metastability between the clicks of Muller's ratchet
}
\author{Mauro Mariani
	\footnote{Faculty of Mathematics, National Research University Higher School Of Economics, 6 Usacheva St., 119048 Moscow, Russia; email~: \texttt{mmariani@hse.ru}}, 
	Etienne Pardoux
	\footnote{I2M, Aix-Marseille Université, CNRS, Marseille, France;
	email~: \texttt{ etienne.pardoux@univ-amu.fr}}, 
	Aur\'elien Velleret
	\footnote{MaIAGE, INRAE,
		Universit\'e Paris-Saclay,
		F-78350, Jouy-en-Josas, France;
	email~: \texttt{ aurelien.velleret@nsup.org}, corresponding author}
}
\maketitle

\section*{Abstract}

We prove the existence and uniqueness 
of a quasi-stationary distribution
for three stochastic processes 
derived from the model of Muller's ratchet.
This model was invented with the aim of evaluating 
	the limitations of an asexual reproduction mode in preventing 
	the accumulation of deleterious mutations 
	through natural selection alone.
The main considered model
is non-classical,
as it is a stochastic diffusion 
evolving on an irregular set of infinite dimension 
with hard killing on an hyperplane.
We are nonetheless able to prove
exponential convergence in total variation to the  quasi-stationary distribution
even in this case.
The parameters in this last convergence result 
are directly related
to the core parameters of Muller's ratchet.
The speed of convergence to the quasi-stationary distribution
is deduced both for the infinite dimensional model 
and for approximations with a large yet finite number of potential mutations.
Likewise, we give uniform moment estimates 
of the empirical distribution of mutations
in the population under quasi-stationarity.
	
\section{Introduction}

\subsection{General presentation}
Since deleterious mutations occur much more frequently
than beneficial ones,
it is crucial to understand 
how the fixation of these deleterious mutations
is regulated.
Notably, 
it is very exceptional 
that a subsequent mutation reverts
a deleterious one,
so that only natural selection
can maintain some purity in the population.
In this respect, 
there is a major distinction to be made 
between sexual and asexual reproduction.
In a purely asexually reproducing population,  
a deleterious mutation can only be purged 
when the lineages carrying it go extinct.
In a sexually reproducing population, 
such a deleterious mutation can be avoided
through recombination,
without getting rid of the whole set 
of other mutations carried by the lineages.
There is actually no strong evidence 
that deleterious mutations
are specifically targeted
during this process of recombination.
For natural selection to effectively reduce the mutational load,
	it appears sufficient 
that this random process of recombination
 prevents some lineages 
from carrying the mutation after a given time.
This ability to better keep the population
purified from deleterious mutations
is one of the main explanations 
of the success of the sexual reproduction mode
(see \cite{MS78} for more details).
Such an advantage for sexual reproduction 
is to be confronted
with the cost (in terms of reproduction efficacy)
of requiring two parents.
The above scheme for purging deleterious mutations
in asexual populations
is the main object of study of the current paper.
\\

We plan to justify the existence and uniqueness 
of a metastable state
in which selective effects are able to maintain 
a subpopulation free from any deleterious mutation.
At the time where this subpopulation goes extinct,
we say that a click occurs in the population.
It has been shown in \cite{AP13}
that clicks happen
in finite time a.s. even for the 
infinite dimensional diffusion model 
(with  infinitely many types of individuals).
Rigorous definitions
of such a metastable state 
(characterized by the absence of click)
can be obtained in a broad generality
by a conditioning of stochastic processes.
We refer to Subsection \ref{M_sec_EQE}
for the definition 
of several crucial characteristic of metastability,
especially the notion of quasi-stationary distribution (abbreviated as QSD).

We treat in this paper three models of Muller's ratchet:
the first one is discrete both in time and space,
the second is a finite dimensional diffusion
and the third one is an infinite dimensional diffusion, 
see Section 1.2 below.
We prove the existence and uniqueness 
of a QSD
for those three stochastic representations
see respectively Theorems \ref{M_ECVdisc}, 
\ref{M_ECVFin} and \ref{M_ECVdInf}.
To our knowledge, 
the existence and uniqueness of a QSD
has not been rigorously proved until now 
except in the case of a finite state space.
This result was nonetheless implicitly exploited
for the approximations provided in \cite{ME13}.

We shall see that these QSD are concentrated on distributions 
with light tails, 
meaning that the proportion of the population 
carrying a large number of mutations
remains negligible under the QSD.
This claim is supported by our Proposition \ref{M_prop_MT}.
\\

We also address the classical issue 
of specifying the conditions 
under which  metastability is observed in practice.
A generally accepted answer 
is to compare the so-called relaxation time $t_R$,
which quantifies the rate at which the dependence in the past conditions vanishes,
and the average clicking time $t_C$ of the system.
Metastability between clicks would be the most common observation 
provided $t_R \ll t_C$,
so that a sequence of i.i.d. exponential law 
provides an accurate description
of the sequence of intervals between clicks.
This is where the comparison 
with the clicks of a ratchet comes from. 
If $t_R$ is of the same order as $t_C$ or larger,
we a priori can not exclude 
that trains of short interdependent intervals 
could alter this observed distribution of interval length.
But already if $t_R$ is of the same order as $t_C$,
there shall still be long realizations of inter-click intervals
after which we can say that the dependence in the past is forgotten.
This discussion is pursued in more details in Section~\ref{sec_Disc}.

The above mentioned theorems provide
a proper definition of these two main quantities.
The typical clicking time $t_C$
is defined as the inverse of the extinction rate of the QSD,
or equivalently as the expected waiting time of the next click 
with the QSD as the initial condition.
On the other hand,
 the QSD is approached
at an exponential rate
by the marginal law of the process conditioned 
upon the fact that the click has not occurred.
We describe the inverse of this exponential rate
(which should be independent of the initial condition besides the QSD)
as the typical relaxation time $t_R$.
\\

%

As compared to the other models 
that we have treated by similar techniques as in the current paper,
the proof of Theorem \ref{M_ECVdInf} is particularly difficult.
It specifically exploits the effect of selection
 to obtain practical bounds 
on the maximal number of accumulated mutations.
The argument is technical
 because at any time 
 an infinitesimal proportion
of heavily counter-selected mutants
cannot be completely neglected.

A simplified version of such bounds 
is already needed for the proof of Theorem \ref{M_ECVdisc}.
This concerns the process
defined in Section \ref{M_Haigh}.
The fact that the process describes a discrete population
 greatly simplifies the argument.
 We then extend the justification of the relaxation time 
 and the clicking rate
 for large population limiting models.
 Note that the results of \cite{ME13}
  or \cite{EPW09}
already largely exploit 
the fact that the population 
 is  large.
 In the diffusive limit defined in Section \ref{M_sec_DiffL},
 an additional difficulty arises
 in that the diffusion is degenerate
on a non-smooth boundary
 that is partly absorbing and partly repulsive.
In order to present a simplified analysis,
we  introduce in a first step a limitation 
in the number of carried mutations 
for the statement of Theorem \ref{M_ECVFin}
 given in Section \ref{M_sec_FinD}.
In the last step given in Section \ref{M_sec_InfD}
with Theorem \ref{M_ECVdInf},
we establish the existence and uniqueness of a QSD
for the more natural infinite dimensional model.
\\

The paper is organized as follows.
In the next Section \ref{M_sec_MMul}, 
we specify the stochastic processes under consideration, 
first the individual-based model in Section \ref{M_Haigh}
and then its diffusive limits in Section \ref{M_sec_DiffL}.
Our results of quasi-stationarity are presented in Section \ref{M_sec_TH}.
Starting in Section \ref{M_sec_EQE} with the general notion of exponential quasi-stationarity
that we aim to establish,
we treat respectively in Sections \ref{M_sec_QSDdisc}, \ref{M_sec_FinD} and \ref{M_sec_InfD}
each of the three stochastic processes mentioned above.
The generic assumptions and theorems on which these proofs rely
are stated in Subsection \ref{M_sec_SeCo}.
Next, we discuss more precisely the interpretation of these results 
in Section \ref{sec_Disc}.
We justify in Subsection \ref{sec_inter}
to \ref{sec_tRtC}
under which conditions quasi-stationarity can be observed.
Finally, we motivate our choice
for not introducing a bound on the number of mutations
in Section \ref{M_sec_MotUB}.
The rest of the paper is dedicated to the proofs.
Sections \ref{M_sec_ecvdisc}, \ref{M_sec_ECVFin} and \ref{M_sec_ecvdinf}
are devoted to the proofs of quasi-stationarity
for each of the three processes,
while Section \ref{M_sec_MT}
is devoted to uniform moment estimates of the QSDs.
The proofs in Section \ref{M_sec_ecvdinf}
take advantage of the lemmas and propositions 
derived in the previous Sections~\ref{M_sec_ECVFin}
and \ref{M_sec_MT}.

\subsection{The mathematical model of Muller's ratchet}
\label{M_sec_MMul}

\subsubsection{The individual-based model as a guideline}
\label{M_Haigh}
For the origin of the models which we study,
we refer to the simplified mathematical model 
which has been proposed by Guess as the multiplicative fitness model \cite{Gu74}
to quantify the regulation of deleterious mutations
in an asexual population.
The interest for this type of simplified models 
stems from general considerations 
on the evolutive advantage of recombination,
as notably advanced by Muller in 1964 \cite{M64}.
Since in any finite population, 
the ultimate fixation of deleterious mutations cannot be avoided
(unless by the extinction of the population),
yet a form of metastability can be observed.
This ``mechanism'' of regulation has been called Muller's ratchet,
notably by Haigh in \cite{H78}.

Assuming a constant deleterious effect of mutations, 
at each time that the fittest individuals disappear, 
the ratchet clicks in the sense 
that the new fittest individuals
carry an additional deleterious effect.
The whole population is doomed
	after this time to carry at least this 
	additional effect. 
Since the population size is constant,
natural selection then acts as if the whole profile of mutations
were translated by this value, 
so that the fittest individuals at that clicking time 
now become the new reference (at mutational burden 0).
If the mutation rate is slow enough
to allow these fittest individuals
to maintain the stability 
of the system for a while,
the dynamics shall rapidly follow the same behavior as
 before the click
 (taking into account that the empirical distribution
 of the number of carried mutations
 is translated).

This first model with discrete generations
and fixed population size $N$ evolves as follows.
Mutations that occur are only deleterious
and they occur at constant rate $\lambda >0$. 
The cost in fitness of each mutation
is quantified by $\alpha\in (0, 1)$.
Assume that the current population is distributed
with $N_i$ individuals carrying $i$ mutations
and consider an individual from the next generation.
Each one chooses its parent independently of the others
according to the same probability distribution,
which is specified by the fact that the chosen parent
carries $i$ mutations
with probability:
\begin{equation*}
	\dfrac{
N_i	(1 - \alpha)^i}{\sum_{k \ge 0} N_k\cdot (1 - \alpha)^k}.
\end{equation*}
Remark  that $\alpha = 0$ corresponds to neutral mutations,
	which are not under purifying selection.
If $\alpha = 1$ on the other hand, 
 no individual could survive the burden of a single mutation.

In addition to the mutations of its parent,
 each newborn gains $\xi$ deleterious mutations, 
where $\xi$ is a Poisson random variable 
with mean $\lambda$, 
specific to the newborn.
$\xi$ is drawn
independently for each newborn and 
of the choice of the parent.

Existence and uniqueness for such a discrete-time Markov chain 
on a countable state-space is a classical result,
which we shall take for granted.
\begin{rem}
Of course,
the situation is  more intricate in reality.
Mutations certainly do not have constant effect,
and combination effects are frequent (i.e. epistasis).
In many asexual populations,
there is evidence of the role of horizontal gene transfers,
for instance with plasmids 
(\cite{KP08}, \cite{MJ10}, \cite{OLG00}, \cite{TR+08}),
which can be seen as a weak form of recombination.
Moreover, the fact that mutations are deleterious 
is due to a change in the physiology 
that may be compensated by other means.
It might even happen that after subsequent mutations, 
the carriers of an initially deleterious mutation 
become more adapted than the wild types~\cite{SB+14}.
Neglecting these effects enables however 
to gain insight on the main regulatory factor.
\end{rem}

\subsubsection{The stochastic diffusion under consideration}
\label{M_sec_DiffL}
In the following,
we also consider a description of the model 
that corresponds to a limit of large population size,
accelerated time-scale 
(for which time is continuous),
thus also small selective effect
and small mutation rate.
In the following statements,
$d\in \II{1, \infty}$
(i.e. $d\in \N \cup \{\infty\})$ 
defines an upper-bound
on the number of deleterious mutations 
that can be carried by an individual.
If $d:=\infty$ in the following expression, $i\in \II{0, d}$ has to be understood as $i\in \Z_+$.

We are interested in the following Fleming-Viot system of 
Stochastic Differential Equations (SDEs)
for the $X\eD_i(t)$’s, $i \in \II{0, d}$,
where
$X\eD_i(t)$ denotes the proportion of individuals in the population 
who carry exactly $i$ deleterious
mutations at time $t$ (with $X\eD_{-1} \equiv 0$):
\begin{equation}
\begin{split}
	&
	\RMd X\eD_i(t) = \alpha\cdot (M\eD_1(t) - i)\cdot X\eD_i(t)\, \RMd t 
	+ \lambda\cdot (X\eD_{i-1}(t)-\idc{i<d} X\eD_{i}(t))\, \RMd t 
	\\&\hcm{2}
	+ \sqrt{X\eD_i(t)} \, \RMd W_i(t)
	- X\eD_i(t) \, \RMd W\sD(s).
\end{split}
\label{Sd}
\end{equation}
In \eqref{Sd},
$(W_i)_{i\ge 0}$ denotes a family of mutually independent standard Brownian motions,
where $W_i$ specifies the demographic fluctuations that are specific 
to the subpopulations $i$.
Secondly, the martingale process $W\sD$ is defined as follows:
\begin{equation}
	W\sD(s):=  \sum_{j=0}^{d} \int_0^t \sqrt{X\eD_j(s)} \RMd W_j(s),
	\label{eq_def_WsD}
\end{equation}
namely as an aggregated component 
according to which the fluctuations in the total population sizes are corrected.
Finally,
the process $M\eD_1$ is defined as follows:
\begin{equation*}
M\eD_1(t):=  \sum_{i=0}^{d} i\cdot X\eD_i(t),
\end{equation*}
namely as the aggregated component 
according to which the variations in the total population sizes 
due to the selective effects are corrected.
Unless otherwise specified, 
the Brownian motions that we introduce 
are all standard
and this precision will possibly be omitted.
\begin{rem}
	\label{rem_martRep}
The martingale term is described by the following representation:
	\begin{equation*}
		\RMd \cN\eD_i(t)
		:=\sqrt{X\eD_i(t)} \, \RMd W_i(t)
		- X\eD_i(t) \, \RMd W\sD(t)
	\end{equation*}
which is actually equivalent to another commonly considered representation, namely:
	\begin{equation*}
		\RMd \cN\eD_i(t) = \Tsum{j\neq i} \sqrt{X\eD_i(t) X\eD_j(t)} \, \RMd W_{i, j}(t)
	\end{equation*}
	for a sequence $(W_{i, j})_{i<j}$ of independent Brownian motions,
	extended to any $i\neq j$
	by the symmetry property $W_{i, j}(t) = W_{j, i}(t)$.
	As proved in Proposition~\ref{prop_Trep} of the appendix, 
	weak existence and uniqueness hold for both systems of SDEs
and	the two representations have actually the same law.
Since $\sum_{i=0}^{d} X\eD \equiv 1$,
and thanks to Lévy's caracterisation,
the martingale $W\sD$ has the law of a Brownian motion. 
Since 
\begin{equation*}
	\RMd \LAg \cN\eD_i\RAg_t = X\eD_i(t)\cdot (1-X\eD_i(t)) \, \RMd t,
\end{equation*}
there exists for any $i\in \II{1, d}$
a Brownian motion $B_i$ (which depends upon $d$) such that:
	\begin{equation}
	\RMd \cN\eD_i(t)
	:=\sqrt{X\eD_i(t)\cdot (1-X\eD_i(t))} \, \RMd B_i(t).
	\label{eq_def_Bi}
\end{equation}
However the $(B_i)$ are not mutually independent.
\end{rem}

In \cite{PSW12}, 
a closely related process with compensatory mutations is considered.
We refer to this article for a detailed presentation 
of the connection to related individual-based models
and only sketch next the interpretation of the parameters.

The selective effect of the deleterious mutations
is the term proportional to $\alpha$ in the drift term.
As we assume that all deleterious mutations carry the same burden
and that the total population size is fixed,
the growth rate of individuals 
carrying $i$ mutations is proportional 
to the difference between $i$ 
and the average number of mutations, i.e. $M\eD_1(t)$.
The appearance of new mutations is modeled 
by the term proportional to $\lambda$ in the drift term.
$\lambda$ corresponds to the rate at which individuals
carrying $i$ mutations 
give birth to individuals
carrying $i+1$ mutations.
Finally, the neutral choice of the individuals
replaced at each birth events gives rise to the martingale term.
Our time-scale 
corresponds to the rescaling of time $t \mapsto t'/N$, 
where $N$ is the  population size.


%

\subsubsection{Notations}
	
For $\cX$ a generic (Polish) space, hereafter $B(\cX)$
(respectively $B_+(\cX)$)
denotes the space of
bounded measurable 
(respectively nonnegative bounded measurable)
functions from $\cX$ into $\bR$ and $\M_1(\cX)$ (respectively $\cM(\cX)$)
the space of Borel probability (respectively signed) measures.
For any $f\in B(\cX)$ and $\mu \in \M(\cX)$, 
we use the abbreviation: 
$\LAg \mu \bv f\RAg:=\int_\cX f(x) \mu(\RMd x)$.
$\bR_+$ (respectively $\bR_+^*$) 
denotes the set of non-negative (respectively positive) reals.	
For any integers $m, n$ such that $m\le n$, $\II{m, n}$ denotes the set of integers from $m$ to $n$ (included).	We recall that $d=\infty$ is included in the expression $d\in \II{1, \infty}$.
$\Z_+$ (respectively $\N$) denote the set of non-negative (respectively positive) integers.

\section{Exponential quasi-stationarity results}
\label{M_sec_TH}

\subsection{Exponential quasi-stationarity}
\label{M_sec_EQE}

The conclusions of the following theorems
are  expressed in terms of the notion
of exponential quasi-stationarity that we borrow from \cite{AV_disc}.

\begin{defnt}
	\label{UEQS}
	For any linear, nonnegative, bounded and sub-conservative  semi-group $(P_t)_{t\ge 0}$
	acting on $\M(\cX)$, 
	we say that $P$ displays a uniform exponential quasi-stationary convergence
	(abbreviated as QSC)
	with characteristics
	$(\nu, h, \rho_0) \in \M_1(\cX)\ltm B_+(\cX)\ltm \bR_+$
	if $\LAg \nu \bv h\RAg = 1$
	and there exists $C, \gamma>0$ such that 
	the following inequality holds
	for any $t>0$
	and for any measure $\mu\in \M(\cX)$ with $\NTV{\mu}\le 1$:
	\begin{equation}
	\NTV{e^{\rho_0 t} \mu P_t- \LAg \mu\bv h\RAg\, \nu}
	\le C e^{-\gamma t}.
	\label{PtCV}
	\end{equation}
\end{defnt}

By the term of characteristics, we mean that the triple $(\nu, h, \rho_0)$
is uniquely defined,
as stated in \cite[Remark 2.7]{AV_disc}.
Thanks to  \cite[Corollary 2.9]{AV_disc}, 
this definition of convergence for $(P_t)$ implies the following convergence result to $\nu$.
\begin{cor}
	\label{D:CVAl}
	Assume \eqref{PtCV}.
	Then for any $t\ge 0$
	and  $\mu \in \M_1(\cX)$
	such that $\LAg \mu\bv \heta \RAg>0:$
	\begin{equation}
	\|\, \LAg \mu P_t, \idg{}\RAg^{-1}
	\cdot \mu P_t(\RMd x)
 - \nu(\RMd x) \, \|_{TV}
	\le C \dfrac{\|\mu - \nu\|_{TV}}
	{\LAg \mu\bv \heta \RAg} \; e^{-\gamma \; \tp}.
	\label{M_CVal}
	\end{equation}	
\end{cor}

\begin{rem}
	\label{rem_nuChar}
Choosing $\mu = \nu$ in \eqref{PtCV}, 
it is not hard to deduce the following relation:
\begin{equation}
\frlq{\tp\ge 0} 
\nu P_\tp = e^{-\rho_0\, \tp}\, \nu,
\label{M_LBz}
\end{equation}
and in particular $\LAg \nu, \idg{}\RAg = e^{-\rho_0\, \tp}$,
cf  \cite[Fact 2.7]{AV_disc}.
This relation is what characterizes 
$\nu$ as a QSD since it implies that for any $t\ge 0$, 
\mbox{$\LAg \nu P_t, \idg{}\RAg^{-1}
	\cdot \nu P_t(\RMd x)$} $=\nu(\RMd x).$
By restricting the convergence stated in \eqref{PtCV}
on the evaluation of the measure on $\cX$, 
we obtain a similar characterization of $h$. 
This latter convergence is what makes us call $h$ 
 the survival capacity.
\end{rem}

There is an additional related notion 
that will be useful to describe the behavior of the process
with the requirement of a long inter-click interval.
This process is generically defined
through the survival capacity $h$,
 on the state space:
$\mathcal{H}:= \{x\in \cX \pv h(x)>0\}$.
In the following, we assume that $(P_t)_{t\ge 0}$ 
is the sub-conservative semi-group 
of a strong Markov process $(\Omega; (\cF_\tp )_{\tp \ge 0};
(X_\tp)_{\tp \ge 0}; (\PR_{x})_{x \in \cX})$
defined up to its extinction time $\ext$, 
which is expresssed as follows for any initial distribution $\mu\in \cM_1(\cX)$:
\begin{equation*}
\mu P_t(\RMd x)
= \bE_{\mu}\big(X_t\in \RMd x;\, t<\ext\big).
\end{equation*}
Note that for such a semi-group $(P_t)$ displaying QSD
	with characteristics
	$(\nu, h, \rho_0)$
	(in $\M_1(\cX)\ltm B_+(\cX)\ltm \bR_+$),
$\rho_0>0$ is equivalent to $\PR_\nu(\ext<\infty)>0$,
	thanks to Remark \ref{rem_nuChar}.
	Furthermore, $\bE_{\mu}\big(X_t\in \RMd x\, \vert\, t<\ext\big) 
	= \LAg \mu P_t, \idg{}\RAg^{-1}
	\cdot \mu P_t(\RMd x)$ in this framework, 
	which converges to $\nu(\RMd x)$ according to \eqref{M_CVal}.

\begin{defnt}
We say that the Q-process exists 
if there exists a family $(\Q_{x})_{x \in \mathcal{H}}$ of probability
measures on $\Omega$ defined by:
\begin{equation}
\limInf{\tp }
\PR_{x}(\Lambda_\spr \bv \tp  < \ext) 
= \Q_{x}(\Lambda_\spr)
\end{equation}
for any $\cF_\spr$-measurable set $\Lambda_\spr$. 
We also require the process $(\Omega; (\cF_\tp )_{\tp \ge 0};
(X_\tp)_{\tp \ge 0}; (\Q_{x})_{x \in \mathcal{H}})$
to be a homogeneous strong Markov process. 
\end{defnt}

With a slight adaptation of the proof of \cite[Theorem 2.3]{AV_QSD},
it was deduced in \cite{AV_disc}
that such a property of existence of the Q-process
is a consequence of QSC.

\begin{rem}
The transition kernel of the Q-process
is given by: 
\begin{equation*}
q(x; \tp ; dy) 
= e^{\rho_0\, \tp } \, \dfrac{\heta(y)}
{\heta(x)}\; p(x; \tp ; dy),
\label{M_qt}
\end{equation*}
where $p(x; \tp ; dy)$ is the transition kernel 
of the Markov process $(X)$ under
\mbox{$(\PR_{x})$}.
Note that $\cX\setminus \mathcal{H}$ is 
generally avoided by the process $X$ 
under $\Q_x$.
In the examples of the current article, 
$h$ is actually positive
while $\nu$ is unique as a QSD. 
No distinction has then to be made between $\cX$ and $\mathcal{H}$ regarding the Q-process.

Thanks to  \cite[Corollary 2.12]{AV_disc},
our justification for the proof 
of QSC
actually implies related results of convergence for the Q-process.
Notably  $\beta(\RMd x):= \heta(x)\, \nu(\RMd x)$ 
is the unique invariant probability measure of this process.
\end{rem}

\begin{rem}
	\label{rem_Omega_path}
The probability space $\Omega$ is usually not made explicit.
	For the purpose of a following condition (namely Assumption $(A3_F)$
	as stated in Subsection~\ref{M_ElAs}),
	we will require  for the analysis of the diffusion processes $X\eD$
	that $\Omega$ is of path type.
Following \cite{Ka02} (to exploit Proposition 8.8 within),
it means that $\Omega$ 
	is the canonical space 
	of a strong Markov process $(\wht{X}_t)_{t\in \bR_+}$
	(with values on a Polish space $\wht \cX$),
	while the filtration $(\cF_t)$ is the one generated by the process $\wht X$.
	In other words, $\wht{X}$ is exactly the identity mapping,
	so that $\wht X_t$ agrees with the
	evaluation map $: \omega \mapsto \omega_t$ for $\omega$ on 
	the set of measurable functions from $\bR_+$ to $\wht{\cX}$, 
	and $(\cF_t)$ is the smallest filtration
	that make this application $\cF_t$-measurable for any $t\ge 0$.
	For our purposes,
	for any $d \in \II{1, \infty}$,
	it will be sufficient to consider for  $\Omega\eD$
	and for the associated filtration $(\cF\eD_t)$
	the canonical choice generated 
	by the process $X\eD$. 
\end{rem}

\begin{rem}
Originally in the quantitative work 
of Champagnat and Villemonais,
the focus is more on estimates of convergence for the conditioned semi-group,
as in \eqref{M_CVal}.
There was notably an efficient and full characterization 
of the case where the upper-bound can be taken as $C e^{-\gamma t}$
uniformly over the initial condition \cite{CV16}.
In a more recent work \cite{CV23}, 
they also relate to the convergence of the semi-group scaled by $e^{\rho_0 t}$
(\cite[Corollary 2.4 and 2.6]{CV23})
and describe it as naturally adapted as well 
(\cite[Remark 5]{CV23}).
As in \cite{AV_disc},
we found the formulation \eqref{M_CVal}
more generally adapted to the dependency over the initial condition
due to the linearity of the semi-group $P_t$.
A more general form of convergence than \eqref{PtCV}
with a weighted norm 
is the object of \cite{BCGM22},
with a full characterization as well.
\end{rem}

\subsection{The discrete population case}
\label{M_sec_QSDdisc}

Let $N\ge 1$ be the number of individuals in the population.
For  $n\le N$ and $t\in \Z_+$,
let $D_n(t)$ be the number of mutations carried by the $n$-th individual.
We consider the empirical measure at time $t>0$
 defined as follow:
\begin{equation}\label{M_ZN}
\cZ^N_t:= \frac{1}{N} \sum\nolimits_{n=0}^N \delta_{D_n(t)},
\end{equation}
so that $\cZ^N_t(i)\in N^{-1}\cdot[\![0, N]\!]$ specifies the proportion of individuals with exactly $i$ mutations
(since everything is discrete, we identify $\cZ^N_t$ as a function from $\Z_+$ to $\bR_+$).
From the rules describing the next generation from the previous one, 
see Subsection 1.2.1, 
$\cZ^N$ is a strong Markov process evolving on $\M_1^N(\Z_+)$, 
where:
\begin{equation}\label{M_MUNZ}
\begin{split}
	\textstyle
\M_1^N(\Z_+)&:= \Lbr\frac{1}{N}\cdot \sum_{i\le N} \delta_{d_i}\pv d_i \in \Z_+
\mVg\sum_{i\in \Z_+} d_i = N\Rbr
\\
&\quad \equiv \Big\{z: \Z_+\mapsto N^{-1} \times [\![0, N]\!]\pv \Tsum{i\in \Z_+} z(i) = 1\Big\}.
\end{split}
\end{equation}
This set $\M_1^N(\Z_+)$ is equipped with the Lévy-Prokhorov metric, 
which makes it a Polish space, that is a separable and complete metric space.
The clicking time under consideration comes from the extinction of the fittest individuals, i.e.:
\begin{equation}\label{M_extDisc}
\ext^N:= \inf\Big\{t\ge 0\pv \cZ^N_t(0) = 0\Big\}
= \inf\Big\{t\ge 0 \pv \cZ^N_t \notin \M_1^{(0), N}(\Z_+)\Big\},
\end{equation}
where $\M_1^{(0), N}(\Z_+) = \Big\{z\in \M_1^N(\Z_+)\pv z(0)\ge \frac{1}{N}\Big\}.$

\begin{rem}
Classical theory on quasi-stationarity
can be exploited by interpreting the clicking time $\ext^N$
as an extinction time.
We implicitly rely on the process $\bar{\cZ}^N_t 
:= \idc{t<\ext^N}\cZ^N_t  + \idc{\ext^N \le t}\partial $
which is Markov 
and lives on $\M_1^{(0), N}(\Z_+)\cup\{\partial\}$.
For the process $\bar{\cZ}^N$, 
$\ext^N$ is the hitting time 
of the absorbing state $\partial$ 
(the cemetery).
\end{rem}

Our main conclusion for this process is the following theorem:
\begin{theo}
	\label{M_ECVdisc}
	Consider for any $\alpha \in(0, 1)$, for any $\lambda >0$ and for any $N\ge 1$ the Markov process $Z^N$
whose transitions are prescribed as in Section \ref{M_Haigh},
	with  clicking time $\ext^N$.
Then, its semigroup $P^N$ 
displays QSC
with characteristics
$(\nu^N, h^N, \rho_0^N) \in \M_1(\M_1^{(0), N}(\Z_+))\ltm B_+(\M_1^{(0), N}(\Z_+))\ltm \bR_+^*$.
Moreover, $h^N$ is uniformly bounded away from 0.

It implies that the  convergence to $\nu^N$  given in \eqref{M_CVal} is uniform
with respect to the initial condition
and that the Q-process exists on the whole state space $\M_1^{(0), N}(\Z_+)$.
\end{theo}

\begin{rem}
The discreteness of the process 
is strongly involved in the proof of Theorem~\ref{M_ECVdisc},
so that the derived convergence rate
strongly depends on $N$.
\end{rem}

\subsection{The finite dimensional case}
\label{M_sec_FinD}

The system of SDEs for finite $d$ evolves on the state space $\cX_{d}$, where:
$$\bar \cX_{d}:= 
\Big\{(x_k)_{k\in \II{0, d}}\in [0, 1]^{d+1}\pv
 \sum\nolimits_{k = 0}^d x_k = 1\Big\}.
$$
With the $L^1$ distance, this set is also a Polish space.
The solution to the system \eqref{Sd} is unique as stated in the next proposition, 
as a corrolary of \cite[Theorem 2.1]{Shi87},
whose proof is deferred for completeness to the appendix, Section~\ref{app_TwoRep}.
\begin{prop}
	\label{prop_ex_uniq}
For any $d\in \N$,
 existence and weak uniqueness 
 of solutions
 hold on the state space $\bar \cX_d$
 for the system \eqref{Sd} of SDEs. 
\end{prop}
We denote by  $\ext\eD$ the clicking time 
of this process $X^{(d)}$, 
that is:
\begin{equation}
	\ext\eD:= \inf\Big\{t\ge 0\pv X^{(d)}_0(t) = 0\Big\}.
	\label{def_extD}
\end{equation}
The most natural state space for the process with extinction at time $\ext\eD$
thus does not include this boundary:
\begin{equation}
\cX_d := \Big\{x\in \bar \cX_d; x_0\in (0, 1]\Big\}.
\label{eq_def_Xd}
\end{equation}

Our main conclusion for this process is the following theorem:
\begin{theo}
\label{M_ECVFin}
Consider the system of SDEs
$\eqref{Sd}$ 
for any $\alpha \in(0, 1)$, for any $\lambda >0$ and
for any $d\in \N$,
with  clicking time $\ext\eD$.
Then, its semigroup $P\eD$ 
displays a QSC
with characteristics
$(\nu\eD, h\eD, \rho\eD_0) \in \M_1(\cX_d)\ltm B_+(\cX_d)\ltm \bR_+^*$.
In addition, 
for any $y_0\in (0,1)$,
$h\eD$ is bounded away from 0
on $\Lbr x\in \cX_d\pv x_0 \ge y_0\Rbr$.
In particular, the associated Q-process exists
on the whole state space $\cX_d$.
\end{theo}

\begin{rem}
The fact that $\alpha$ is non-zero is actually not exploited in this proof.
As in \cite{CV21}, 
we rely mainly on Harnack's inequality.
The dependency in the dimension $d$ is roughly considered.
Nonetheless, 
we have here to be cautious in the way 
we handle jointly the absorbing 
and repulsive boundary conditions.
\end{rem}

Moreover, 
we prove the following controls on the moments of the QSDs $\nu\eD$, for $d\in \N$:

\begin{prop}
	\label{M_prop_MT}
For any $\alpha \in(0, 1)$, for any $\lambda >0$ and 
for any $k\ge 1$,
	we have uniform tightness in $d$ 
	over the moments of order $k$ 
	of the unique QSDs $\nu\eD$
	associated with the solution to $\eqref{Sd}$,
	which means that the following supremum tends to 0
	as $m$ tends to infinity:
	\begin{equation*}
	\textstyle 
	\sup\Big\{\int_{\cX_d}  \idc{M_k(x) \ge m}\nu\eD(\RMd x)\pv 
	d\in \N\Big\},
	\end{equation*}
where $M_k(x):= \sum_{i\in \II{0, d}} i^k\, x_i$.
	In particular,
	the sequence $\hat{\nu}_k\eD$,
	where the values for the coordinates larger 
	than $k+1$ 	are put to 0,
	is tight in $\M_1(\bR^\Z_+)$.
\end{prop}

\begin{rem}
	This control on the moments is actually crucial for the proof of uniqueness when $d=\infty$.
Given the above theorem,
we expect the sequence $(\nu\eD)$ to converge 
as $d\ifty$ to 
the unique QSD $\nu\eI$ 
for the infinite system (for which the control extends),
though it is beyond the scope of the current paper.
\end{rem}

\subsection{The infinite dimensional case}
\label{M_sec_InfD}

We consider now the infinite dimensional case,
for which we require the existence of moments. Let 
us consider the following definition for any $\eta\in (2, \infty)$:
\begin{equation*}
\bar	\cX^{\eta}:= 
\Lbr(x_k)_{k\in \Z_+}\in [0, 1]^{\infty}\pv
 \sum_{k = 0}^{\infty} x_k = 1\mVg \sum_{k = 0}^{\infty} k^{\eta} x_k < \infty\Rbr.
\label{def_cXeta}
\end{equation*} 

As in \cite{AP13},
we consider for $\bar \cX^\eta$ the topology under which
a probability $x^n = (x^n_k, k \ge 0)$ on $\Z_+$
converges to $x = (x_k, k \ge 0)$ if both it converges weakly, 
and $\sup_{n} \sum_{k\ge 0} k^\eta x_k^n < \infty$.
This topology is actually generated by the following metric,
defined for any $x, y\in \cX^\eta$,
which makes $\bar \cX^\eta$ Polish, that is separable and complete:
$$d_\eta(x, y) = |x_0-y_0| + \Tsum{k\ge 1} k^\eta\cdot |x_k-y_k|.$$
Thanks to \cite[Theorem 3]{AP13}
(and to Proposition~\ref{prop_Trep} in the appendix to show that
the process under consideration is the same as ours),
we know that for any $\eta > 2$
and for any initial condition $x$ that belongs to $\bar\cX^{\eta}$, 
\eqref{Sd} has a unique
weak solution $X\eI$
which is a.s. continuous with values in $\bar \cX^{\eta}$.
This process has been introduced in \cite{EPW09}
and it has also been shown in \cite{AP13}
that clicks occur a.s. in finite time, 
with the same definition of $\ext\eI$ as in \eqref{def_extD}.
We consider the state space without the boundary $\{x_0= 0\}$
as the state space with $\eta=6$ for the process with extinction at time $\ext\eI$:
\begin{equation}\label{eq_def_Xinf}
	\cX_\infty := \{x\in \bar \cX^6; x_0\in (0, 1]\}.
\end{equation}

We now state the main theorem of the current paper:
\begin{theo}
	\label{M_ECVdInf}
	Consider for any $\alpha \in(0, 1)$ and for any $\lambda >0$ 
	the system of SDEs \eqref{Sd}
with  $d = \infty$,
	defined on $\cX_\infty$ 
	with  extinction at time $\ext^{(\infty)}$.
	Then, its semigroup $P^{(\infty)}$ 
	displays a QSC
	with characteristics
	$(\nu\eI, h\eI, \rho_0\eI) \in \M_1(\cX_\infty)\ltm B_+(\cX_\infty)\ltm \bR_+^*$.
		In addition, 
	for any $y_0 \in (0,1)$,
	$h\eI$ is lower-bounded by a positive constant
	on $\Lbr x\in \cX_\infty\pv x_0 \ge y_0\Rbr$.	
In particular, the Q-process exists on the whole state space $\cX_\infty$.


	Besides, there exist $C, \gamma, d_\vee>0$ 
	such that
the convergences stated in \eqref{PtCV} and \eqref{M_CVal}
hold true with the constants $C, \gamma$ 
for the processes given by \eqref{Sd} 
on $\cX_d$
for any $d\in \II{d_\vee, \infty}$.
\end{theo}

By the notation, $\II{d_\vee, \infty}$,
remark that we include $\infty$ besides all integers larger than $d_\vee$.

\begin{rem}
The core of our proof is based on the intuition 
that the slower the decay of the tail
the more rapidly it gets erased and renewed.
So we do not expect large tails to play a significant role.
In practice, we exploit the finiteness of moments of order $2 \eta'$
to control moments of order $\eta'$,
and that $\eta'$ is strictly larger than 2
plays a significant role.
For simplicity, we restrict ourselves to 
initial condition within $\cX^{6}$ 
although our proof could likely be generalized
to $\cX^\eta$ 
provided at least that $\eta >4$.
\end{rem}

\subsection{Some crucial sets of conditions ensuring exponential quasi-stationarity}
\label{M_sec_SeCo}

The proof of QSC
are exploiting the criteria 
given in  \cite[Subsection 2.3.1]{AV_disc},
with a trajectorial approach. 
The methods and statements taken from \cite{AV_QSD}  
have actually been adjusted in \cite{AV_disc}
 with the current paper in mind.
Thanks to  Theorem \ref{M_th:ECV}
in Subsection \ref{M_genSt}
that summarizes these conclusions,
it will remain to prove one of the following two sets of assumptions,
$\mathbf{(A)}$ or $\mathbf{(AF)}$,
to complete
the proofs of 
Theorems \ref{M_ECVdisc},
\ref{M_ECVFin} and
\ref{M_ECVdInf}.
Assumptions $\mathbf{(A)}$ and $\mathbf{(AF)}$
 are made up of four basic assumptions,
 three being common to both.
 So first, we present 
 these five basic assumptions
 in Subsection \ref{M_ElAs}
 in the general context of a strong Markov process $X$
 with extinction at time $\ext$
 that makes it exit the state space $\cX$. 
 
$\cX$ can here simply be thought to be a Polish space.
	The state spaces we consider for our theorems, 
	namely $\M_1^{(0), N}(\Z_+)$
	for any $N\ge 1$ (see \eqref{M_extDisc}),
	$\cX_d$ for any $d\ge 1$ (see \eqref{eq_def_Xd}),
	and $\cX_\infty$ (see \eqref{eq_def_Xinf}),
are all Polish spaces.

\subsubsection{Basic assumptions}
\label{M_ElAs}

The first assumption is on a sequence $(\cD_\ell)_{\ell\ge 1}$
that shall be exploited for the following assumptions.
 $int(\cD)$ there denotes the interior of any set $\cD$.
\\

\noindent
$(A0)[(\cD_\ell)]$:  \textbf{Specification on the state space}
with the sequence $(\cD_\ell)_{\ell\ge 1}$:
\begin{center}
	\begin{minipage}{0.9\textwidth}
		\textsl{For any $\ell$, it holds both that $\cD_\ell$
is a closed subset of $\cX$
and that $\cD_\ell \subset 	int(\cD_{\ell+1})$.}
\end{minipage}
\end{center}

\begin{rem}
Originally in \cite{AV_QSD},
this assumption is strengthened with $\medcup_{\ell\ge 1} \cD_\ell= \cX$.
It was shown in \cite{AV_disc}
how to adapt the conclusions 
without this additional restriction.
The state spaces $\cX$ that we consider
include boundaries
in the vicinity of which the process $X\eD$ is a degenerate diffusion
(for instance the boundary corresponding to $x_1 = 0$).
It is thus more convenient for some of our proofs 
not to assume $\medcup_{\ell\ge 1} \cD_\ell= \cX$
and if possible to keep the diffusion process non-degenerate on any $\cD_\ell$.
\end{rem}

The sequence $\cD_\ell$ will serve as a reference in the following, 
and we also denote:
\begin{equation}
	\mathbf{D}:= \big\{ \cD\subset \cX\pv \overline{\cD} = \cD\mVg  
	\Ex{\ell\ge 1} \cD \subset \cD_\ell\big\},
	\label{D:Dps}
\end{equation}
where $\overline{\cD}$ denotes the closure of the subset $\cD$ of $\cX$, 
so that the elements of $\mathbf{D}$ are closed subsets of $\cX$
that are countained in $\cD_\ell$, for $\ell$ sufficiently large.

For this trajectorial approach, 
we strongly rely on the representation of the semi-group $(P_t)$
in terms of a strong Markov process $(X_t)_{t\in [0, \ext)}$
defined up to some extinction time $\ext$.
Recall $P_t f(x)
= \bE_x[f(X_t);\, t<\ext]$
for any $x\in \cX$ and $f\in B(\cX)$.
For the next statements, we will exploit the following notations 
for the exit and entry times of any subset $\cD$ of $\cX$:
$$T_{\cD}:= \inf\Lbr  \tp \ge 0 \pv X_\tp \notin \cD \Rbr
\mVg\quad
\tau_\cD:= \inf\Lbr \tp \ge 0\pv X_\tp \in \cD \Rbr.$$
Besides this dependency on the sequence $(\cD_\ell)$
the next assumptions share as commom parameters
the probability measure $\zeta$ on $\cX$,
the real number $\rho>0$ and the set $E \in \mathbf{D}$.
They are recalled in square brackets in the notations of the basic assumptions.
\\

\noindent
$(A1)[\zeta, (\cD_\ell)]$ : \textbf{Mixing property}
for the reference probability measure $\zeta \in \cM_1(\cX)$
according to the sequence $(\cD_\ell)_{\ell\ge 1}$:
\begin{center}
\begin{minipage}{0.9\textwidth}
\textsl{For any integer $\ell\ge 1$, 
	there exist both an integer $L>\ell$ and  $c, t>0$ such that
	the following holds
	for any initial condition $y \in \cD_{\ell}$:}
	\begin{equation*}
		\PR_y \Big( {X}_{\tp}\in \RMd x \pv
		\tp < \ext \wedge T_{\cD_L}\Big)
		\ge \cp\; \alc(\RMd x).	
	\end{equation*}
\end{minipage}
\end{center}

\noindent
$(A2)[\rho, E]$ : \textbf{Escape from the Transitory domain}
with penalty-rate $\rho>0$,
where the complementaty of this transitory
domain is $E\in \mathbf{D}$:
\begin{equation*}
	\Tsup{x\in \cX} \;\bE_{x} \big(
	\exp\big[\rho\cdot (\ext\wedge \tau_E) \big] \big) < \infty.
\end{equation*}
$\rho$ in the previous exponential moment
is required to be strictly larger
than the following \textbf{``survival estimate"}
$\rho_S$, which a priori depends on $\zeta$: 
\begin{equation*}
	\rho_S[\zeta]
	:= \sup\Big\{\gamma \ge 0;\, 
	\textstyle
	\Tsup{L\ge 1} \liminf_{\{t>0\}} \;
	e^{\gamma t}\,\PR_\alc\big(t < \ext\wedge T_{\cD_L}\big) 
	= 0
	\Big\}\vee 0.
\end{equation*}

\begin{rem}
It is proved in  \cite[Theorem 2.16]{AV_disc} 
that $\rho_S[\zeta]$ coincides with the extinction rate $\rho_0$
(and is thus independent of $\zeta$),
provided the semi-group displays QSC.  
\end{rem}


The next two  assumptions are proposed as alternatives
and each alternative will be exploited in the current paper.
The former is the assumption first introduced in \cite[Section 2.1]{AV_QSD}.
The latter provides a way 
to ensure the former given $(A0)$, $(A1)$ and $(A2)$ 
as proved in \cite[Theorem 2.3]{AV_disc}.
\\
$(A3)[\zeta, E]$~: \textbf{Asymptotic comparison of survival}
on the set $E\in \mathbf D$ of initial conditions,
with reference probability measure $\zeta$:
\begin{equation*}
		\underset{t\ifty}{\limsup} \;
		\underset{x\in E}{\sup} \;
		\dfrac{\PR_{x} (\tp< \ext)}
		{\PR_{\alc} (\tp < \ext)} 
		< \infty.
	\end{equation*}
\\
$(A3_F)[\zeta, \rho, E]$~: \textbf{Almost perfect harvest}
on the set $E\in \mathbf D$ of initial conditions,
with reference probability measure $\zeta$
and penalty-rate $\rho$:
\begin{center}
	\begin{minipage}{0.9\textwidth}
		\textsl{For any $\epsilon\in (0,\, 1)$,
	there exist $t_F, \cp >0$ 
	such that for any $y \in E$
	there exist two stopping times
	$U_H$ and $V$
with the following properties:} 
\begin{equation*}
	\PR_{y} \big(X(U_H) \in \RMd x \pv U_H < \ext \big) 
	\le \cp \,\PR_{\alc} \big(X(V) \in \RMd x
	\pv V < \ext\big).
\end{equation*} 
\textsl{
	including the next conditions on $U_H$:}
 \begin{equation*}
 	\Lbr\ext \wedge t_F < U_H \Rbr
 	= \Lbr U_H = \infty\Rbr
 	\quad  \text{ and }\quad
 	\PR_{x} (U_H = \infty, \,  t< \ext) 
 	\le \epsilon\, \exp(-\rho\, t_F).
 \end{equation*}
\textsl{Furthermore, $\Omega$ is of  path type.}
\end{minipage}
\end{center}

As stated in \cite[Proposition 2.2]{AV_disc},
the assumption that $\Omega$ is of path type 
ensures sufficient regularity properties of stopping times 
with respect to iterative procedures exploiting the strong Markov property.

\subsubsection{General theorems of convergence}
\label{M_genSt}

We exploit the following definitions of Assumptions $\mathbf{(A)}$ and $\mathbf{(AF)}$
from \cite{AV_disc}.
\begin{description}
	\item[\textbf{Assumption $\mathbf{(A)}$:}]
``There exists $\zeta\in \mathcal{M}_1(\cX)$ such that $(A1)[\zeta, (\mathcal{D}_\ell)]$ holds 
for a specific sequence $(\mathcal{D}_\ell)$ satisfying $(A0)$.
Moreover, there exists $\rho  > \rho_{S}[\zeta]$ 
and $E \in \mathbf{D}$ such that
Assumptions $(A2)[\rho, E]$
and $(A3)[\zeta, E]$  hold."
\item[\textbf{Assumption} $\mathbf{(AF)}$] is exactly Assumption $\mathbf{(A)}$
with $(A3)[\zeta, E]$ replaced by $(A3_F)[\zeta, \rho, E]$.
\end{description}
\noindent
Theorem  2.8
in \cite{AV_disc}
can be restated for our purpose as follows:
\begin{theo}
	\label{M_th:ECV}
	Assume that either $\mathbf{(A)}$ or $\mathbf{(AF)}$ holds.
	Then, the semigroup $P$ 
	displays QSC
	with characteristics
	$(\nu, h, \rho_0) \in \M_1(\cX)\ltm B_+(\cX)\ltm \bR_+$,
	$h$ is bounded away from $0$ on $\cD_\ell$ for any $\ell\ge 1$
	and the Q-process exists on $\cH:= \{x\in \cX \pv h(x)>0\}$.	
\end{theo}

\begin{rem}
	\label{rem_th_ECV}
The proof of Theorem \ref{M_th:ECV}
is originally stated in continuous time.
A careful check of the arguments given 
 in \cite[Section 3]{AV_disc} 
and in \cite[Section 3]{AV_QSD} 
shows that their extension to the discrete-time setting 
(exploited in the following Theorem \ref{M_ECVdisc})
does not raise any issue.
\end{rem}

 Since the exploited sequence $(\cD_\ell)_{\ell\ge 1}$ 
 usually does not cover the whole state space,
we shall exploit \cite[Proposition 2.10]{AV_disc}
to deduce some lower-bounds of $h$.
The next proposition recalls its statement.
\begin{prop}
	\label{th:UCV}
	Assume that $\mathbf{(A)}$ or  $\mathbf{(AF)}$  holds.
	Then, the survival capacity $h$ 
	is uniformly lower-bounded on any set $H\subset  \cX$
that satisfies the following condition:\\
	$(H_0) :$	there exists $t>0$, $\ell\ge 1$ such that 
	$\quad 
	\Tinf{x\in H}  
	\PR_x(\tau_{\cD_\ell} \le t \wedge \ext) > 0.$
	\\
	It implies the following identification of $\cH:= \{x\in \cX \pv h(x)>0\}$:
	$$\cH = \Big\{ x \in \cX\pv \Ex{\ell\ge 1}\PR_x(\tau_{\cD_\ell} < \ext)>0\Big\}.$$
\end{prop}

Thanks to Proposition \ref{th:UCV},
we shall prove in our models that $h$ is actually positive.
Note that this property
entails the uniqueness of the QSD thanks to Corollary \ref{D:CVAl}.

\section{Outlook}
\label{sec_Disc}

\subsection{Interpretation of crucial parameters}
\label{sec_inter}

For the infinite-dimensional process, no parameter 
other than $\alpha$ and $\lambda$ is introduced.
We deduce from Theorem \ref{M_ECVdInf}
that the QSD and the survival capacity depend only on $\alpha$ and $\lambda$,
as well as the real numbers $C, \gamma>0$ in \eqref{PtCV} and \eqref{M_CVal}.

As already noted by Haigh in \cite{H78},
$\alpha/\lambda$ is the average number of deleterious mutations
that are established in the deterministic limit 
(neglecting neutral fluctuations).
The deterministic distribution of mutations
is a function of $\alpha/\lambda$,
and actually follows a Poisson distribution with this mean,
as shown in \cite{EPW09}.
To infer the level of fluctuations 
around this deterministic equilibrium,
we shall look at the coefficient 
in front of the martingale term 
in a new time-scale such that 
the mutation rate is set to 1.
This gives $1/\lambda$,
which we recall to scale as $\sqrt{1/N_e}$, 
where $N_e$ is the population size.
A large population size thus corresponds 
to letting $\lambda$ go to infinity,
making  the deviations away 
from the deterministic distribution more unlikely.
\\

A natural scale for the time between clicks 
can be directly derived 
from the notion of QSC, 
with the definition $t_C:= \rho_0^{-1}$.
On the other hand, 
we can propose the following definition 
for the relaxation time:
\begin{multline}
t_R:= \inf\Big\{t_r>0\pv 
\Ex{C>0} \frl{\mu\in \cM_1(\cX)}
\frlq{t>0} 
\\ \NTV{\mu A_t - \nu} \le
(C/ \LAg \mu \bv \heta\RAg) \cdot e^{-t/t_r}  \Big\}.
\label{M_tR}
\end{multline}
Our results justify that this definition,
which involves no ad-hoc parameters,
leads to a finite quantity
(upper-bounded by $1/\gamma$,
where $\gamma$ is the rate deduced 
in the proof of Theorem \ref{M_ECVdInf}).
Remark also that the convergence to $\heta$ and $\beta$ 
also occurs at a  rate that is quicker than $1/t_R$,
as one can check by adjusting the justification in \cite{AV_QSD}.

By relying on the arguments 
of Theorem \ref{M_ECVdInf} and Proposition \ref{M_prop_MT},
we expect that truncating the number of accumulated mutations 
is not likely to alter much this value of $t_R$ 
provided the threshold is sufficiently large.
Since we cannot evaluate $t_R$ precisely 
and are only able to provide an upper-bound,
this is still conjectural.
But substantial increase of these last components 
are proved to be rare thanks to Proposition \ref{M_prop_MT}
and not so significant 
when we look at Section \ref{sec_AFdI}.

\begin{rem}
The dependence on the initial condition in \eqref{M_tR}
is expected from the linearity of the semi-group $(P_t)$, 
as observed in \cite{AV_disc}.
More general dependencies could nonetheless be imagined,
relying for instance on Lyapunov functions
 as in \cite{CV23} or in \cite{BCGM22}.
We simply do not think it would change the value of $t_R$
because the confinement is mainly due to extinction
and immediate repulsion from the boundaries.
\end{rem}

\subsection{Previous estimations}

The study of this quasi-stationary regime 
arises naturally when one wishes to estimate the rate
at which the ratchet clicks.
To obtain quantitative estimates,
several authors have justified their approach 
by assuming that the typical clicking time $t_C$
is much larger than 
the typical relaxation time $t_R$ of the system,
usually with an empirical reference for the latter
(\cite{EPW09}, \cite{ME13}).

In \cite{ME13},
an estimation of $t_C$ in the context where $t_R \ll t_C$
is obtained through the characteristic equation of a certain QSD $\nu_\star$,
of the form $\cL \nu_\star = -\lambda\cdot \nu_\star$, 
with $\cL$ is a certain infinitesimal generator 
and $\lambda$ its main eigenvalue.
This QSD $\nu_\star$ that they study 
is not the general QSD $\nu\eI$ that we describe.
The detailed description of the latter 
is reasonably argued to be too intricate.
The former is in fact derived
from a one-dimensional approximation 
of the process under metastability.
It is argued that in the context of large populations,
and given the number of fittest individuals,
one can approximate the rest of the distribution 
as an almost deterministic profile.
The dependence in this number of fittest individuals
only occurs in the normalizing factor of this distribution.
This latter argument of concentration 
could probably be made  rigorous  
by using Large Deviation theory. 
Such results are  beyond the scope of the present paper.

Note that the validity of this approximation 
relies upon the fact that $t_R\ll t_C$,
where $t_R$ is to be related to the QSD $\nu\eI$.
The relaxation rate of $\nu_\star$ is only a partial indication,
although presumably carrying most of the information.

\subsection{The quasi-stationary regime 
	is generally observed for $t_R \ll t_C$}
\label{sec_tRtC}
Provided $t_R \ll t_C$,
we  expect to generally observe 
the quasi-stationary regime between clicks.
It is classical that with the QSD as an initial condition,
the extinction time and extinction state are independent,
the former being exponentially distributed,
as it is stated in  \cite[Theorem 2.6]{CMS13}.
Assuming that we start the analysis at a new click
after a long time-interval without click,
it implies that the profile of mutations just after the click
is distributed as the restriction of the QSD 
to the hyperplane $\{X_0 = 0\}$.

Since having large values of $M_1$ 
makes it actually harder for the process to reach the hyperplane,
we expect that,
under the QSD restricted to $\{X_0 = 0\}$,
$M_1$ tends to be  smaller
than the prediction $1+\lambda/\alpha$
derived from the deterministic limit
(under the constraint that $x_0 =0$).
Besides, the fittest individuals are 
altered by first changing 
into the type with only one mutation.
So we expect also that under the QSD restricted to $\{X_0 = 0\}$, 
there is an over-representation 
of the proportion 
of individuals carrying a single mutation
(the new optimal trait).
Thus, 
we expect the distribution just after the click 
to be less prone to a future click than would be the QSD itself.
Since $t_R \ll t_C$, 
the quasi-stationary regime is then rapidly reached.

Let us also imagine a dramatic situation 
where some clicks would rapidly follow each others.
Then, 
it would imply that these fittest classes of individuals
are rapidly wiped off, 
while not letting much time for the others to change much.
Since we have seen that we have very strong controls of moments 
under the QSD, cf notably Proposition \ref{M_prop_MT},
such succession of clicks cannot hold for long.
A class that is not prone to a quick extinction
should be reached quite early
and generate a new quasi-stationary regime.
Such dramatic situations,
which are very rare,
are thus expected to be very isolated
and of limited impact.

Expecting an exponential law for the inter-click intervals
and the independence between them 
should be in conclusion a good approximation provided $t_R \ll t_C$.
\\

As we discuss in \cite[Subsection 2.3]{AV_Ada},
one can also conclude whether or not 
the QSD profile is likely to be observed without conditioning
by comparing $\nu$ to the survival capacity $\heta$.
If quasi-stationarity is stable, 
we do not expect that the conditioning 
on having a click in the far future
shall substantially alter the dynamics. 
In most trajectories, 
the Q-process shall thus behave as the original process.
So $\heta$ should be mostly constant on the support 
of $\beta(\RMd x) = \heta(x)\,\nu(\RMd x)$,
implying $\heta \approx 1$ 
where the density of $\nu$ is large. 

In practice, the QSD and the survival capacity
are certainly quite difficult to specify with simulations
because they live on a large dimensional space.
Likewise,
the convergence in total variation exploited in \eqref{M_tR}
is probably not very practical for numerical estimation.

\subsection{Motivation for an unbounded number of deleterious mutations}
\label{M_sec_MotUB}
In order to prove quasi-stationarity results,
the case where $d<\infty$ can be treated more easily
and provides an introduction to the case $d=\infty$.
Nonetheless,
the arguments for having a convergence at a given rate 
becomes more and more artificial as $d$ tends to infinity.
The constant involved in the Harnack inequalities
goes to zero as the dimension increases.
By considering the case $d=\infty$,
we actually handle as a whole the case where $d$ is sufficiently large.
By these means,
we are able to prove that 
the rate of convergence can be upper-bounded
by a quantity that does not depend on the specific value of $d$.
This is to be expected since, even when a large number of deleterious mutations is permitted,
we expect individuals 
carrying a large number of mutations 
to remain negligible.

Referring for instance to \cite{EPW09}, 
it is not difficult to prove that 
in the deterministic limit of a large population,
the empirical measure of the number of mutations in the population
tends to a Poisson distribution.
The tail of the distribution is quickly disappearing.
This deterministic limit corresponds
to a limiting time-change of equation $\eqref{Sd}$ 
of the form $t' = t/\eps$
with $\alpha = \alpha'/\eps$, $\lambda = \lambda'/\eps$
as $\eps$ tends to 0.
The Poisson distribution has a mean of $\lambda'/\alpha' = \lambda/ \alpha$
so that it may be possible to quantify 
much more precisely than we do
the threshold in the number of deleterious mutations
after which differentiating individuals
is not so crucial.
This could make it possible 
to obtain quantitative bounds from our arguments
in the context of very large populations
(in the vicinity of the deterministic limit).


\section*{Proofs}

All the following properties strongly depend on the values of $\alpha>0$
($\alpha \in (0, 1)$ for the discrete state-space)
and $\lambda>0$. 
For brevity since the line of proof holds for any such values, 
this expected dependency is not recalled in the following statements.
The dependencies in $N$ for the discrete state-space 
and in $d$ for the finite dimensional SDE 
are generally recalled 
(in the state spaces and the random times notably) 
because of the interest of having  statements that are uniform over these parameters.
They may still be omitted to avoid too heavy notations.


\section{Proof of Theorem \ref{M_ECVdisc}}
\label{M_sec_ecvdisc}

The proof of Theorem \ref{M_ECVdisc}
relies on the proof of Assumption $\mathbf{(A)}$
as stated in Subsection~\ref{M_genSt}.
For any $z\in \M_1^{(0), N}(\Z_+)$,
we denote by $\PR^N_z$ 
the law of the process $Z^N$ with initial condition $z$,
as defined in Subsection \ref{M_Haigh}.
This process is associated to the extinction time $\ext^N$
in \eqref{M_extDisc}.
The first step deals with the mixing estimate,
while we focus on the persistence of large mutational burdens
in the second step 
(for the estimate on the escape from the transitory domain).


\paragraph{Step 1: access to any focal state}
 \hfill\\
We first prove  that 
any focal state of the population 
can be reached uniformly in the initial condition 
with a non-neglible probability, as stated in the upcoming Proposition~\ref{M_PosTrans}:
\begin{prop}
	\label{M_PosTrans}
	For any integer $N\ge 1$ 
	 and any element $z \in \M_1^{(0), N}(\Z_+)$:
	\begin{equation*}
		\inf \Lbr \PR^N_{z_0}(Z^N(1) = z)\bv z_0\in \M_1^{(0), N}(\Z_+)\Rbr >0.
	\end{equation*}
\end{prop}

\bpf:
We impose that all the individuals
of the next generation are the offspring
of an individual without any mutation,
and prescribe the number of mutations that they get
from the profile of $z$.
For any initial condition $z_0 \in \M_1^{(0), N}(\Z_+)$,
the probability
of choosing a fittest individual as a parent is:
$z_0(0)/(\Tsum{i\ge 0} z_0(i) \cdot (1-\alpha)^i).$
This probability is uniformly  lower-bounded
by $1/N$.
The number of mutations is then chosen independently of $z_0$,
and there is indeed a positive probability 
that the sequence of independent Poisson distributed random variables
has an empirical law distributed as $z$.
This concludes the proof of Proposition \ref{M_PosTrans}.
\epf
%

\paragraph{Step 2: disappearance of any large mutational burden}\hfill\\
With a probability close to 1,
the sub-population of individuals 
carrying a large number of mutations 
leave no progeny, as stated in the upcoming Proposition~\ref{M_Pers}:
\begin{prop}
	\label{M_Pers}
	For any $N\ge 1$ and $\eps>0$,
	there exists $K\ge 1$ such that
	the following inequality 
	holds with $E^N:= \big\{z \in \M_1^{(0), N}(\Z_+)\pv z([\![K, \infty[\![)=0\big\}$
	for any $z \in \M_1^{(0), N}(\Z_+)$:
	\begin{equation*}
		\PR^N_z(Z^N(1) \notin E^N) \le \eps.
	\end{equation*}
\end{prop}

\bpf:
Let $k_1\ge 1$ for the threshold in the number of mutations.
The probability that an individual chooses a parent with more than ${k_1}$ mutations
is upper-bounded by $N\cdot(1-\alpha)^{k_1}$, since $z(0) \ge 1$.
For any $\eps>0$, we thus choose ${k_1}\ge1$
such that,
with a probability greater than $1-\eps/2$, 
no individual in the next generation descends 
from an individual with more than ${k_1}$ mutations.
Likewise, there exists $k_2\ge 1$
such that,
with a probability greater than $1-\eps/2$, 
the number of additional mutations is less than $k_2$
(for any individual, independently of the initial condition $z$).
We can then conclude the proof of Proposition \ref{M_Pers} in that:
\begin{equation*}
\frlq{z \in \M_1^{(0), N}(\Z_+)}
\PR^N_z(Z^N(1) \notin E^N) \le \eps,
\end{equation*}
where $E^N$ is defined as in the Proposition \ref{M_Pers}
with $K= k_1 + k_2$.
\epf

\paragraph{Concluding the proof of Theorem \ref{M_ECVdisc}}
\hfill\\
We simply set $\cD^N_\ell$ to be the whole space $\M_1^{(0), N}(\Z_+)$ for any $\ell$.
Note that $(A0)$ is satisfied even for this degenerate case.
Actually, the exit time are just infinite and the entry times in $\cD^N_\ell$ always equal zero.
Secondly, Proposition \ref{M_PosTrans}
 implies Assumption $(A1)$.
Concerning $(A2)$,
we exploit Proposition \ref{M_Pers} 
 inductively over $k\ge 1$
 to deduce an upper-bound of the following form,
 thanks to the Markov property
 and with $\tau_E^N$ the hitting time of $E^N$
 by the process $Z^N$:
\[\frlq{z \in \M_1^{(0), N}(\Z_+)}
\PR^N_z(k < \ext^N\wedge \tau^N_E)
\le 2^{-k}\cdot\exp(-\rho k).\]
It implies the upper-bound on the exponential moment for any $z$
by splitting the expectation 
depending on the interval of the form $[k, k+1)$ 
that contains $\ext^N\wedge \tau^N_{E}$, i.e.:
\begin{equation*}
\begin{split}
	\bE^N_z(\exp[\rho\cdot (\ext^N\wedge \tau^N_{E})])
&\le \sum_{k\ge 0} \exp[\rho\cdot(k+1)]\cdot \PR^N_z\big(\ext^N\wedge \tau^N_{E}\in [k, k+1)\big)
\\&\le e^\rho\; {\textstyle \sum_{k\ge 0} 2^{-k}}
= 2 e^\rho <\infty.
\end{split}
\end{equation*}

For the last criterion $(A3)$,
we remark that $E^N$ as defined in Proposition \ref{M_Pers}  
is finite.
Thanks to Proposition \ref{M_PosTrans} and to the Markov property, 
we can thus choose
$c= c(N)>0$ such that the following comparisons of survival hold for any $t\ge 1$:
\begin{equation}
\PR^N_{\delta_0}(t< \ext^N)
\ge c\, \sup_{z\in E^N} \PR^N_{z}(t -1 < \ext^N)
\ge c\, \sup_{z\in E^N} \PR^N_{z}(t < \ext^N).
\label{A3_disc}
\end{equation}
This concludes $(A3)$ and that Assumption $\mathbf{(A)}$ is satisfied.

Thanks to Theorem \ref{M_th:ECV}, 
the semigroup $P^N$ 
displays QSC
with characteristics
$(\nu^N, h^N, \rho_0^N)$ $\in \M_1(\M_1^{(0), N}(\Z_+))\ltm B_+(\M_1^{(0), N}(\Z_+))\ltm \bR_+$.
Moreover, $h^N$ is uniformly bounded away from 0.
Besides, $\rho_0^N = -\log\big[\PR^N_{\nu^N}(1 < \ext)\big]>0$
	 because $\PR^N_z(\ext = 1)>0$ holds for any $z\in \M_1^{(0), N}(\Z_+)$.
This concludes the proof of Theorem \ref{M_ECVdisc}.
\epf

\begin{rem}
	If we were to replace the law of $\xi$ by a Bernoulli distribution (mutations occurring one by one),
	Proposition \ref{M_PosTrans} would still hold with the restriction of $z = \delta_0$,
	which is the only case we need.
	It would extend to any $z$ provided we change the time $1$ 
	by the maximal number of mutations in $z$.
	The proof would not be much more difficult with overlapping generations,
	except that individuals should then be removed one by one.
	The proof of the equivalent of Proposition \ref{M_Pers}
	would merely be slightly more difficult.
	The details are left to the interested reader.
\end{rem}

\section{Proof of Theorem \ref{M_ECVFin}}
\label{M_sec_ECVFin}

The proof of Theorem \ref{M_ECVFin}
also relies on the set $(\mathbf{A})$ of criteria,
as stated in Subsection~\ref{M_genSt}.
The proofs of two of these criteria 
(especially the mixing estimate
and the asymptotic comparison of survival)
are very much inspired by those of \cite[Subsection 4.2.2]{AV_QSD}.
Likewise, they exploit Harnack's inequality
--the following Property $(H)$--
classically deduced for elliptic diffusions,
see Subsection~\ref{sec_Harnack}.
The estimate of escape
on the other exploits several comparisons with one-dimensional diffusions
to deal with the behavior of the process near the boundary of the domain,
by exploiting the classical results presented in Subsection~\ref{sec_bound}.

The mixing estimate is then proved as the first step 
(in Subsection \ref{sec_Mix_fin})
and the asymptotic comparison of survival as the second step
(in Subsection \ref{sec_ACS_fin}). 
We turn next to the estimate of the escape from the transitory domain
(in Subsection \ref{M_sec_etD}),
then to estimations of lower-bound on the survival capacity
(in Subsection \ref{M_SHPos})
before we can conclude the proof of Theorem~\ref{M_ECVFin}
in Subsection~\ref{M_sec_ECvfin}.

We consider the following increasing closed subsets of $\cX_d$ as a reference:
\begin{equation}
\cD\eD_\ell:= \Lbr x = (x_i)_{i\in \II{0, d}}\in 
\lc \frac{1}{2 \ell d}\mVg 1 - \frac{1}{2 \ell d}\rc ^d 
\pv \sum_{i= 0}^d
x_i = 1 \Rbr.
\label{M_Dn}
\end{equation} 

Remark that $T\eD_{\cD_{\ell}}\wedge \ext\eD = T\eD_{\cD_{\ell}}$
holds for any $\ell$, where $T\eD_{\cD_{\ell}}$
denotes  the exit time of the process $X\eD$ out of $\cD\eD_{\ell}$.
Similarly and for simplicity, the superscript $(d)$ for the set 
$\cD\eD_{\ell}$ is dropped in the corresponding entry time, namely $\tau\eD_{\cD_{\ell}}$.

For any $d\in \II{1, \infty}$ and any $x\in \cX_d$,
the process $X\eD$ is solution under $\PR_x$
of the system~\eqref{Sd} 
with initial condition $X\eD_i(0) = x_i$.

There is no real ambiguity in the dependency in $d$,
especially since this dependency shall be made visible 
for the process $X\eD$ and the associated processes and stopping times,
By extension, for any probability measure $\zeta$ on $\cX_d$,
$\PR_\zeta(\RMd \omega) = \int_{\cX_d}  \PR_x(\RMd \omega) \zeta(dx)$.
Note that $\PR$ can be seen as specifying the law of the sequence $(W_i)_{i\in \Z_+}$
while $x\in \cX_d$ can be seen as an element of $\cX_\infty$ (with $x_i = 0$ for any $i\ge d+1$),
so that a common notation makes sense as well. 
\\

\subsection{First crucial properties}

\subsubsection{Harnack's inequality}
\label{sec_Harnack}

Property $(H)$ is defined as follows for any
process $(Y(t))_{t\ge 0}$ on $\Y$
with generator $\cL$
(including possibly an extinction rate $\rho_c$):\\
\textsl{
	For any two connected open relatively compact sets $\mathfrak{K}^\wedge, \mathfrak{K}^\vee \subset \Y$ 
	with $\cC^\infty$ boundaries such that $\overline{\mathfrak{K}^\wedge}\subset \mathfrak{K}^\vee$
	(where $\overline{\mathfrak{K}}$ denotes the closure of the set $\mathfrak{K}$),
and any	$0<t_1<t_2$,
there exists $C = C(\cL, t_1, t_2, \mathfrak{K}^\wedge, \mathfrak{K}^\vee)>0$
such that the following properties hold
for any positive $\cC^\infty$ constraints:
	$u_{\partial \mathfrak{K}^\vee}: (\{0\} \times \mathfrak{K}^\vee) \cup ([0, t_2] \times \partial \mathfrak{K}^\vee)
	\rightarrow \bR_+$.
There exists a unique positive strong solution $u(t,x)$ to the following Cauchy problem:}
\begin{equation*}
	\begin{aligned}
&\partial_t u (t, x) = \cL u (t, x) 
\\&u(t, x) = u_{\partial \mathfrak{K}^\vee}(x) 
	\end{aligned}
\qquad 
	\begin{aligned}
	&\text{ on } [0, t_2] \times \mathfrak{K}^\vee,
	\\&\text{ on } (\{0\}\times \mathfrak{K}^\vee) \cup ([0, t_2] \times \partial \mathfrak{K}^\vee),
\end{aligned}
\end{equation*}
\textsl{and $u$ satisfies the following inequality:}
\begin{equation*}
	{\textstyle \inf_{x\in \mathfrak{K}^\wedge} u(t_2, x) 
		\ge C\, \sup_{x\in \mathfrak{K}^\wedge} u(t_1, x).}
\end{equation*}

For any $d\in \N$ and $\ell\in \N$,
thanks to Proposition \ref{prop_Trep}
in the appendix, we identify the generator $\cL^{(d)}$ 
of the finite dimensional process $X^{(d)}$
in restriction to the set $\cD\eD_\ell$ 
to the following nondivergence form, for $u\in C^2(\cD\eD_\ell)$ and $x \in \cD\eD_\ell$:
\begin{equation*}
\cL^{(d)} u(x)
= \frac{1}{2}\sum_{i, j = 1}^d \sigma^{(d)}_{i, j}(x) \partial^2_{i, j} u(x) 
+ \sum_{i= 1}^d b^{(d)}_{i}(x) \partial_{i} u(x).
\end{equation*}
\begin{lem}
	\label{lem_Har}
Property $(H)$ 
holds for any $d\in \N$ and any $\ell\in \N$
for  the finite dimensional process $X^{(d)}$
with generator $\cL^{(d)}$ 
in restriction to the set $\cD\eD_\ell$
as defined in \eqref{M_Dn}.
\end{lem}
\begin{proof}
For any $d\in \N$ and any $\ell\in \N$,
the diffusion matrix $\sigma^{(d)}$ is uniformly elliptic on $\cD\eD_\ell$
while $\sigma^{(d)}$ and the drift term $b^{(d)}$ are $\cC^\infty$ on $\cD\eD_\ell$.
As noted in \cite[Section 4.2.2]{AV_QSD},
this entails Property $(H)$.
We sketch the argument for completeness,
rather referring to \cite{Ev10} for the clarity
of its presentation.
The existence and uniqueness of the solution $u$,
with the fact that $u\in \cC^\infty$,
is a consequence of \cite[Theorem 7]{Ev10}.
Thanks to  \cite[Theorem 8]{Ev10},
this solution is positive.
We can then apply \cite[Theorem 10]{Ev10} on $u$
to deduce Harnack's inequality
and complete the proof of Property $(H)$.
\end{proof}

\subsubsection{Boundary classification for one-dimensional diffusions}
\label{sec_bound}
The following proofs rely on comparison principles
with one-dimensional diffusions, 
for which the boundary classification is well-described.
 We first present briefly in the upcoming Lemma~\ref{bound_oneD} the conclusions from 
 \cite[Section~6, Chapter~15]{KT81}
for the specific cases of solutions $Y$ 
on the state-space $(0, 1)$
to SDEs of the following form:
\begin{equation}
\RMd Y(t) = b\cdot Y(t) \RMd t + \sqrt{Y(t)\cdot (1-Y(t))}\RMd B(t)
\label{eq_SDE_bound}
\end{equation}
where $B$ is a standard Brownian motion
and $b\in \bR$.

To motivate this form of the diffusion coefficient,
note that the martingale term of each coordinate $X_i$ taken separately
takes actually this form, as one can check from Proposition~\ref{prop_Trep}
in the appendix.
The vicinity of $0$ will inform on the behavior of the boundary
for a linear drift
and the vicinity of $1$ for a nearly constant drift.

\begin{lem}
	\label{bound_oneD}
$(i)$	For any $b\in \bR$,
	the left-boundary $0$
is accessible to the solution $Y$ of \eqref{eq_SDE_bound},
which entails the following convergence for any $t>0$:
\begin{equation*}
	\lim_{y \rightarrow 0} \PR_y(\tau^Y_0 < t) = 1.
\end{equation*}
$0$ is actually an exit boundary,
so that $Y$ gets absorbed at $0$.

$(ii)$ 	For any $b>-\tdiv{2}$,
the right-boundary $1$
is accessible,
which entails the following convergence for any $t>0$:
\begin{equation*}
\lim_{y \rightarrow 1} \PR_y(\tau^Y_1 < t) = 1.
\end{equation*}

1 is actually an exit boundary iff $b\ge 0$.
If $b\in (-\tdiv{2}, 0)$, 
1 is a regular reflecting boundary,
which implies that $Y(t)<1$ holds a.s. for any $t>0$.

$(iii)$ 	For any $b\le -\tdiv{2}$,
the right-boundary $1$
is inaccessible
and an entrance boundary, in the following sense:
for any $y\in (0, 1)$,
$$\lim_{t\ifty} \Tinf{z\in (y, 1)} \PR_z(\tau^Y_y < t) = 1.$$
\end{lem}
\begin{rem}
In the case of one-dimensional diffusions,
there exists a classical reformulation of the process,
which involves the scale function and the speed measure,
that especially helps for the analysis of the boundary.
This reformulation has also been exploited for the study of quasi-stationarity 
of one-dimensional diffusions,
be it with spectral techniques that are specific and more classical for such diffusions
(as in \cite{CC+09})
or recent extensions more closely related to our approach
(notably in \cite{CV18}).
Since we do not exploit further
the specific properties of an exit boundary and a regular reflecting boundary,
these definitions are not specified and
the interested reader is deferred to \cite[Section~6, Chapter~15]{KT81}.
\end{rem}

As a particular case of \cite[Proposition 3.12]{PR14},
we have the upcoming Lemma~\ref{lem_comp}, 
that will often be exploited for our comparison estimates
on the boundaries:
\begin{lem}
	\label{lem_comp}
Let $\hat b, \check b: \Omega\times \bR_+\times [0, 1]\mapsto \bR$
and $L>0$
be such that $\hat b(., ., y), \check b(., ., y)$
are measurable for any $y\in [0, 1]$
and that $\hat b(\omega, t, .), \check b(\omega, t, .)$ 
are L-Lispchitz continuous for all $\omega, t\in  \Omega\times \bR_+$.
Assume that $\hat b(\omega, t, y)\ge \check b(\omega, t, y)$ 
holds for any $(\omega, t, y)\in \Omega\times \bR_+\times [0, 1]$
and that $\hat y, \check y\in [0, 1]$ are such that $\hat y\ge \check y$.
Then the inequality $\hat Y(t) \ge \check Y(t)$ holds for any $t\in \bR_+$
for the solution $\hat Y$ and $\check Y$ to the following SDEs:
\begin{equation*}
\begin{split}
	\RMd \hat Y(t) 
&:= \hat b_t(\hat Y(t))\RMd t + \sqrt{\hat Y(t)\cdot (1-\hat Y(t))}\RMd B(t)
\quad \hat Y(0) = \hat y,
\\	\RMd \check Y(t) 
&:= \check b_t(\check Y(t))\RMd t + \sqrt{\check Y(t)\cdot (1-\check Y(t))}\RMd B(t)
\quad \check Y(0) = \check y,
\end{split}
\end{equation*}
where $B$ is a standard Brownian motion.
\end{lem}
In practice, we will exploit the following corollary
instead of Lemma~\ref{bound_oneD}$(ii)$ several times 
(the corollary being implied by this lemma together with Lemma~\ref{lem_comp}).
\begin{cor}
	\label{cor_bound_oneD}
Consider for any $\varphi, \psi\in \bR$ the solution $Z$
to the following SDE:
\begin{equation*}
	\RMd Z(t) = \big[\varphi + \psi\cdot Z(t)\big] \RMd t + \sqrt{Z(t)\cdot (1-Z(t))}\RMd B(t),
\end{equation*}
where $B$ is a standard Brownian motion.
For any $\varphi\in (0, 1/2)$ and $\psi\in \bR$,
0 is a regular reflecting boundary for $Z$.
\end{cor}

\subsection{Step 1: mixing estimate}
\label{sec_Mix_fin}
The aim of this subsection is exactly $(A1)$, 
as stated in the upcoming Proposition~\ref{M_Mixd}:
\begin{prop}
	\label{M_Mixd}
For any $d\in \N$, there exists $\alc\eD\in \M_1(\cX_{d})$
	with support in $\cD\eD_2$
	such that the following property holds for any $\ell \ge 1$.
	There exists $c>0$ such that:
	\begin{equation*}
		\frlq{x\in \cD\eD_\ell} 
		\PR_x\lp X\eD(1) \in \RMd y \pv 1 < T\eD_{\cD_{\ell +1} }\rp
		\ge c\; \alc\eD(\RMd y).
	\end{equation*}
\end{prop}
For the comparison with $(A1)$,
recall that 
$T\eD_{\cD_{\ell}}= \ext\eD \wedge T\eD_{\cD_{\ell}}$ 
holds for any $\ell$.

\bpf: Let $\ell\ge 1$.
We choose two connected bounded open sets $\mathfrak{K}^\wedge, \mathfrak{K}^\vee\subset \cD\eD_{\ell+1}$ 
with $\cC^\infty$ boundaries such that 
$\cD\eD_\ell \subset \mathfrak{K}^\wedge$ and
$\overline{\mathfrak{K}^\wedge}\subset \mathfrak{K}^\vee$.
Thanks to Lemma~\ref{lem_Har},
we apply Property $(H)$ 
 (see Subsection~\ref{sec_Harnack})
to $u(t,x):= \bE_x\lp f(X\eD(t))\pv t< T\eD_{\mathfrak{K}^\vee} \rp$,
where $f$ is any non-negative $\cC^\infty$ function with support in $\cD\eD_2 \subset \mathfrak{K}^\wedge$,
and 
$ T\eD_{\mathfrak{K}^\vee} 
:= \inf\{t\ge 0 \pv
X^{(d)}(t) \notin \mathfrak{K}^\vee \}$.
Since in addition $T\eD_{\cD_{2}} \le T\eD_{\mathfrak{K}^\vee} \le T\eD_{\cD_{\ell +1}}$,
there exists a real $C_\ell\eD >0$ 
such that the following inequality holds 
for any $x\in \cD\eD_\ell$,
 $z\in \cD\eD_1$:
\begin{equation*}
\bE_x\lp f(X\eD(1))\pv 1< T\eD_{\cD_{\ell+1}}\rp
\ge C_\ell\eD\; \bE_{z} \lp f(X\eD(\tdiv{2}))\pv \tdiv{2}< T\eD_{\cD_{2}}\rp.
\end{equation*}
With the arbitrary choices of $z\eD$ as the barycenter of $\cD\eD_1$,
we then define the probability measure $\alc\eD$ as follows:
\begin{equation*}
\alc\eD(\RMd x) := \PR_{z\eD} \lp X\eD(\tdiv{2}) \in \RMd x\bv \tdiv{2}< T\eD_{\cD_{2}}\rp,
\end{equation*}
which  has support on $\cD\eD_2$ and is independent of $\ell$.
Since the constant $C_\ell\eD$ does not depend on $f$,
we deduce with $c_\ell = C_\ell\eD\cdot \PR_{z\eD} \lp \tdiv{2}< T\eD_{\cD_{2}}\rp>0$
the following inequality between measures 
for any initial condition $x\in \cD\eD_\ell$:
\begin{equation*}
	\bE_x\lp X\eD(1) \in \RMd y \pv 1 < T\eD_{\cD_{\ell+1}}\rp
	\ge c_\ell\; \alc\eD(\RMd y),
\end{equation*}
which concludes the proof of Proposition~\ref{M_Mixd}.
\epf

\subsection{Step 2: asymptotic comparison of survival}
\label{sec_ACS_fin}
The aim of this subsection is to prove the upcoming Proposition~\ref{M_Cpld},
thanks to which we will deduce $(A3)$:
\begin{prop}
	\label{M_Cpld}
	The following boundedness property holds for any $d\in \N$
	and any $\ell\ge 1$:
	\begin{equation*}
		\limsup_{\tp \ifty}\sup_{x, x'\in \cD\eD_\ell}
		\dfrac{\PR_x\lp \tp < \ext\eD\rp}
		{\PR_{x'}\lp \tp < \ext\eD\rp} < \infty.
	\end{equation*}
\end{prop}

The proof of Proposition \ref{M_Cpld}
 is  similar to the one of Proposition \ref{M_Mixd}.
 It is nonetheless
more technical because  
 we can no longer neglect trajectories 
exiting $\cD\eD_{\ell+1}$.\\

\bpf:
We can  find two connected open relatively compact sets $\mathfrak{K}^\wedge, \mathfrak{K}^\vee\subset \cX_d$ 
with $\cC^\infty$ boundaries such that $\cD\eD_\ell \subset \mathfrak{K}^\wedge$
and $\overline{\mathfrak{K}^\wedge}\subset \mathfrak{K}^\vee\subset int(\cD\eD_{\ell+1})$.
We want to approximate the function:
 $$u(t, x):= \bE_{x}\lp  f(X\eD(t)) \pv t < \ext\eD \rp,
  \with \; t  \ge 1,\; x \in \mathfrak{K}^\vee$$
 defined for some non-negative $f \in \cC^\infty(\cX_d)$.
 Thanks to  \cite[Theorem 5.1.15]{L95},
 $u$ is continuous. 
 It is clearly non-negative.
 However, it is a priori not regular enough 
 to apply Harnack's inequality directly.
 Thus, we approximate it
 on the parabolic boundary \mbox{$[1,\, \infty) \times \partial \mathfrak{K}^\vee\,$}
  \mbox{$\bigcup\, \{1\} \times \mathfrak{K}^\vee$}
  by some family $(U_k)_{k\ge 1}$ of non-negative smooth  functions.
 We then deduce approximations of $u$ in $[0,\, \infty) \times \mathfrak{K}^\vee$ 
 by solutions $u_k$ to
 the following Cauchy Problem:
 \begin{equation*}
 \begin{aligned}
 &\partial_t u_k (t, x) - \cL u_k (t, x)= 0,
 \\& u_k(t, x) = U_k(t+1, x),
 \end{aligned}
\qquad 
 \begin{aligned}
	&\text{ for } t \ge 0,\; x \in int(\mathfrak{K}^\vee)
	\\& \text{ for } t \ge 0,\; x \in \partial \mathfrak{K}^\vee, \quad 
	\text{ or } t = 0,\, x \in \mathfrak{K}^\vee.
\end{aligned}
 \end{equation*}
  
 Thanks to Property $(H)$
 (in Subsection~\ref{sec_Harnack}), 
 the constant involved in Harnack's inequality 
 does not depend on the values of $U_k$ on the boundary.
 Thus, it applies with the same constant 
 for the whole family of approximations $u_k$.
 With $t_1:= 1$ and $t_2:= 2$, 
we can thus choose $C_\ell\eD >0$ such that for any $k$ and any $x, x'\in \cD\eD_\ell$:
 $$ u_k(1, x) \le C_\ell\eD u_k(2, x'),$$
 where the constant $C_\ell\eD$ does not depend on $f$ either.
Thanks to the proof in \cite[Section~4, step 4]{CV21},
the Harnack inequality extends to the approximated function $u$,
 with the convergence of $U_k$ on the parabolic boundary
 (and the time-shift of 1 taken into account).
It means that the following inequality
holds for any non-negative $f \in \cC^\infty(\Y)$
and any $x, x'\in \cD\eD_\ell$:
\begin{equation*}
\bE_{x}\lp  f(X\eD(2)) \pv 2< \ext\eD \rp
\le C_\ell\eD\, \bE_{x'}\lp   f(X\eD(3))
\pv 3 < \ext\eD \rp
\end{equation*}
The inequality then extends to any measurable and bounded $f$.
We now fix $t\ge 2$ and apply this result to 
the function
$f_t(x):= \PR_x(t-2 < \ext\eD)$,
so that,
thanks to the Markov property,
 the following comparison of survival holds for any 
$x, x'\in \cD\eD_\ell$
and any $t\ge 2$:
\begin{equation*}
\begin{split}
	\PR_{x}\lp t < \ext\eD \rp
&\le C_\ell\eD\, \PR_{x'}\lp t + 1 < \ext\eD \rp
\\&\le C_\ell\eD\, \PR_{x'}\lp t< \ext\eD \rp.
\end{split}
\end{equation*}
Note that, since $C_\ell\eD$ does not depend upon $f$, it does not depend upon $t$.
This   concludes the proof of Proposition \ref{M_Cpld}.
 \epf

\subsection{Step 3: escape from the transitory domain}
\label{M_sec_etD}

The aim of this subsection is to prove the upcoming Proposition~\ref{M_Tdy},
that entails $(A2)$:
\begin{prop}
	\label{M_Tdy}
	For any $d\in \N$ and $\rho>0$, 
	there exists $\ell\ge 1$ such that:
	\begin{equation*}
		\sup_{x\in \cX_d} \bE_x \exp\Big[\rho\cdot \Big(\tau\eD_{\cD_\ell}\wedge \ext\eD\Big)\Big] \le 16.
	\end{equation*}
\end{prop}

\paragraph{Two elementary steps}
\label{M_sec_Lxke}
\hfill\\
The proof 
is achieved with two forthcoming  lemmas
as intermediate steps.
We first prove that the click is very likely to happen
when the size of the optimal subpopulation is small,
as stated in the upcoming Lemma~\ref{M_Lx0e}:
\begin{lem}
	\label{M_Lx0e}
	For any $d\in \N$ and any time $t>0$,
	the following supremum tends to 0 as $y_0$ tends to $0$:
	\begin{equation*}
		\sup\Lbr \PR_x\lp t<\ext\eD\rp
		\Big \vert x\in \cX_d,
		x_0 \le y_0 \Rbr.
	\end{equation*}
\end{lem}
\bpf : 
Provided the initial condition $x$ is such that $x_0\le y_0$,
the process $X\eD_0$,
namely the initial component of the solution to \eqref{Sd}-\eqref{eq_def_Bi},
 is upper-bounded
by the solution $Y$ to the following SDE,
thanks to Lemma \ref{lem_comp}:
\begin{equation*}
	\RMd Y(t) = (\alpha d) \cdot Y(t)\, \RMd t + \sqrt{Y(t)\cdot (1-Y(t))}\, \RMd B_0(t)
	\mVg \quad Y(0) = y_0.
\end{equation*}
Thanks to Lemma \ref{bound_oneD}(i),
 $0$ is an exit boundary of $Y$,
 which concludes the proof of 
Lemma~\ref{M_Lx0e}.
\epf\\

We then deal iteratively with each subclass size
to prove that these sizes escape the vicinity of 0,
as stated in the upcoming Lemma~\ref{M_Lxke}.
\begin{lem}
	\label{M_Lxke}
	For any integer $J\in \II{1, d}$, any $y\in (0,1)$,
	and any time $t>0$,
	the next infimum tends to 1 as the constant $y'\in (0, y)$ tends to $0$:
	\begin{equation*}
		\inf\Lbr \PR_x\lp \tau\ekD J_{y'}< t \wedge \ext\eD\rp
		\Big \vert x\in \cX_d,
		\frl{j\le J-1} x_j \ge y \Rbr,
	\end{equation*}
where
$\tau\ekD J_{y'}:= \inf\{s\ge 0\pv \frl{j\le J} X\eD_j(s) \ge y'\}$.
\end{lem}

\bpf : 
Let $J\in \II{1, d}, y\in (0, 1), \eps,  t>0$.
The process $X\eD_j$, 
namely the $j$-th component of the solution to \eqref{Sd}-\eqref{eq_def_Bi}
for any $j\in \II{0, J-1}$,
 is lower-bounded
 by solutions $Y_j$
to SDEs of the following form,
thanks to Lemma \ref{lem_comp}:
\begin{equation*}
	\RMd Y(t) = - (\lambda +\alpha J) \, \RMd t + \sqrt{Y(t)\cdot (1-Y(t))}\, \RMd B(t)
	\mVg \quad Y(0) = y,
\end{equation*}
where $B$ is a standard Brownian motion.
We then choose $t'\in (0, t)$ such that $Y$ 
stays above $y/2$ on the time-interval $[0, t']$
with probability greater than $1-\eps/(2 J)$.
We then exploit this property on $Y_j$ for any $j\in  \II{0, J-1}$,
so as to deduce the following inequality 
for any $d\in \N$ and any $ x\in \cX_d$ such that $x_j \ge y$
holds for any $j\le J-1$:
\begin{equation}
\PR_x\Big( t < T\ekD{J-1}_{y/2}\Big)
\ge 1- \eps/2,
\label{eq_ydiv2}
\end{equation}
where $T\ekD{J-1}_{y/2}:= \inf\{s\ge 0\pv \Ex{j\le J-1} X\eD_j(s) \le y/2\}.$

Let $y_1:= \lambda\, y / (4\lambda  + 4\alpha J)$
and $\tau\ekD J_{y_1}$ defined as in Lemma~\ref{M_Lxke}.
The following inequalities thus hold for any  $s\le T\ekD{J-1}_{y/2}\wedge \tau\ekD J_{y_1}$:
\begin{equation*}
\begin{split}
	(\alpha\cdot M\eD_1(s) - \alpha J- \lambda) \cdot X\eD_J(s) + \lambda\cdot X\eD_{J-1}(s)
&\ge - (\alpha J + \lambda)\cdot y_1 +\lambda\, y / 2 
\\&\ge \lambda\, y / 4.
\end{split}
\end{equation*}
Thanks to Lemma \ref{lem_comp},
$X\eD_J$ is thus a.s. lower-bounded
on the time interval $[0, T\ekD{J-1}_{y/2}\wedge \tau\ekD J_{y_1}]$
 by the solution $Y_J$ to the following SDE:
\begin{equation*}
\RMd Y_J(s) =  \frac{(\lambda\, y)\wedge 1}{4}
	\, \RMd t + \sqrt{Y_J(t)\cdot(1-Y_J(t))}\, \RMd B_J(t)
\mVg \quad Y_J(0) = 0,
\end{equation*}
where $B_J$ is a standard Brownian Motion.
Thanks to Corollary~\ref{cor_bound_oneD},
the left-boundary $0$ is regular reflecting for $Y_J$.
Therefore, there exists $0<y'\le y_1\wedge (y/2)$ such that:
\begin{equation}
\PR(\textstyle{\sup_{\{s\le t'\}}}  Y_J(s)< y') \le \eps / 2.
\label{eq_yP}
\end{equation}
Provided that the initial condition satisfies 
that $x_j\ge y$ for any $j\in \II{0, J-1}$,
the property $\tau\ekD J_{y'} < t\wedge \ext\eD$ holds a.s.
on the event $\Lbr \sup_{s\le t'} Y_J(s)\ge  y'\Rbr 
\cap \Lbr t' < T\ekD{J-1}_{y/2} \Rbr$, 
which occurs with probability greater than $1-\eps$
thanks to \eqref{eq_ydiv2} and \eqref{eq_yP}.
This ends the proof of Lemma~\ref{M_Lxke}.
$\hfill \square$
\\

\paragraph{Concluding the proof of Proposition \ref{M_Tdy}}
\hfill\\
Given $\rho >0$,
let $t_0:= \log(2)/\rho$. 
We can choose $y_0\in (0,1)$
thanks to Lemma \ref{M_Lx0e}  such that 
the probability of survival up to time $t_0$ is small as follows
for any initial condition $x = (x_i)$ such that $x_0< y_0$:
\begin{equation}
\PR_x\lp t_0<\ext\eD\rp 
\le \exp(-\rho\,  t_0) / 2 = \tdiv{4}.
\label{eq_def_y0}
\end{equation}

Thanks to Lemma \ref{M_Lxke},
we iteratively obtain lower-bounds for the different components of $X\eD$,
so that for any $1\le J\le d$
and any $y_{J-1}\in (0, 1)$,
there exists $y_{J} \in (0, y_{J-1})$
such that the following inequality holds:
	\begin{equation*}
	\inf\Lbr \PR_x\lp \tau\ekD J_{y_J}< (t_0/d) \wedge \ext\eD\rp
	\Big \vert x\in \cX_d,
	\frl{j\le J-1} x_j \ge y_{J-1} \Rbr
	\ge 1 - 1/(4 d).
\end{equation*}
This property defines iteratively the sequence 
$(y_J)_{J\in \II{1, d}}$ in terms of $y_0$.

Thanks to the strong Markov property
at times $\tau\ekD J_{y_J}$
and by induction
on $0\le J\le d$, 
we deduce that the following inequality
holds for any $x\in \cX_d$ such that $x_0\ge y_0$:
%
\begin{equation*}
\PR_x\lp \tau\ekD J_{y_J} \le (J\cdot t_0/d)\wedge \ext\eD \rp \ge 1 - J/(4 d).
\end{equation*}
Let $E\eD = \cD\eD_\ell$ for some $\ell \ge y_d/(2 d)$,
so that $E\eD\in \mathbf D\eD$ and $\tau\eD_E\le \tau\ekD d_{y_d}$.
Thanks to the previous inequality, the following one 
holds for any $x\in \cX_d$ such that $x_0\ge y_0$:
\begin{equation*}
\begin{split}
		\PR_x\lp t_0 < \tau\eD_E\wedge \ext\eD\rp
	&\le 1 - 
	\PR_x\lp\tau\ekD d_{y_d} \le t_0\wedge \ext\eD \rp 
	\\&\le \tdiv{4} = \exp(-\rho\,  t_0) / 2,
\end{split}
\end{equation*}
which extends the inequality stated in \eqref{eq_def_y0}
for any $x$ such that $x_0<y_0$.

By induction over $k\ge 1$	thanks to the Markov property at times $k\, t_0$, 
 the following inequality holds for any $x\in \cX_d$ and any $k\ge 1$:
	\begin{equation*}
		\PR_x (k\, t_0\le \ext\eD\wedge \tau\eD_{E})\le  2^{-k}\exp(-\rho\, k\, t_0).
	\end{equation*}
It entails the following upper-bound on the exponential moment:
	\begin{align*}
		\sup_{\{x\in \cX_d\}} 
		\bE_x \lc \exp\lp \rho\cdot (\ext\eD\wedge \tau\eD_{E})\rp\rc
		&\le \sum_{k\ge 0} \exp(\rho\cdot (k+1)\cdot t_0 ) \,
		\sup_{\{x\in \cX_d\}} \PR_x (k\, t_0\le \ext\eD\wedge \tau\eD_{E})
		\\&\textstyle
		\le  e^{\rho t_0}\sum_{k\ge 0} 2^{-k} = 4 < \infty.
	\end{align*}
 This concludes the proof of Proposition \ref{M_Tdy}.
\epf

%
\subsection{Step 4: Lower-bound of the survival capacity}
\label{M_SHPos}

The aim of this subsection is to prove the upcoming Lemma~\ref{M_HPos},
which implies 
some lower-bound of the survival capacity:
\begin{lem}
	\label{M_HPos}
	For any $y_0>0$, 
		there exists $t>0$ and $\ell\ge 1$ such that 
		the following holds:
		\begin{equation*}
\inf\Big\{\PR_x\big(\tau\eD_{\cD_\ell} \le t \wedge \ext\big)\bv
x\in \cX_d \pv x_0\ge y_0\Big\} >0.
		\end{equation*}	
\end{lem}
In other words, the sets
$H\eD_{y_0}:= \{x\in \cX_d \pv x_0\ge y_0\}$ 
 satisfy $(H_0)$
as stated in Proposition~\ref{th:UCV},
so that, for any $y_0>0$,
 $h\eD$ is uniformly bounded away from 0 on $H\eD_{y_0}$.
\\

\bpf:
Let any $y_0>0$. 
The process $X\eD_0$,
namely the initial component of the solution to \eqref{Sd}-\eqref{eq_def_Bi},
is  lower-bounded for any time
a.s. under $\PR_x$ 
by the solution $Y$ 
to the following SDE,
thanks to Lemma \ref{lem_comp}:
\begin{equation*}
\RMd Y(s) = -\lambda\, \RMd s
	+ \sqrt{Y(s)\cdot(1-Y(s))}\, \RMd B_0(s)\mVg \quad Y(0) = y_0.
\end{equation*}
We consider the extinction time for $Y$
as follows:
$\ext^Y:=\inf\{s \pv Y(s) = 0\}$.
We deduce that the probability of survival up to time $1$
is uniformly lower-bounded as follows
for any $d\ge 1$ and any $x \in H\eD_{y_0}$:
\begin{equation}\label{eq_min_X0}
\PR_{x}(1< \ext\eD) \ge c_S:= \PR_{y_0}(1< \ext^Y)>0.
\end{equation}
The choice of 1  for the time is arbitrary.
Thanks to Markov's inequality and to the exponential moment 
of $\tau\eD_{\cD_\ell}\wedge \ext\eD$
derived in Proposition~\ref{M_Tdy},
there exists $\ell$ such that for any $d\in \N$ and $x\in \cX_d$:
\begin{equation*}
	\PR_{x}(1< \tau\eD_{\cD_\ell}\wedge \ext\eD) \le c_S/2.
\end{equation*}
This  implies $(H_0)$ for $H\eD_{y_0}$ in the sense that
the following holds
 for any $x\in H\eD_{y_0}$:
\begin{equation*}
	\PR_x(\tau\eD_{\cD_\ell} \le 1\wedge \ext\eD)
	\ge \PR_{x}(1< \ext\eD) - \PR_{x}(1< \tau\eD_{\cD_\ell}\wedge \ext\eD)
	\ge c_S / 2.
	\SQ
\end{equation*}
\begin{rem}
Though this property $(H_0)$ holds uniformly in $d$,
it remains unchecked that the sequence $(h\eD)$ of functions
is uniformly bounded away from zero 
on the sequence $(H\eD_{y_0})$.
\end{rem}

\subsection{Concluding the proof of Theorem \ref{M_ECVFin}}
\label{M_sec_ECvfin}

For this proof,
we plan to exploit Theorem \ref{M_th:ECV}
and first ensure Assumption $\mathbf{(A)}$ 
(see Subsection \ref{M_genSt}). Let $d\in \N$.
The sets $\cD\eD_\ell$,
defined in  \eqref{M_Dn},
 satisfy $(A0)$.

Propositions \ref{M_Mixd}, \ref{M_Cpld} and \ref{M_Tdy}
ensure respectively $(A1)$, $(A3)$ and $(A2)$.
Notably, for $(A3)$,
since $\alc\eD$ has support in $\cD\eD_2$,
for any $\ell\ge 2$,
thanks to Proposition $\ref{M_Cpld}$:
\begin{equation*}
	\limsup_{\tp \ifty}\sup_{x\in \cD\eD_\ell}
	\dfrac{\PR_{x}\lp \tp < \ext\eD\rp}
	{\PR_{\alc\eD}\lp \tp < \ext\eD\rp} < \infty.
\end{equation*}
Thanks to Theorem \ref{M_th:ECV},
the semi-group therefore displays QSC
with some characteristics $(\nu\eD, h\eD, \rho\eD_0)$.
Thanks to Lemma \ref{M_HPos} and to Proposition \ref{th:UCV},
the survival capacity $h\eD$ is actually positive.
Any QSD $\nu'\in \M_1(\cX_d)$ must then satisfy 
both $\LAg \nu' \bv h\eD\RAg >0$.
Since $\nu'$ is a QSD,  $\LAg \nu' P\eD_t, \idg{}\RAg^{-1}\cdot
\nu' P\eD_t = \nu'$ holds for any $t$.
Thanks to  Corollary~\ref{D:CVAl},
this implies that $\nu'= \nu\eD$,
so that $\nu\eD$ is in fact the unique QSD.
With the upper-bound considered in Lemma \ref{M_Lx0e},
	we deduce that
$\PR_y(\ext\eD\le 1) >0$ holds for any $y\in \cX_d$.
Thus, 
$\rho\eD_0= -\log[\PR\eD_{\nu\eD}(1<\ext\eD)]>0$.
This concludes the proof of Theorem \ref{M_ECVFin}.\epf

\section{Proof of Proposition \ref{M_prop_MT}}
\label{M_sec_MT}

The proof of Proposition \ref{M_prop_MT},
concluded in Subsection~\ref{M_sec_pMT},
relies on two main steps,
handled uniformly over $d$.
The first step in Subsection~\ref{sec_augm} 
is to ensure that descent from large values of the moment
 quickly occurs with probability close to one;
the second step in Subsection~\ref{sec_Pm4}  prove that a too large increase of the moment 
is unlikely to occur.



\subsection{Step 1: descent of the moment}
\label{sec_Pm4}
The aim of this subsection is the upcoming Proposition~\ref{M_Pm4}. 
The descent of the $k$-th moment $M\eD_k = M_k(X\eD)$ of the process $X\eD$
 is stated in terms of the hitting time $\tau\eDk_{m}$:
\begin{equation}
\tau\eDk_{m}:= \inf\Big\{t\ge 0\pv M\eD_k(t)\le~m\Big\}.
\label{eq_def_tauMk}
\end{equation}
\begin{prop}
	\label{M_Pm4}
	For any time $t>0$ and
	any  $k \ge 1$,
	the following supremum tends to 0 as $m$ tends to infinity:
	\begin{equation*}
		\sup\Lbr \PR_x\lp t<\tau\eDk_{m} \wedge \ext\eD \rp
		\Big \vert d\in \N, x\in \cX_d\Rbr.
	\end{equation*}
\end{prop}

The proof of Proposition \ref{M_Pm4},
as achieved in Subsection~\ref{sec_fPm4},
relies on three steps detailed in the three forthcoming subsections.
The first step is focused on 
the vicinity of $\{X_0 = 0\}$,
the second on the vicinity of $\{M_1 = \infty\}$,
the last one being iterated for each moment between 2 and k.


\subsubsection
{Step 1.1: rare survival for small optimal subpopulation}
\label{sec_Lx0}
The aim of this subsection is the upcoming Lemma~\ref{M_Lx0},
which states, 
provided that the first moment is initially lower-bounded,
that a click is very close to occurring when the initial size $x_0$ 
of the optimal subpopulation is very small:
\begin{lem}
	\label{M_Lx0}
	For any time $t>0$,
	the following supremum tends to 0 as $\delta$ tends to 0:
	\begin{equation*}
		\sup\Lbr \PR_x\lp t<\ext\eD\rp
		\Big \vert \, d\in \N,\, x\in \cX_d,\,
		M_1(x) \ge 1,\,
		x_0\cdot M_1(x) \le \delta \Rbr.
	\end{equation*}
\end{lem}

\bpf: This proof is an extension of the one of  \cite[Proposition 3.8]{AP13}.

Let $t>0$ be fixed and define $\dlw$ as follows:
	\begin{equation*}
		\dlw:= \frac{1}{16\,\alpha}\wedge \frac{1}{4}.	
\end{equation*}
We exploit two parameters $m_1\ge 1$ and $\delta\in (0, \dlw)$:
the upper-bound shall hold for initial conditions $x$ 
such that both $m_1 \ge M_1(x)$
and $x_0\, m_1\le \delta$.
$m_1$ is freely chosen
 and $\delta \le \dlw$ is to be fixed below, 
according to \eqref{eq_delta}. 
We will see nonetheless that the choice of $\delta$ 
and the upper-bound can be stated independently 
of $m_1$
provided $m_1\ge 1$ 
(the larger is $m_1$, the better is the estimate).

\noindent
\textbf{Step 1:} We introduce a crucial minoration of the process $X\eD_0$ on the following event $\cE\eD_0$:
\begin{equation*}
\cE\eD_0 
:= \Big\{\Tsup{s\le t} X\eD_0(s)\cdot M\eD_1(s) \le 2\, \dlw\Big\}
\cap \Big\{\Tsup{s\le t} X\eD_0(s) \le \tdiv{2}\Big\}.
\end{equation*}
On the event 
$\cE\eD_0$,
the next inequalities hold for any $s\le t$: 
\begin{equation*}
	(\alpha\cdot M\eD_1(s) - \lambda)\cdot X\eD_0(s) \le 2\,\alpha\, \dlw
	\le  \frac{1-X\eD_0(s)}{4}.
\end{equation*}
A.s.  on the event $\cE\eD_0$,
the process $X\eD_0$ is thus upper-bounded
on the time-interval $[0, t]$
 by the solution $Y$ to
the following SDE,
with $y_0:= \delta /m_1$
and thanks to Lemma~\ref{lem_comp}:
\begin{equation*}
\RMd Y(s) = \frac{1-Y(s)}{4}\,  \RMd s
+ \sqrt{Y(s)\cdot (1-Y(s))}\, \RMd B_0(s)\mVg \quad Y(0) = y_0.
\end{equation*}
The main interest of this upper-bound is that it is explicitly given.
For any $y>0$, we denote $\tau^Y_y:= \inf\{t\ge 0\pv Y(t) = y\}$,
i.e. the hitting time of $y>0$ by the process $Y$.
We consider the following martingale process,
defined for any $s\in [0, \tau_0^Y)$:
\begin{equation*}
\RMd	\cN(s):= \sqrt{\frac{1-Y(s)}{Y(s)}}\, \RMd B_0(s),
\quad \cN(0) = 0.
\end{equation*}
Thanks to Itô's lemma,
$Y(s)$ is expressed 
for any $s\in [0, \tau_0^Y)$ as follows in terms of $\cN$:
\begin{equation*}
	Y(s):= y_0\, \exp\Big( \cN(s)
	- \frac{\LAg \cN\RAg_s}{4} \Big).
	\label{eq_def_Y}
\end{equation*}
Let us consider the random time change $(\xi_u)_{u>0}$ 
as follows:
\begin{equation*}
\xi(u)
:= \inf\big\{s\in [0, \tau_0^Y);\, \LAg \cN\RAg_s > u\big\},
\quad \mathcal G_u := \cF_{\xi_u},
\end{equation*} 
with value $\infty$ if the above set is empty.
Thanks to \cite[Theorem 18.4]{Ka02},
there exists a Brownian motion $W$ with respect to a standard extension of $\cG$
such that a.s. $W = \cN\circ  \xi$ on $[0, \LAg \cN\RAg_{\tau_0^Y})$
and $\cN(s) = W(\LAg \cN\RAg_s)$ for any $s\in [0, \tau_0^Y)$.
Note that $\LAg \cN\RAg_{\tau_0^Y}$ is here defined as the left limit at $\tau_0^Y$
of the non-decreasing process $\LAg \cN\RAg$.
Since the process $\LAg \cN\RAg$ is actually increasing in the time-interval $[0, \tau_0^Y)$,
$\LAg \cN\RAg_{\xi(u)} = u$ holds for any $u\in [0, \LAg \cN\RAg_{\tau_0^Y})$.
The following identity thus holds for any $u\in [0, \LAg \cN\RAg_{\tau_0^Y})$:
\begin{equation*}
Y\circ \xi(u)
= y_0\, \exp\big(W_u - u/4\big).
\end{equation*}

\noindent
\textbf{Step 2:} We prove that $\LAg \cN\RAg_{\tau_0^Y} = \infty$ holds a.s.
on the event $\{\tau_0^Y < \tau_1^Y\}$.

Assume by absurdum that there exists $A>0$ 
such that $\PR(\LAg \cN\RAg_{\tau_0^Y} \le A; \tau_0^Y < \tau_1^Y) >0$ and let us denote this quantity $\eps \in (0, 1)$.
Note that $0$ is a regular reflecting boundary for $Y$,
while $1$ is an exit boundary,
thanks to Lemma~\ref{bound_oneD}$(ii)$ and Corollary~\ref{cor_bound_oneD}.
There exists $r>0$ such that 
$\PR(r\le \tau_0^Y < \tau_1^Y) \le \eps/2$.
We consider the following sequence of stopping times,
in terms of an integer $\ell$:
\begin{equation*}
T_\ell := \inf\{s\ge 0;\, \LAg \cN\RAg_s\ge A \text{ or } |\cN(s)|\ge \ell\}.
\end{equation*} 
By continuity of $\LAg \cN\RAg$ that starts from 0, 
we know that the following property holds for any $s>0$ and any $\ell\ge 1$:
\begin{equation*}
\bE[\LAg \cN\RAg_{s\wedge T_\ell}] \le A.
\end{equation*}
On the other hand, the local martingale $(\cN(s))_{s}$
induces for any $\ell$ a continuous martingale $(\cN(s\wedge T_\ell))_{s}$
 whose predictible quadratic variation
is exactly $(\LAg \cN\RAg_{s\wedge T_\ell})_{s}$.
Thanks to Doob's inequality, 
the following inequalities thus hold 
for any $s>0$, $\ell\ge 1$ and $B>0$:
\begin{equation*}
\begin{split}
\PR(\Tsup{u\le s} |\cN(u\wedge T_\ell)|\ge B)
&\le 4\frac{\bE[\LAg \cN\RAg_{s\wedge T_\ell}]}{B^2}
\le \frac{4 A}{B^2}.
\end{split}
\end{equation*}
With $B = 4\sqrt{A/\eps}$ and any $\ell \ge B$,
it yields the following upper-bound:
\begin{equation*}
		\PR\Big(\Tsup{u\in \xi([0, A])} |\cN(u)|\ge B\Big)
		\le \eps/4.
\end{equation*}
By virtue of the properties of $\eps$ and $r$,
it entails that the following event has a positive probability, 
namely larger than $\eps/4$:
\begin{equation*}
\cE :=\{\LAg \cN\RAg_{\tau_0^Y} \le A\}
\cap \{\tau_0^Y < \tau_1^Y\wedge  r\}
\cap \{\Tsup{u\in \xi([0, A])}|\cN(u)|< B\}.
\end{equation*}
Yet, on this event $\cE$, it holds for any $s<\tau_0^Y$
that $\LAg \cN\RAg{s}\le A$ and $|\cN(s)|< B$, 
so that $Y(s) \ge y_0\exp(-A/4 - B).$
Yet, since $\tau_0^Y<\infty$, 
the continuity of the process $Y$
implies that $Y(\tau_0^Y-) = 0$, which is in contradiction
with the previous statement.
The event $\cE$ is therefore empty, which contradicts 
that $\cE$ has a positive probability.
It concludes that $\LAg \cN\RAg_{\tau_0^Y} = \infty$ holds a.s.
on the event $\{\tau_0^Y < \tau_1^Y\}$.

By the way, 
since $1$ is an absorbing boundary for $Y$
and by definition of $\LAg \cN\RAg$,
$\LAg \cN\RAg_{\tau_0^Y} = \LAg \cN\RAg_{\tau_1^Y}<\infty$
holds a.s. on the complementary event $\{\tau_1^Y < \tau_0^Y\}$.
\medskip

\noindent
\textbf{Step 3:} We justify the choice of the constants involved
in the definition of an event $\A\eD$ 
on which survival of $X\eD$ holds a.s. up to time $t$.

Let us consider some real number $\mu>0$, to be fixed below according to \eqref{eq_mu}.
On the event $\{\tau^Y_0< \tau^Y_{y_0+\mu}\}$, 
$\LAg \cN\RAg_{\tau_0^Y} = \infty$ as a consequence of step 2,
so that the following inequality holds 
for any $u\ge 0$: 
\begin{equation*}
	y_0\cdot \exp(W(u) - u/4)
= Y(\xi(u)) < y_0 + \mu.
\end{equation*}
It implies that the following inequality holds a.s. 
on the event $\{\tau^Y_0< \tau^Y_{y_0+\mu}\}$:
\begin{equation*}
	\int_0^\infty 
	\dfrac{y_0\, \exp\lc W(r) - r/4\rc}{1- y_0\, \exp\lc  W(r) - r/4\rc} \RMd r
	\le \dfrac{y_0}{1- y_0 - \mu} \int_0^\infty \exp\lc W(r) - r/4\rc \RMd r.
\end{equation*}
Note also that $\tau^Y_0 = \lim_{u\ifty} \xi(u)$
and the following expression for the derivative of $\xi$:
\begin{equation*}
	\xi'(r)
	= \dfrac{1}{\LAg \cN\RAg'\circ \xi(r)}
	= \dfrac{y_0\, \exp\lc W(r) - r/4\rc}{1- y_0\, \exp\lc W(r) - r/4 \rc}.
\end{equation*}
Therefore, the next  inequality holds 
for any $t\ge 0$:
\begin{equation}
	\PR_{y_0} (t< \tau^Y_0< \tau^Y_{y_0+\mu})
	\le \PR \lp \frac{t\cdot (1-y_0-\mu)}{y_0} <  \int_0^\infty \exp\lc W(r) - r/4\rc \RMd r \rp.
	\label{M_tY0}
\end{equation}
On the event $\{\tau^Y_{y_0+\mu} < \tau^Y_0\}$,
$y_0 + \mu= y_0\cdot \exp(W(r) - r/4)$
 for $r = \LAg \cN\RAg_{\tau^Y_{y_0+\mu}}\in (0, \infty).$
We thus deduce the following upper-bound:
\begin{equation}
\PR_{y_0}\Big(\tau^Y_{y_0+\mu} < \tau^Y_0\Big)
= \PR \lp (y_0+\mu)/y_0 \le \sup_{r\ge 0} \exp\lc W(r) - r/4\rc\rp.
\label{M_tYmu}
\end{equation}

Let $\eps >0$.
Since the law of $W$ is the one of a Brownian motion,
$W(r) / r \cvifty{t} 0$ holds a.s.
We can thus choose $c_1, c_2> 1$ such that:
\begin{equation*}
\PR \lp c_1 < \int_0^\infty \, \exp\big( W(r) - r/4\big), \RMd r \rp \le \eps \mVg
\quad
\PR \lp c_2 < \sup_{r\ge 0} \exp\lc W(r) - r/4\rc \rp \le \eps.
\end{equation*}
Thanks to  \cite[Lemma 3.2]{AP13}, we can choose $c_3>0$
such that 
the following inequality holds for any $x\in \cX_d$:
\begin{equation*}
\PR_x\Big(\Tsup{s\le t} M\eD_1(s) - M\eD_1(0) \ge \lambda\, t + c_3\Big) \le \eps.
\end{equation*}
This motivates the following choices for $m_1'$ and $\mu$:
\begin{equation}
m'_1:= m_1 + \lambda\, t + c_3>1,\qquad
\mu:= \frac{\dlw }{m'_1}.
\label{eq_mu}
\end{equation}
We then choose $\delta\le \dlw$ sufficiently small to ensure
the following two inequalities:
\begin{equation}
\frac{t\cdot (1- \dlw)}{\delta} - t
\ge c_1,\quad
\frac{\dlw}{\delta}\cdot (1+\lambda\, t + c_3)^{-1} \ge c_2.	
	\label{eq_delta}
\end{equation}
These conditions are prescribed in order to ensure 
the following two inequalities, recalling that $y_0 = \delta/m_1$
and $m'_1\ge m_1\ge 1$: 
\begin{equation*}
\frac{t\cdot (1-y_0-\mu)}{y_0} 
= \frac{t m_1}{\delta}\cdot (1- \frac{\dlw}{m'_1}) - t
\ge  c_1,\quad	
\frac{y_0+\mu}{y_0} \ge 
\frac{\dlw}{\delta}\cdot\frac{m_1}{m'_1} 
\ge  c_2.
\end{equation*}
Therefore, thanks to \eqref{M_tY0}, to \eqref{M_tYmu} and to the above choices 
of the constants,
$\PR_x(\A\eD) \ge 1- 3\, \eps$, where:
\begin{equation*}
\A\eD
:= \Lbr \Tsup{s\le t} M\eD_1(s) \le m_1'\Rbr
\cap \Lbr \tau^Y_0 < t\wedge \tau^Y_{y_0+\mu}\Rbr.
\end{equation*}
\textbf{Step 4:} We check that $t<\ext\eD$ holds a.s. on the event $\A\eD$.

To check the upper-bound by $Y$, 
let $T\eD_B:= \inf\{s\ge 0\pv X\eD_0(s)\cdot M\eD_1(s) \ge 2\, \dlw\}$.
Then, a.s. on the event $\A\eD$, 
both $M\eD_1(s) \le m'_1$
and $X\eD_0(s) \le Y(s) \le y_0+\mu$
hold 
for any $s\le t\wedge T\eD_B$.
Since $\mu/y_0\ge c_2> 1$ and $\mu:= \dlw/m'_1$,
the process $X\eD_0\cdot M\eD_1$ is upper-bounded as follows,
with the constant $c_4 := (y_0+\mu)\cdot m'_1$ 
and for any $s\le T\eD_B$:
$$X\eD_0(s)\cdot M\eD_1(s) \le c_4<2 \dlw.$$
By continuity of $X\eD_0\cdot M\eD_1$, the event $\{T\eD_B < t\}$ 
has an empty intersection with $\A\eD$.
Therefore $X\eD_0(s) \le Y(s)$ holds for any $s\le t$,
and in particular $\ext\eD\le \tau^Y_0 \le t$. 
So the following inequality holds
for any $x\in \cX_d$ such that $m_1\ge M_1(x)$ and $x_0\cdot m_1 \le \delta$:
\begin{equation*}
\PR_x (\ext\eD \le t) \ge \PR_x(\A\eD) \ge 1- 3\, \eps.
\end{equation*}
Though the event $\A\eD$ shall be adjusted 
depending on $m_1$,
 the choice of $\delta$ as a function of $\eps$
 (in \eqref{eq_delta}) is made independently of $m_1$,
 which concludes the proof of Lemma~\ref{M_Lx0}.
 \epf

%

\subsubsection
{Step 1.2: quick descent of the first moment}
\label{sec_Lm1} 
Provided that the optimal subpopulation is still significant,
 the first moment is unlikely to stay high for a significant time,
 as stated in  the upcoming Lemma~\ref{M_Lm1},
 whose proof is the purpose of this Subsection \ref{sec_Lm1}:
\begin{lem}
	\label{M_Lm1}
	Given any $t, y_0>0$,
the following supremum tends to 0 
as $m_{1}>0$ tends to infinity:
	\begin{equation*}
		\sup\Lbr \PR_x\lp t \le \tau\euD_{m_1}\wedge \ext\eD\rp
		\Big \vert d\in \N, x\in \cX_d,
		x_0\ge y_0\Rbr.
	\end{equation*}
\end{lem}

\bpf:
A.s. on the event $\{\inf_{\{s\le t\}} M\eD_1(s) \ge m_1\}$,
the process $X\eD_0$,
namely the initial component of the solution to \eqref{Sd}-\eqref{eq_def_Bi},
is  lower-bounded on $[0, t]$ 
by the solution $Y$ to the following SDE, thanks to Lemma \ref{lem_comp}:
\begin{equation*}
\RMd Y(s) = r(m_1)\, Y(s)\, \RMd s 
+ \sqrt{Y(s)\cdot(1-Y(s))}\, \RMd B_0(s)\mVg \quad Y(0) = y_0,
\end{equation*}
where $r(m_1):= \alpha\, m_1 -\lambda \cvifty{m_1}\infty.$

Since $M\eD_1(s) = 0$ whenever $X\eD_0(s) = 1$, this lower-bound of $X\eD_0$ by $Y$
must have stopped before $T^Y_1:= \inf\{ t\ge 0\pv Y(t) \ge 1\}$,
so $\tau\euD_{m_1}\wedge \ext\eD \le T^Y_1$.
We thus only have to prove that 
$\PR(t < T^Y_1)$ tends to 0 as $m_1$ tends to infinity.

Let $\eps, t_1 >0$.
Let us denote by $\cN^Y$ the following martingale process
naturally associated to $Y$:
\begin{equation*}
\cN^Y(t) := \int_0^t \sqrt{Y(s)\cdot(1-Y(s))}\, \RMd B_0(s).
\end{equation*}
The quadratic variation of $\cN^Y$ until time $t_1\le t$
is upper-bounded by $t_1/4$,
so that Doob's inequality implies the following upper-bound:
\begin{equation}\label{M_supMP}
\textstyle
\PR_{y_0}\Big(\sup_{s\le t_1} |\cN^Y(s)| > \dfrac{y_0}{2}\Big)
\le \dfrac{16 t_1}{y_0^2}.
\end{equation}
By choosing $t_1$ sufficiently small, we assume $16 t_1/y_0 ^2 \le \eps$.
On the complementary event $\Lbr \Tsup{s\le t_1} |\cN^Y(s)| \le y_0/2\Rbr$,
$Y$ stays above $y_0/2$ on the time-interval $[0, t_1]$.
The drift term can thus be lower-bounded by $ s\cdot r(m_1) \cdot y_0/2$
for any $s\le t_1 \wedge T^Y_1$.
Since it cannot exceed $1- y_0/2$ before $T^Y_1$,
it necessarily implies that for $r(m_1)$ sufficiently large
(that is $m_1$ sufficiently large),
we must have $T^Y_1 < t_1$ on the event $\Lbr \sup_{s\le t_1} |\cN^Y(s)| \le y_0/2\Rbr$.
Thanks to \eqref{M_supMP} and since $t_1\le t$,
this implies  that $\PR(t < T^Y_1)$ tends to 0 as $m_1$ tends to $\infty$
and concludes the proof of Lemma~\ref{M_Lm1}.
\epf


\subsubsection
{Step 1.3: quick descent of the next moments}
\label{sec_LmK}
Provided that one of the moment is initially upper-bounded,
it is unlikely for the next moment to stay high
on a significant time-interval afterwards,
 as stated in  the upcoming Lemma~\ref{M_LmK}
 (recall the notation \eqref{eq_def_tauMk}),
whose proof is the purpose of this Subsection \ref{sec_LmK}:
\begin{lem}
	\label{M_LmK}
	Given any integer $k\ge 1$ and any $t,\, m>0$,
the following supremum tends to 0 
as $m'>0$ tends to infinity:
	\begin{equation*}
		\sup\Lbr \PR_x\lp t\le \tau\eDuk_{m'}\wedge \ext\eD \rp
		\Big \vert d\in \N,\, x\in \cX_d,\,
		M_k(x)\le m\Rbr.
	\end{equation*}
\end{lem}

For the proof of Lemma \ref{M_LmK},
we exploit the following properties 
on the semi-martingale decomposition of the process $M\eD_k$,
summarised in the upcoming Lemma~\ref{lem_MkSM}.
\begin{lem}
\label{lem_MkSM}
For any $d\in \II{1, \infty}$, $k\in \Z$,
and any $x$ which belongs either to $\cX_d$ if $d\in \N$
or to  $\cX^{2k}$ if $d = \infty$,
the process $M\eD_k$ can be decomposed as follows
under $\PR_x$:
\begin{equation*}
	\RMd M\eD_k(t) 
	= V\eD_k(t)\, \RMd t + \RMd \M\eD_k(t),
	\quad M\eD_k(0) = M_k(x),
\end{equation*}
where $\M\eD_k$ is a continuous local-martingale starting from 0,
whose quadratic variation is 
$$\langle\M\eD_k\rangle_t
=\int_0^t (M\eD_{2\,k}(s)- M\eD_{k}(s)^2)\, \RMd s,$$
and $V\eD_k$ is a bounded variation process.
$(M\eiD_{2k}(t))_{t\ge 0}$
is a.s. locally upper-bounded.
In addition, there exists a universal constant $C_k>0$ 
such that the following upper-bounds hold a.s. 
for  any time-interval:
\begin{equation*}
\begin{split}
	V\eD_k &\le  \alpha \cdot(M\eD_1\cdot M\eD_k - M\eD_{k+1})
+ \lambda\cdot (C_k M\eD_{k} + 1)
\\&\le \lambda\cdot (C_k M\eD_{k} + 1)
\end{split}
\end{equation*}
\end{lem}
Note that given the quadratic variation of $\M\eD_k$,  
it is a martingale provided $d\in \N$.
\\

\bpf{of Lemma \ref{lem_MkSM}}:
Given the definition of $M\eD_k$
in terms of the solution $(X\eD)$ to the system of SDEs \eqref{Sd},
we obtain the following expression:
\begin{equation}
	V\eD_k
	:= \alpha\cdot (M\eD_1\cdot M\eD_k - M\eD_{k+1})
	+ \lambda \sum_{\ell= 0}^{d-1} (\ell + 1)^{k} X\eD_{\ell} 
	- \lambda\cdot (M\eD_k - \idc{d<\infty} d^k X\eD_d).
	\label{M_KVd}
\end{equation}
Thanks to Hölder's inequality, 
the next inequalities hold for any $x\in \cX_d$:
\begin{equation*}
	M_1(x) \le (M_{k+1}(x))^{\sfrac{1\!}{\!(k+1)}},
	\qquad M_k(x) \le (M_{k+1}(x))^{\sfrac{k}{\!(k+1)}},
\end{equation*}
thus $M_1(x)\cdot M_k(x) \le M_{k+1}(x).$
It thus implies the inequality $M\eD_1\cdot M\eD_k \le M\eD_{k+1}$
between the stochastic processes involved in \eqref{M_KVd}.

Exploiting that  $(\ell + 1)^k \le 2^k\cdot \ell^k$ for $\ell \ge 1$,
and that $X\eD_0\le 1$ for $\ell = 0$,
it yields the following inequality, with $C_k  =  2^k$:
\begin{equation*}
	V\eD_k
	\le \lambda\cdot (C_k\cdot M\eD_k + 1).
\end{equation*}
On the other hand, 
the local martingale term is defined as follows:
\begin{equation*}
\RMd \M\eD_k(t)
:= \sum_{i = 1}^{d} i^k \sqrt{X\eD_i(t)} \RMd W_i(t)
- M\eD_k(t) \RMd W\sD(t),
\quad \M\eD_k(0) = 0.
\end{equation*}
Thanks to the independence between the Brownian motions $(W_i)_i$
and to \eqref{eq_def_WsD},
its quadratic variation
satisfies the following identities:
\begin{equation*}
\begin{split}
	\RMd \LAg\M\eD_k\RAg_t
&= (M\eD_{2k}(t)
+ [M\eD_k(t)]^2)\, \RMd t
- 2 M\eD_k(t) \sum_{i = 1}^{d} i^k \sqrt{X\eD_i(t)} \, \RMd \LAg W_i, W\sD\RAg_t
\\&= (M\eD_{2k}(t)
- [M\eD_k(t)]^2)\, \RMd t.
\end{split}
\end{equation*}
In the case where $d=\infty$,
since $M\eiD_{2k}(x)<\infty$, 
we know thanks to \cite[Theorem 3]{AP13} 
that $(M\eiD_{2k}(t))_{t\ge 0}$
is a.s. locally upper-bounded, i.e.
that  $\sup_{s\le t} M\eiD_{2k}(s)<\infty$ holds a.s.
To be precise, 
	this is not stated directly in \cite[Theorem 3]{AP13},
	yet visibly obtained in the proof, 
	first in the case where $\alpha = 0$,
	then exended for any $\alpha > 0$
	because the Girsanov change of density 
	is uniformly upper-bounded
	(see Proposition \ref{pr_Fin_mom} in the appendix for completeness).

This expression for the quadratic variation is thus well-defined.
This concludes the proof of Lemma \ref{lem_MkSM}.\epf
\\

\bpf{of Lemma \ref{M_LmK}}:
Let $k\ge 1$, $t>0$,  $m>0$,
and $\eps >0$. Note that $M_1(x) \le M_k(x)$ holds for any $d\in \N$ and any $x\in \cX_d$.
Thanks to  \cite[Lemma 3.2]{AP13}, 
we can thus choose  $m_{1}>0$
 such that  the following inequality holds for any $d\in \N$ and any $x\in \cX_d$ such that  $M_k(x) \le m$:
\begin{equation}
	\PR_x(T\euD_{m_{1}} < t)
	\le \eps,
	\label{eq_defm}
\end{equation}
where $T\euD_{m_{1}}
:= \inf\{t\ge 0\pv M\eD_1(t) \ge m_{1} \}$,
namely
$T\euD_{m_{1}}$ is the hitting time of $m_{1}$ by $M\eD_1$. 

Thanks to the analysis of the process $M\eD_k$ conducted in Lemma \ref{lem_MkSM},
the following inequality holds then with $\check C_k=\alpha\cdot m_{1} + \lambda C_k$ 
for any $d\in \N$ and any $x\in \cX_d$:
\begin{equation*}
\begin{split}
	\bE_x[M\eD_{k}(t\wedge T\euD_{m_{1}})]
&\le M_k(x) 
- \alpha \cdot\bE_x\lp\int_0^{t}
\idc{s \le T\euD_{m_{1}}}
 M\eD_{k+1}(s) \, \RMd s\rp
\\&\quad +\check C_k\cdot \bE_x\lp \int_0^{t}
\idc{s \le T\euD_{m_{1}}} M\eD_{k}(s) \, \RMd s\rp.
\end{split}
\end{equation*}
Since $M\eD_{k}$ is a non-negative process,
with $\bar C_k= \check C_k/\alpha$, 
the following upper-bound on the $k+1$-th moment
thus holds for any $d\in \N$ and any $x\in \cX_d$ such that  $M_k(x) \le m$:
\begin{equation*}
\bE_x\lp \int_0^{t}\idc{s \le T\euD_{m_{1}}} M\eD_{k+1}(s) \, \RMd s\rp
\le \frac{m}{\alpha}
+ \bar C_k\cdot \bE_x\lp\int_0^{t} \idc{s \le T\euD_{m_{1}}} M\eD_{k}(s) \, \RMd s\rp.
\end{equation*}
By immediate induction over $k$, 
there exists $\hat C_k>0$ such that the following inequality holds 
for any $d\in \N$ and any $x\in \cX_d$ such that  $M_k(x) \le m$:
\begin{equation*}
\begin{split}
	\bE_x\lp  \int_0^{t}\idc{s \le T\euD_{m_{1}}} M\eD_{k+1}(s) \, \RMd s  \rp 
&\le \frac{(k-1)\cdot m}{\alpha}
+ \hat C_k\cdot \bE_x\lp \int_0^{t}\idc{s \le T\euD_{m_{1}}} M\eD_{1}(s) \, \RMd s  \rp
\notag\\&
\le \frac{(k-1)\cdot m}{\alpha}
+ \hat C_k\cdot t\cdot  m_{1}.
\end{split}
\end{equation*}
Thanks to Markov's inequality, 
we can thus choose $m'>0$ such that the following inequality holds 
for any $d\in \N$ and any $x\in \cX_d$ such that  $M_k(x) \le m$:
\begin{equation}
\PR_x\lp \int_0^{t}\idc{s \le T\euD_{m_{1}}} M\eD_{k+1}(s) \, \RMd s 
\ge t\cdot m' \rp
\le \eps.
\label{eq_majMkp}
\end{equation}
Recalling \eqref{eq_defm},
this concludes the proof of Lemma \ref{M_LmK},
since:
\begin{equation*}
\Lbr t> \tau\eDuk_{m'} \Rbr
\cap \Lbr t\le T\euD_{m_{1}}\Rbr
\subset \Lbr  \int_0^{t}\idc{s \le T\euD_{m_{1}}}  M\eD_{k+1}(s) \, \RMd s 
\ge t\cdot m'\Rbr.
\SQ
\end{equation*}

\subsubsection
{Concluding the proof of Proposition \ref{M_Pm4}}
\label{sec_fPm4}
Let $t, \eps>0$. Thanks to Lemma \ref{M_Lx0}, 
we choose an upper-bound $\delta\in (0, \tdiv{(16\alpha)})$ 
for the initial condition of the process $X_0\eD\cdot M\eD_1$,
so that the following inequality holds, for any $d\in \N$ and any $x\in \cX_d$
such that both $M_1(x) \ge 1$ and $x_0\cdot M_1(x) \le \delta$:
\begin{equation}
 \PR_x\lp t<\ext\eD\rp
	\le \eps.
	\label{M_Pmort}
\end{equation}
We consider the following exit time:
\begin{equation}\label{M_T01}
	T\ekD B_\delta:=  \inf\{t\ge 0\pv X\eD_0(t)\cdot M\eD_1(t)\le \delta \}.
\end{equation}
Let $m_1^\vee = (\sfrac{2\, \lambda}{\alpha})\vee 1$.
We consider $\tau\euD_{m_1}:=  \inf\{t\ge 0\pv M\eD_1(X_t)\le m_1\}$
for any $m_1\ge m_1^\vee$
(so that $\tau\euD_{m_1}\le \tau\euD_{m_1^\vee}$).
%
Thanks to \eqref{M_Pmort} and to the strong Markov property
at time $T\ekD B_\delta$, 
the following inequality holds for any $x\in \cX_d$ and $d\ge1$:
\begin{equation}
	\PR_x \lp T\ekD B_\delta \le t \le \tau\euD_{m_1^\vee}
	\mVg 2\, t<\ext\eD\rp \le \eps.
	\label{M_ZM}
\end{equation}
On the event $\{t \le T\ekD B_\delta\wedge \tau\euD_{m_1^\vee} \wedge \ext\eD\}$,
recalling that $m_1^\vee \ge 2\, \lambda/\alpha$,
the following inequality holds
 for any $s\le t$:
\begin{equation*}
	(\alpha\cdot M\eD_1(s) -\lambda)\cdot X\eD_0(s) \ge \frac{\alpha\cdot \delta}{2}.
\end{equation*}

Thus, thanks to Lemma \ref{lem_comp},
$X\eD_0$ is  lower-bounded
by the solution $Y$  to the following SDE,
 a.s. on $[0, t]$ on the event $\{t \le T\ekD B_\delta\wedge \tau\euD_{m_1^\vee} \wedge \ext\eD\}$:
\begin{equation}
	\RMd Y(s) = \frac{\alpha\cdot \delta}{2}\cdot \RMd s
	+ \sqrt{Y(s)\cdot (1-Y(s))}\, \RMd B_0(s)\mVg \quad Y(0) = 0. 
	\label{M_Ydef}
\end{equation}
Thanks to Corollary~\ref{cor_bound_oneD}
the left-boundary $0$
 is regular reflecting for $Y$.
Thus we can choose
$y_0>0$ such that:
$\PR(Y(t) \le y_0) \le \eps.$
Recalling \eqref{M_ZM},
 this implies  the following inequalities for any $d\in \N$ and $x\in \cX_d$:
\begin{equation}\label{eq_min_xy}
\begin{split}
&\hcm{-0.3}\PR_x\Big(X\eD_0(t) \le y_0\mVg t\le \tau\euD_{m_1^\vee}\mVg 2\, t<\ext\eD\Big)
	\\& \le \PR_x\Big(T\ekD B_\delta \le t \le \tau\euD_{m_1^\vee}\mVg 2\, t<\ext\eD\Big)
	+ \PR_x\Big(Y(t) \le y_0\mVg  t \le T\ekD B_\delta\wedge \tau\euD_{m_1^\vee}\wedge \ext\eD\Big)
	\\&\le 2\, \eps.
\end{split}
\end{equation}

Thanks to Lemma \ref{M_Lm1}, we can then choose $m_1\ge m_1^\vee$ 
associated with $y_0$.
Thanks to~\eqref{eq_min_xy}
and to the Markov property at time $t$,
we deduce the following inequalities:
\begin{equation}
\begin{split}
		\PR_x\Big(2\,t<  \tau\euD_{m_1} \wedge \ext\eD\Big) 
	&\le \PR_x\Big(X\eD_0(t) \le y_0\mVg t\le \tau\euD_{m_1^\vee}\mVg 2\, t<\ext\eD\Big)
	\\&\hcm{1}
	+ \bE_x\Big[ \PR_{X\eD(t)}\big(t \le  \tau\euD_{m_1}\wedge \ext\eD\big) \pv X\eD_0(t) \ge y_0\Big]
	\\&\le 3\, \eps.
\end{split}
	\label{M_moneD}
\end{equation}

Thanks to Lemma \ref{M_LmK}, 
we can choose $m_2>0$ associated with $m_1$
and  $m_3>0$  associated with $m_2$ such that
the following inequalities hold:
\begin{equation*}
	\PR_x\Big(3\,t<  \tau\edD_{m_2} \wedge \ext\eD\Big) 
	\le 4\, \eps\mVg \quad
	\PR_x\Big(4\,t< \tau\etD_{m_3} \wedge \ext\eD\Big) 
	\le 5\, \eps.
\end{equation*}
Applying the same argument inductively
over $k\ge 3$, 
we can choose $m_k>0$ such that the following inequality holds
for any $d\in \N$ and any $x\in \cX_d$:
$$\PR_x\Big((k+1)\cdot t<  \tau\eDk_{m_k} \wedge \ext\eD\Big) 
\le (k+2)\cdot \eps,$$
so as to treat any moment.
This concludes the proof of Proposition \ref{M_Pm4}.
\epf

\subsection
[Prop:augm]{Step 2: rare large increase of the moment}
\label{sec_augm}

With a probability close to 1,
the increase of the moments
after their descent
can be upper-bounded uniformly over a given time-interval,
as stated in the next proposition,
whose proof is the purpose of this Subsection \ref{sec_augm}. 
For any $k\ge 1$ and  $m\ge 1$,
let: 
\begin{equation}
	T\eDk_{m}:=\inf\Big\{t\ge 0;\, M\eD_k(t) \ge m\Big\}.
	\label{eq_def_TkM}
\end{equation}
\begin{prop}
	\label{M_augm}
	For any $k>1,$\,  $t, m>0$,
	the following supremum tends to 0 
	as $m'>0$ tends to infinity:
	\begin{equation*}
\sup\Lbr
	\PR_x \lp T\eDk_{m'} \le t\rp 
	\bv d\in \N, x \in \cX_d, M_k(x) \le m\Rbr.
\end{equation*}	
\end{prop}

\bpf:
Let $k>1, t>0$, $m'\ge m$, $d\in \N$ 
and $x\in \cX_d$ such that $M_k(x)\le m$¨.


We aim at exploiting Doob's inequality on a non-negative sub-martingale
that is an upper-bound of $M\eD_k$.
Given the semi-martingale decomposition of $M_k$
as stated in Lemma~\ref{lem_MkSM},
we consider the solution $\widehat M\eD_k$ to the following SDE:
\begin{equation}
	\textstyle
	\widehat M\eD_k(t)
	:= m + \lambda t + \lambda C_k \int_0^t \widehat M\eD_k(s)\, \RMd s 
	+ \M\eD_k(t).
	\label{M_hMK}
\end{equation}
Thanks to \cite[Proposition 3.12]{PR14}
(as for Lemma \ref{lem_comp}),
$\widehat M\eD_k(t) \ge M\eD_k(t)$ holds for any $t\ge 0$.
Since $M\eD_k$ is non-negative,
$\widehat M\eD_k$ is also non-negative.
As the solution to equation~\eqref{M_hMK},
$\widehat M\eD_k$ is a sub-martingale.
Since $\bE_x\lc  \widehat M\eD_k(s)\rc$ is upper-bounded by $d^k$
for any $s$, 
Grönwall's lemma 
implies the following inequality for any $x\in \cX_d$ such that $M_k(x)\le m$:
\begin{equation}\label{M_GL}
	\Tsup{s\le t}\bE_x\lc  \widehat M\eD_k(s)\rc
	\le (m+ \lambda t)\, e^{C t}.
\end{equation}

Exploiting Doob's inequality on $\widehat M\eD_k$,
then \eqref{M_GL} 
with $C_M:= (1+ \lambda t)\, e^{C t}$,
we obtain the following inequality for any $x\in \cX_d$ such that $M\eD_k(x)\le m$:
\begin{align*}
	\PR_x\Big(\Tsup{s\le t} M\eD_k(s) > m'\Big)
	&\le \PR_x\Big(\Tsup{s\le t} \widehat M\eD_k(s) > m'\Big)
	\\&
	\le \dfrac{\bE_x[\widehat M\eD_k(t)]}{m'}
	\le \dfrac{C_M\, m}{m'}.
\end{align*}
This concludes the proof of Proposition \ref{M_augm}.
\epf

\subsection[Prop MT]{Concluding the proof of Proposition \ref{M_prop_MT}}
\label{M_sec_pMT}
First of all, we show that we have a uniform upper-bound on the extinction rate $\rho_0^{(d)}$ associated with the system 
\eqref{Sd}:
\mbox{$\Tsup{d\in \N} \rho_0^{(d)}<\infty$}.
Indeed, whatever $d\in \N$,
we can find $x\eD\in \cX_d$ such that $x\eD_0\ge \tdiv{2}$
so that,
thanks to Lemma \ref{lem_comp},
 $X\eD_0$ 
is a.s. lower-bounded under $\PR_x$ 
for any time by the solution $Y$ 
to the following SDE:
\begin{equation*}
\RMd Y(s) = -\lambda\cdot Y(s) \RMd s
	+ \sqrt{Y(s)\ (1-Y(s))}\, \RMd B_0(s)\mVg \quad Y(0) = \tdiv{2}. 
\end{equation*}
Thanks to Lemma \ref{bound_oneD}$(i)$,
	 the left-boundary $0$ 
	is an exit boundary.
The semi-group governing $Y$, 
with extinction at $\tau_0^Y$,
corresponds exactly to the system $\eqref{Sd}$
with $d = 1$, $\alpha = 0$, $X'_0 = Y$ and $X'_1 = 1 - Y$.
Thanks to Theorem \ref{M_ECVFin},
the semigroup thus displays QSC
with extinction rate $\rho_\vee$.
Denoting $\PR^Y_{\tdiv{2}}$ the law of $Y$,
it entails the following inequality from the convergences of the survival capacities:
\begin{equation*}
	\rho_0^{(d)} = \lim_{t\rightarrow \infty} \tfrac{- 1}{t}\, 
	\log \PR_{x^{(d)}} (t< \ext\eD) 
	\le \lim_{t\rightarrow \infty}\tfrac{ - 1}{t}\, 
	\log \PR^Y_{\tdiv{2}} (t< \tau_0^Y):=  \rho_\vee.
\end{equation*}
Thanks to Proposition \ref{M_Pm4}
 (similarly to the way Proposition~\ref{M_Tdy} was deduced),
we can choose  $m>0$ such that $\tau\euD_{m}$ satisfies
the following inequality for any $d\ge1$ and any $x \in \cX_d$:
\begin{equation*}
\bE_x\exp\Big[ (\rho_\vee + 1)\cdot (\tau\euD_{m}\wedge \ext\eD)\Big] \le C< \infty.
\end{equation*}
In particular, it implies the following inequality for any $t>0$ and any $d\in \N$:
\begin{equation}
	\PR_{\nu\eD} \Big(t <  \tau\euD_{m} \wedge  \ext\eD\Big) 
	\le C\, \exp\Big[-(\rho_\vee+1)\cdot t\Big].
	\label{M_TM}
\end{equation}
Then, for any $\eps>0$, consider $t:= -\log(\eps/(2\, C))$.
Thanks to Proposition \ref{M_augm},
the probability of large increase of the process $M\eD_k$ can be 
made negligible. 
So we can choose a constant $m'>0$ such that 
the following inequality holds
for any initial condition $x\in \cX_d$ such that $M_k(x) \le m$:
\begin{equation}
\PR_x \lp\Tsup{s\le t} M\eD_k(s) \ge m'\rp
	\le \eps/2\, \exp[-\rho_\vee\, t].
	\label{M_MkB}
\end{equation}
By virtue of the definition of $t$,
recalling (\ref{M_TM}) and (\ref{M_MkB}),
the following inequalities 
thus hold for any $d\in \N$: 
\begin{align*}
	\nu\eD\Big(\Big\{M_k \ge m'\Big\}\Big)
	&= \exp[\rho_0^{(d)}\, t]\cdot 
	\PR_{\nu\eD} \Big(M\eD_k(t) \ge m' \,;\, t\le \ext\eD\Big)
	\\&
	\le \exp[\rho_\vee\, t]\cdot
	\Big( \PR_{\nu\eD} \Big(\tau\euD_{m} > t\,;\, t< \ext\eD\Big) 
	\\&\hspace{1cm}
	+ \bE_{\nu\eD} \Big[\PR_{X\eD(\tau\euD_{m})} \Big(\Tsup{s\le t} M\eD_k(s) \ge m'\Big)
	\,;\, \tau\euD_{m} < t\wedge \ext\eD\Big]
	\Big)
	\\&
	\le C\cdot \eps/(2\, C) + \eps/2 \le \eps.
\end{align*}
This concludes the proof of Proposition \ref{M_prop_MT}.
\hfill $\square$

\begin{rem}
In fact, thanks to Lemma \ref{M_LmK},
 we are thus able to upper-bound any moment,
	with probability close to 1
	and under the QSD $\nu\eD$
	uniformly on $d\in \N$.
	These upper-bounds will extend to the limiting QSD on $\cX_\infty$.
\end{rem}

\section{Proof of Theorem \ref{M_ECVdInf}: the infinite dimensional case}
\label{M_sec_ecvdinf}


The proof of Theorem \ref{M_ECVdInf}
is achieved in Subsection \ref{M_sec_TECV}
by ensuring Assumption $({\bf AF})$
as stated in Subsection~\ref{M_genSt}.
We will treat both the case of large yet finite values of $d$
and $d= \infty$, for which we recall that any $x\in \cX_\infty$
has a finite sixth moment (see \eqref{eq_def_Xinf}).

	As one can imagine, 
	the proof of Theorem \ref{M_ECVdInf}
	is much more technical than the ones of Theorem \ref{M_ECVdisc},
	Theorem \ref{M_ECVFin} and Proposition \ref{M_prop_MT}. 
	For instance, 
	there is no explicit reference measure 
	that seems to be exploitable as $\alc^{(\infty)}$:
	the Lebesgue measure cannot be extended 
	on an infinite dimensional space!
	Many elements of these previous proofs 
	are however to be exploited 
	with only slight adaptations, 
	so that the reader is really encouraged to read them before the next proof.
	
The core idea behind the proof is 
still that the individuals carrying many mutations
	are actually wiped out very rapidly,
	implying rapid shuffle of the last coordinates.
	Quite unexpectedly,
	the criteria we developed 
	to deal with jump events 
	has proved to be very effective in this context.
Notably, we could exploit the Girsanov transform to relate to 
the finite dimensional problem and
deal with moments increasing too largely
	as exceptional events.

\subsection{Outline of the proof}
\label{sec_}

We now consider $d\in \II{1, \infty}$, i.e. including the case $d=\infty$.
For the purpose of Theorem~\ref{M_ECVdInf}
in this Section \ref{M_sec_ecvdinf},
we replace the notation given in \eqref{M_Dn}
by the following one:
\begin{equation}
\cD\eD_\ell:= \Big\{x\in \cX_d\pv  
x_0 \ge \frac{1}{2 \ell}\Big\}.
\label{def_DL}
\end{equation}
In this approach, the family  $\cD\eD_\ell$ covers the whole state space:
\begin{equation*}
\medcup_{\ell\ge 1} \cD\eD_\ell
= \cX_d.
\end{equation*}

%

Because it is close to the previous proof of Proposition \ref{M_prop_MT},
we will first focus on the result of Theorem \ref{eT_dInf},
which can be interpreted as the statement of escape from the transitory domain
(though the issue of the dependency in the parameter $\eps$
led us to integrate the result into the following estimate of almost perfect harvest).
The sets $E\eD$ that we will consider are defined through 
three parameters $m$, $y, \eta>0$ as follow:
\begin{equation}
	E\eD 
	:= \Big\{x\in \cX_d\pv M_3(x) \le m\mVg
	\frl{j\le \Lfl m/ \eta \Rfl+1} x_j \ge y\Big\}.
	\label{M_cDE}
\end{equation}

The reference probability measure $\zeta\eD$ on $\cX_d$ is chosen 
to be especially adapted for our arguments,
in a way that makes it actually complex to express. 
Its specific definition, stated in \eqref{M_zeta}, is given 
in Subsection~\ref{sec_mixAcc} that is dedicated to the mixing property and the accessibility of the subsets of $\cX_d$.
It relies on the notations 
and properties introduced in Subsections \ref{M_sec_Agg},
where we justify a close connection with the finite dimensional SDEs.
The estimate of almost perfect harvest 
is then conducted in Subsection~\ref{sec_AFdI},
before we finally conclude the proof of Theorem~\ref{M_ECVdInf}
in Subsection~\ref{M_sec_TECV}.

As stated at the beginning of  Section \ref{M_sec_ECVFin},
for any $d\in \II{1, \infty}$ and 
probability measure $\zeta$ on $\cX_d$,
the process $X\eD$ is solution under $\PR_\zeta$
of the system \eqref{Sd} 
with initial condition $X\eD(0)$ distributed as $\zeta$.

\subsection{Escape from the Transitory domain}

We prove in this subsection that the process rapidly escape 
the sets $E\eiD$ as proposed in \eqref{M_cDE}, 
provided that the upper-bound on the moment is sufficiently large 
and the lower-bound on the optimal subpopulation size is sufficiently small,
as stated in the next theorem.
We recall the notation $\tau\eiD_{E}$ as the entry time of $E\eiD$.
\begin{theo}
	\label{eT_dInf}
For any  $t, \eta, \eps>0$, 
there exist a couple $\mpr, y>0$ such that
the following inequality holds 
with $E\eiD=E\eiD(\mpr, y, \eta)$ as in \eqref{M_cDE}
for any $d \in \II{1, \infty}$ and $x\in \cX_d$: 
\begin{equation}
\PR_x\lp t<\ext\eiD\wedge \tau\eiD_{E}\rp
	\le \eps.
	\label{eq_eTdInf}
\end{equation}	
\end{theo}
\begin{rem}\label{rem_eTdInf}
With the same reasoning as in the proof of Proposition \ref{M_Tdy},
\eqref{eq_eTdInf} has the following implication 
in terms of exponential moments.
For any $\rho, \eta>0$,
there exist a couple $\mpr, y>0$ such that
the following inequality holds
for any $d \in \II{1, \infty}$ and $x\in \cX_d$:
\begin{equation*}
\bE_x\lp \exp[\rho\cdot (\ext\eiD\wedge \tau\eiD_{E})]\rp
	\le 4.
\end{equation*}
\end{rem}

The proof of Theorem \ref{eT_dInf} 
relies on the six forthcoming lemmas,
mostly adapted from the uniform escape 
and the uniform descent of the moments
in the finite dimensional systems.
They are given in the order at which they will be exploited 
to conclude the proof in Subsection~\ref{sec_eT_dInf}.

We first show in the upcoming Lemma~\ref{M_Lx0i} that the click is very likely 
when the growth of the optimal subpopulation size is initially very small,
in the situation where the first moment is large 
so the initial size itself is small.
\begin{lem}
\label{M_Lx0i}
For any $t>0$,
the following supremum tends to 0
as 	$\delta$ tends to $0$:
\begin{equation*}
\sup\Lbr \PR_x\lp t<\ext\eiD\rp
\Big \vert  d \in \II{1, \infty}\mVg
x\in \cX_d\mVg
M_1(x) \in (1, \infty)\mVg
x_0\cdot M_1(x) \le \delta \Rbr.
\end{equation*}
\end{lem}
Then, provided the optimal subpopulation size is non-negligible,
we show that the first moment is unlikely to stay very high, 
as stated in the upcoming Lemma~\ref{M_Lm1i}.
\begin{lem}
\label{M_Lm1i}
For any two real numbers $t, y_0>0$,
the following supremum tends to 0 
as $m_1$ tends to $\infty$:
\begin{equation*}
\sup\Lbr \PR_x \lp t \le \tau\euD_{m_1}\wedge \ext\eiD\rp
\Big \vert  d \in \II{1, \infty}\mVg
x\in \cX_d\mVg
x_0\ge y_0\Rbr,
\end{equation*}
where we recall $\tau\euD_{m_1}:=  \inf\{t\!\ge\!0\pv M\eiD_1(t)\le m_1\}$.
\end{lem}
The proofs of Lemmas 
\ref{M_Lx0i} and 
\ref{M_Lm1i},
can be adapted mutatis mutandis from the ones 
of respectively Lemmas 
\ref{M_Lx0} (in Subsection~\ref{sec_Lx0})
and \ref{M_Lm1} (in Subsection~\ref{sec_Lm1}).

The next step is to handle the situation 
where both the first moment and the optimal subpopulation size
are initially small, as stated in the upcoming Lemma~\ref{M_LX0i}.
\begin{lem}
\label{M_LX0i}
For any two real numbers $t, m_1>0$, the following supremum tends to 0 
as 	$\delta$ tends to $0$:
\begin{equation*}
\sup\Lbr \PR_x\lp t<\ext\eiD\rp
\Big \vert  d \in \II{1, \infty}\mVg
x\in \cX_d\mVg
M_1(x) \le m_1\mVg
x_0 \le \delta \Rbr.
\end{equation*}
\end{lem}
\bpf: As a generalization of  Lemma \ref{M_Lx0e},
Lemma \ref{M_LX0i} is a consequence of the fact 
that $X\eiD_0$ is upper-bounded
on the event $\{\sup_{s\le t} M\eiD_1(s) \le m'_1\}$,
thanks to Lemma \ref{lem_comp},
by the solution $Y$ to the following SDE:
\begin{equation*}
	\RMd Y(t) = \alpha\, m'_1\cdot Y(t) \, \RMd t + \sqrt{Y(t)\cdot (1-Y(t))}\, \RMd B_0(t)
	\mVg \quad Y(0) = \delta.
\end{equation*}
Thanks to Lemma \ref{bound_oneD}$(i)$,
$0$ is an exit boundary for $Y$.
Thanks to \cite[Lemma 3.2]{AP13},
we know an upper-bound 
of $\PR_x(\sup_{s\le t} M\eiD_1(s) \ge m'_1)$
that tends to 0 as $m'_1$ goes to $\infty$,
uniformly in the $x\in \cX_d$ such that $M_1(x) \le m_1$.
The combination of these two facts concludes the proof.\hfill
\epf 
\\

As the next step, we justify with the upcoming Lemma~\ref{M_Lmki} 
that, once a moment has descended,
it is unlikely for  
the next moment to stay high on a significant time-interval afterwards:
\begin{lem}
	\label{M_Lmki}
	Given any integer $k\in \{1, 2, 3\}$ and any two real numbers $t,\, m>0$,
	the following supremum tends to 0 as $m'$ tends to $\infty$:
	\begin{equation*}
		\sup\Big\{ \PR_x \lp t \le \tau\eDuk_{m'}\wedge \ext\eiD\rp
		\Big \vert\; d \in \II{k, \infty}\mVg
		x\in \cX_d\mVg
		M_{2\, k}(x) < \infty\mVg
		M_k(x)\le m\Big\}.
	\end{equation*}
\end{lem}
\bpf: The proof of Lemma \ref{M_Lmki}
generalizes
the one of Lemma~\ref{M_LmK}
in Subsection~\ref{sec_LmK}.
We just sketch the localization argument
for the case $d=\infty$,
that is similar yet simpler 
than the one presented for the proof of the forthcoming Lemma \ref{M_moment_control}.

Thanks to Lemma \ref{lem_MkSM}, recall that 
$M\eI_{2k}$ is a.s. locally upper-bounded.
We can thus introduce a sequence $T_\ell$, $\ell\ge 1$, 
of stopping times such that $M\eI_{2k}$
is upper-bounded by $\ell$ on $[0, t\wedge T_\ell]$ and that goes to infinity 
as $\ell$ tends to infinity.
With the same arguments as for Lemma~\ref{M_LmK}:
\begin{equation*}
		\bE_x\lp  \int_0^{t}\idc{s \le T\euD_{m_{1}}\wedge T_\ell} M\eD_{k+1}(s) \, \RMd s  \rp 
		\le \frac{(k-1)\cdot m}{\alpha}
		+ \hat C_k\cdot t\cdot  m_{1}.
\end{equation*}
Taking $\ell$ to infinity by monotone convergence,
we can then proceed as previously
and conclude the proof of Lemma \ref{M_Lmki}.
\epf
\\

With the upcoming Lemma~\ref{M_Lymi},
we state that, on the event of its survival,
the process is bound to reach non-negligible 
subpopulation sizes 
for the $(X\eiD_j)_{j\in \II{0, J}}$, 
for any $J\le d$ (that are the $J$-optimal classes).
\begin{lem}
\label{M_Lymi}
Given any integer $J\in \N$, and any three real numbers $t,\, m_1, y_0>0$,
the following supremum tends to 0 as $y$ tends to 0:
\begin{equation*}
\sup\Lbr \PR_x\lp t\wedge \ext\eiD < \tau\ekD J_{y} \rp
\Big \vert  d \in \II{J, \infty}\mVg
x\in \cX_d
\mVg x_0\ge y_0
\mVg M_1(x)\le m_1\Rbr,
\end{equation*}
where we recall the notation
$\tau\ekD J_{y}:= \inf\big\{s\ge 0\pv \frl{j\le J} X\eiD_j(s) \ge y\big\}$
for any $y>0$ and any integer $J \ge 0$.
\end{lem}
The above notation is to be understood 
as $\tau\ekD J_{y} = \tau\ekD d_{y}$ for any $J\ge d$ (for the case $d\in \N$).
The proof of Lemma~\ref{M_Lymi} is an adaptation mutatis mutandis of the one of Lemma~\ref{M_Lxke},
as stated in Subsection~\ref{M_sec_etD}.



\subsubsection[Lem:moment control]{Upper-bound on the probability of moment increase}
\label{sec_mom_incr}
Subsection~\ref{sec_mom_incr} is devoted to the upcoming Lemma~\ref{M_moment_control}.
In the time-interval between the descent of the moment 
and the increase of the $J$-optimal population sizes $(X\eiD_j)_{j\le J}$,
we show that the control on the corresponding moment stays tight.
Recall the definition of $T\eDk_{m}$ from \eqref{eq_def_TkM}.
\begin{lem}
	\label{M_moment_control}
	For any two real numbers $k>1$ and $t>0$, 
	there exists $C\ge 1$
	such that 
	the following inequality holds
	for any $m, m'>0$, any $d\in \II{1, \infty}$ 
	and any $x\in \cX^{2k}$
	($x\in \cX_d$ for $d\in \N$)
	 such that $M_k(x)\le m$,
	$$\PR_x\lp  T\eDk_{m'}\le t \rp
	\le \dfrac{C m}{m'}.$$
\end{lem}
\bpf: The proof of Proposition \ref{M_augm}
already implies the result provided $d<\infty$.
We justify in the following 
 that it extends to the case where $d=\infty$
 in which the martingale part $\M\eI_k$ in Lemma \ref{lem_MkSM} is a priori only local. 

The initial condition $x\in \cX^{2k}$ is such that $M_k(x) \le m$.
The expression of $V\eI_k$ in \eqref{M_KVd} 
implies the following inequalities:
\begin{equation*}
	- \alpha M\eI_{k+1}
\le V\eI_k
\le \lambda\cdot (C_k M\eI_{k} + 1).
\end{equation*}
Thanks to Lemma \ref{lem_MkSM}
with the fact that $x\in \cX^{2k}$,
the process $M\eI_{k+1}$  is a.s. locally upper-bounded
(thus also $M\eI_{k}$).
Since $\M\eI_k(t) = M\eI_k(t) - M_k(x) - \int_0^t V\eI_k(s) \RMd s$,
$\M\eI_k$ is thus also 
a.s. locally upper-bounded.
Thanks to Duhamel's formula, 
this entails that the process 
$\widehat M\eI_k$ 
is well-defined as a solution to the following SDE,
similar to \eqref{M_hMK}:
\begin{equation*}
	\textstyle
	\widehat M\eI_k(t)
	:= m + \lambda t + \lambda C_k \int_0^t \widehat M\eI_k(s)\, \RMd s 
	+ \M\eI_k(t).
\end{equation*}
We localize this process thanks to the following sequence of stopping times
$T_\ell$, for any $\ell\ge 1$:
\begin{equation}
	T_\ell:= \inf\big\{ s \ge 0\pv \LAg\M\eI_k\RAg_s \ge \ell\mVg 
	M\eI_k(s) \ge \ell\big\}.
	\label{def_Tell}
\end{equation}
The process $(\widehat M\eI_k(t\wedge T_\ell))_{t\ge 0}$
defines a non-negative submartingale 
such that the following inequality holds for any $t\ge 0$:
\begin{equation*}
	\textstyle
	\bE_x[\widehat M\eI_k(t\wedge T_\ell)]
	:= m + \lambda t + \lambda C_k \int_0^t \bE_x[\widehat M\eI_k(s\wedge T_\ell)]\, \RMd s.
\end{equation*}
Thanks to Grönwall's lemma
(see for instance \cite[Proposition 6.59]{PR14}),
since our localization procedure entails an upper-bound on $(\widehat M\eI_k(t\wedge T_\ell))_{t\ge 0}$,
the following inequality holds for any $t\ge 0$:
\begin{equation}
	\bE_x\lc \widehat M\eI_k(t\wedge T_\ell)\rc
\le (m+ \lambda t)\cdot e^{\lambda C_k t}\le C\cdot m,
\label{eq_maj_MeIk}
\end{equation}
with $C := (1+ \lambda t)\cdot e^{\lambda C_k t}$
(recall that $m\ge 1$).
Thanks to \cite[Proposition 3.12]{PR14},
 $\widehat M\eI_k \ge M\eI_k$.
Thanks to Doob's inequality, 
the following inequalities thus hold
 for any $t\ge 0$ and $m'>0$:
\begin{equation*}
\begin{split}
		\PR_x\Big(\Tsup{s\le t\wedge T_\ell} M\eI_k(s) \ge m'\Big)
	&\le \PR_x\Big(\Tsup{s\le t\wedge T_\ell} \widehat M\eI_k(s) \ge m'\Big)
	\\&
	\le \dfrac{\bE_x[\widehat M\eI_k(t\wedge T_\ell)]}{m'}.
\end{split}
\end{equation*}
Recalling \eqref{eq_maj_MeIk}
and thanks to the monotone convergence theorem
letting $\ell$ tend to infinity,
it entails the following inequality
for any $m, m'>0$
and any $x\in \cX^{2k}$ such that $M_k(x) \le m$:
$$\PR_x\lp  T^{(k|\infty)}_{m'}\le t \rp
\le \dfrac{C m}{m'},$$
which concludes the proof of Lemma \ref{M_moment_control}.
\epf

\subsubsection[Th:eT.dInf]{Concluding the proof of Theorem \ref{eT_dInf}}
\label{sec_eT_dInf}
Let us consider any three real numbers $t, \eta, \eps >0$.
We consider the following event as a function of $m_1>0$
that describes a failure in the descent of the first moment:
\begin{equation*}
	\cE\eiD_1:= \big\{\tau\euD_{m_1} > 2t\big\}\cap
	\big\{2\, t < \ext\eiD\big\}.
\end{equation*}
With exactly the same reasoning as for Proposition  \ref{M_Pm4},
exploting Lemmas \ref{M_Lx0i} and \ref{M_Lm1i} instead of Lemmas \ref{M_Lx0} and \ref{M_Lm1},
we can choose $m_1>0$ such that
$\PR_x(\cE\eiD_1) \le 3\, \eps$
holds for any $x\in \cX_d$.

The following event is stated as a function of  $m'_1>0$
and describes a failure in having the first moment contained
on a significant time-interval:
\begin{equation*}
	\cE\eiD_2:=\Big\{\tau\euD_{m_1}\le 2t\Big\}
	\cap\Big\{2t < \ext\eiD\Big\}
	\cap\Big\{\wtd{T}\euD_{m'_1}  \le 5t\Big\},
\end{equation*}
where $\wtd{T}\euD_{m'_1}
:= \inf\Big\{s\ge \tau\euD_{m_1}\pv
M\eiD_1(s) \le m'_1\Big\}$.
Thanks to \cite[Lemma 3.2]{AP13}
and the strong Markov property at time $\tau\euD_{m_1}$, we can choose $m'_1>0$ such that
$\PR_x(\cE\eiD_2) 
\le \eps$ holds for any $x\in \cX_d$.

We then consider the following event as a function of $y_0\in (0, 1)$,
that describes a failure in having $X\eiD_0$ bounded away from 0:
\begin{equation*}
	\cE\eiD_3:=\Big\{\tau\euD_{m_1}\le 2t\Big\}
	\cap\Big\{\wtd{T}\euD_{m'_1}  > 5t\Big\}
	\cap\Big\{\wtd T\ekD 0_{y_0} \le 5t\Big\}\cap
	\Big\{6\, t < \ext\eiD\Big\},
\end{equation*}
where $\wtd T\ekD 0_{y_0}
:= \inf\Big\{s\ge \tau\euD_{m_1}; X\eiD_0(s) < y_0\Big\}$.
Thank to Lemma \ref{M_LX0i} and the strong Markov property at time 
$\wtd T\ekD 0_{y_0}$ 
(within the time-interval $[\tau\euD_{m_1}, \wtd{T}\euD_{m'_1}]$), 
we can choose $y_0>0$ such that
$\PR_x(\cE\eiD_3) \le \eps$ holds for any $x\in \cX_d$.

The following event is stated as a function of $m_3>0$
and describes a failure in the descent of the third moment:
\begin{equation*}
	\cE\eiD_4
	:=\Big\{ \tau\euD_{m_1} \le 2t \Big\}
	\cap\Big\{2t < \ext\eiD\Big\}\cap
	\Big\{\wtd{\tau}\etD_{m_3} > 4\,t\Big\},
\end{equation*}
where 
$\wtd{\tau}\etD_{m_3} 
:= \inf\Big\{s\ge \tau\euD_{m_1}+t\pv
M\eiD_3(s) \le m_3\Big\}$.
Thanks to Lemma \ref{M_Lmki} and the strong Markov property
at time $\tau\euD_{m_1}$,
we can choose $m_3>0$ 
such that
$\PR_x\lp \cE\eiD_4\rp
\le \eps$ holds for any $x\in \cX_d$
(with an implicit step for the second moment).

The failure in the containment of the third moment 
is stated in terms of the following event,
as a function of $m'_3>0$:
\begin{equation*}
	\cE\eiD_5
	:=\Big\{ \tau\euD_{m_1} \le 2t \Big\}\cap
	\Big\{\wtd{\tau}\etD_{m_3} \le 4\,t\Big\}\cap
	\Big\{\wtd{T}\etD_{m'_3} \le 5\,t\Big\}
	\cap\Big\{5\, t< \ext\eiD\Big\},
\end{equation*}
where $\wtd{T}\etD_{m'_3} := \inf\Big\{s\ge \tau\etD_{m_3}\pv
M\eiD_3(s) \le m'_3\Big\}$.
Thanks to Lemma \ref{M_moment_control}, 
we  can choose $m'_3>0$ such that
$\PR_x\lp \cE\eiD_5  \rp
\le \eps$ holds for any $x\in \cX_d$.

Now, we can define $J:= \Lfl \sfrac{m'_3}{\eta}\Rfl +1$
($\eta$ being an imposed parameter in the statement of Theorem \ref{eT_dInf}).
The failure in having the $J$-optimal subpopulation sizes 
bounded away from 0
is stated in terms of the following event,
as a function of $y\in (0, 1)$:
\begin{multline*}
	\cE\eiD_6:=	\Big\{ \tau\euD_{m_1} \le 2t \Big\}\cap
	\Big\{\wtd{\tau}\etD_{m_3} \le 4\,t\Big\}
	\cap\Big\{\wtd{T}\euD_{m'_1}  > 5t\Big\}
	\cap\Big\{\wtd T\ekD 0_{y_0} > 5t\Big\}
\\ \cap\Big\{ 5t < \ext\eiD \Big\}
	\cap\Big\{\wtd{\tau}\ekD J_y > 5t\Big\},
\end{multline*}
where  $\wtd{\tau}\ekD J_y
:= \inf\Big\{s\ge \wtd\tau\etD_{m_3}\pv
\frl{j\le J} X\eiD_j(s) \ge y\Big\}$.
Thanks to Lemma \ref{M_Lymi},
we  can choose $y\in (0, 1)$ such that
$\PR_x\lp \cE\eiD_6  \rp
\le \eps$ holds for any $x\in \cX_d$.

Let $E\eiD$ take the following form, which agrees with \eqref{M_cDE}:
\begin{equation}
	E\eiD 
	:=\Big\{x\in \cX_d\pv M_3(x) \le m'_3\mVg
	\frl{j\le \Lfl \sfrac{m'_3}{\eta} \Rfl+1} x_j \ge y\Big\}.
\end{equation}
Firstly, $\tau\euD_{m_1} \le 2\, t$ holds a.s. on the event $\{6t<\ext\eiD\}\setminus \cE\eiD_1$.
$\wtd{T}\euD_{m'_1}  \ge 5t$ thus holds a.s. on the event
$\{6t<\ext\eiD\}\setminus \cup_{i = 1}^2 \cE\eiD_i$.
Consequently, $\wtd T\ekD 0_{y_0} \ge 5t$
holds a.s.  on the event $\{6t<\ext\eiD\}\setminus \cup_{i = 1}^3 \cE\eiD_i$.
On the other hand, 
$\wtd \tau\etD_{m_3} \le 4\,t$
holds a.s.  on the event $\{6t<\ext\eiD\}\setminus \cup_{i = 1}^4 \cE\eiD_i$.
$\wtd{T}\etD_{m'_3} \ge 5\,t$
then holds a.s. on the event
$\{6t<\ext\eiD\}\setminus \cup_{i = 1}^5 \cE\eiD_i$.
Finally, $\wtd{\tau}\ekD J_y \le 5t$
holds a.s. on the event $\{6t<\ext\eiD\}\setminus \cup_{i = 1}^6 \cE\eiD_i$.
By definition of $\wtd{\tau}\ekD J_y$ and $\wtd{T}\etD_{m'_3}$,
and since $\wtd{\tau}\ekD J_y \in [\wtd \tau\etD_{m_3}, \wtd{T}\etD_{m'_3}]$,
$X\eD(\wtd{\tau}\ekD J_y)$
belongs to $E\eiD$ a.s. on the event $\{6t<\ext\eiD\}\setminus \cup_{i = 1}^6 \cE\eiD_i$.
This concludes the following inclusion $\{6t<\ext\eiD\wedge \tau\eiD_E\} \subset 
\cup_{i = 1}^6 \cE\eiD_i$,
which entails the following upper-bound in probability for any $d\in \II{1, \infty}$
and any $x\in \cX_d$: 
\begin{equation*}
\PR_x\Big(6\, t<\ext\eiD\wedge \tau\eiD_E\Big) \le 8\, \eps.
\end{equation*}
This concludes the proof of Theorem~\ref{eT_dInf}
(by adjusting the initial choices of $t$ and $\eps$).
\epf


\subsection{Aggregation of the last coordinates}
\label{M_sec_Agg}

The changes in the description of the system
specified in this  subsection
will be crucial for the proofs of
both the estimates of mixing (in Subsection \ref{sec_mixAcc}) 
and of almost perfect harvest
(in Subsection \ref{sec_AFdI}).
Up to a multiplicative constant in the probabilities,
they make it possible to gather the last coordinates 
in one specific block while keeping a Markovian description.
Our aim is then to prove that the dependency
in the initial values of these last coordinates
vanishes very quickly.

More precisely, 
the current subsection is dedicated to the study 
of the law $\PR\ejD_x$,
for any $d\in \II{1, \infty}$, any integer $J \in \II{1, d}$
and any $x\in \cX_d$,
 of the solution to the following set of equations,
stated for any $i\in \II{0, d}$:
\begin{equation}
\begin{split}
	&\RMd X\eD_i(t) = \alpha\cdot (M\eijD_1(t) - i\wedge J )\cdot X\eD_i(t)\, \RMd t 
+ \lambda \cdot(X\eD_{i-1}(t)-\idc{i<d}\,X\eD_{i}(t))\, \RMd t 
\\&\hcm{2}
+ \sqrt{X\eD_i(t)} \, \RMd W_i(t)
- X\eD_i(t) \, \RMd W\sD(t),
\qquad X\eD_i(0) = x,
\end{split}
\label{Skd}
\end{equation}
with $(W_i)_{i\ge 0}$ is still a family of mutually independent Brownian motions,
$W\sD$ is expressed as in \eqref{Sd}:
\begin{equation}
\RMd	W\sD(t):= \textstyle \sum_{j\in \II{0, d}}  \sqrt{X\eD_j(t)} \RMd W_j(t),
\quad W\sD(0) = 0,
\label{eq_def_WsD2}
\end{equation}
and the process $M\eijD_1$ is the first moment saturated at value $J$:
\begin{equation}	
M\eijD_1:= \textstyle  \sum_{i\in \II{0, d}} (i\wedge J)\, X\eD_i
 = \sum_{i\le J-1} i\, X\eD_i + J\,  \sum_{i\ge J} X\eD_i.
 \label{eq_def_MJD}
\end{equation}
The reason we do not include the obvious dependency in $J$ 
in the solution $X\eD$ to the system \eqref{Skd}
is that we want to connect this solution to the one to the system
\eqref{Sd} under  $\PR_x$
with a change of probability density given by the Girsanov transform.

The main interest of this law lies in that it can be efficiently projected on a finite-dimensional system, as stated and proved in Subsection~\ref{sec_Kproj}
(see Proposition~\ref{M_Kcomp}).
In the next Subsection~\ref{sec_Girs},
we will obtain comparison estimates relating $\PR\ejD$ and $\PR$
thanks to the Girsanov transform
(see Proposition~\ref{M_Girs}).
Finally, in Subsection~\ref{M_sec_momCk},
we handle the probability of large increase of the moments
(see Lemma~\ref{M_moment_controlk}).

\subsubsection{Connexion between $\PR\ejD$ and $\PR$}
\label{sec_Kproj}
Since the selection effect is identical on all the classes 
larger than $J$ under  $\PR\ejD$,
this law
is naturally associated with 
the following projection $\pi_J$ from $\cX_d$ to $\cX_J$:
\begin{equation}
	\textstyle 
	\pi_J(x)_i
	=\left\{
	\begin{aligned}
		& \quad x_i, &\text{ if } i\le J-1, 
		\\&\sum\nolimits_{j = J}^d x_j = 1-  \sum\nolimits_{j = 0}^{J-1} x_j, &\text{ if } i = J. 
	\end{aligned}
	\right.
	\label{M_piK}
\end{equation}

The renormalized sequence of the tail classes
under the projection
will be described in terms of the law
of the solution $X\eF$ to the following set of equations,
where $F$ stands for ``Final".
We first denote its corresponding total size as 
the process $X\eD_{(J)}$, with initial condition $x_{(J)}:= [\pi_J(x)]_J$:
\begin{equation}
	X\eD_{(J)}:= 1- \Tsum{i\le J-1} X\eD_i.
\label{eq_def_XJ}
\end{equation}
This process affects both the entrance flux on the class $J$,
the associated correction term due to the renormalisation
and the level of demographic fluctuations.
For $i\in \II{J, d}$,  for any $t\ge 0$:
\begin{equation}
	\label{M_XL}
\RMd X\eF_i(t)
= V\eF_i(t)\RMd t + \RMd \cN\eF_i(t),
\quad X\eF_i(0) = \dfrac{x_i}{x_{(J)}},
\end{equation}
with the following definitions of the process $V\eF_i$:
\begin{equation}
	V\eF_i
	= \lambda \cdot \lc \frac{X\eD_{J-1}}{X\eD_{(J)}}\cdot (\idc{i = J}  - X\eF_i)
	+  \idc{i \ge  J+1}X\eF_{i-1}- \idc{i<d}\,X\eF_i\rc,
	\label{eq_def_VFi}
\end{equation}
and of the martingale $\cN\eD_{(J)}$:
\begin{equation}
	\RMd\cN\eF_i(t)
	= \sqrt{\frac{X\eF_i(t)}{X\eD_{(J)}(t)}} \RMd W\eF_i(t) 
	- \frac{X\eF_i(t)}{\sqrt{X\eD_{(J)}(t)}} \RMd W\eF_{(d)}(t),
	\quad \cN\eF_i(0) = 0,
	\label{eq_def_NFi}
\end{equation}
in terms of the sequence $(W\eF_i)_{i\in \II{J, d}}$,
which defines a mutually independent family of  Brownian motions 
that are mutually independent of the family $(W_i)_{i\in \II{0, d}}$
and in terms of the martingale $W\eF_{(d)}$:
\begin{equation*}
\RMd	W\eF_{(d)}(t)
	:= \sum\nolimits_{i= J}^{d} \sqrt{X\eF_i(t)} \RMd W\eF_i(t),
	\quad W\eF_{(d)}(0) = 0.
\end{equation*}
We can then define the process $\bar{X}\eD_i$ 
as $X\eD_i$ for any $i\in \II{0, J-1}$
and as $X\eD_{(J)} \cdot \bar{X}\eF_i$ for any $i\in \II{J, d}$.
\begin{prop}
	\label{M_Kcomp}
	For any $J\ge 1$, 
	$\pi_J(X\eD)$ is by itself a Markov process
	under any $\PR\ejD_x$, 
	whose law is independent of  $d\in \II{J,\, \infty}$
	and  depends on $x\in \cX_d$ only through $\pi_J(x)$.
	Under $\PR\ejD_x$,
	the process $\bar X\eD$  on $\cX_d$
	has the same law as the process $X\eD$.
\end{prop}


\bpf: 
By virtue of \eqref{eq_def_MJD}
and of \eqref{eq_def_XJ}:
$$M\eijD_1 
:= \sum\nolimits_{i= 0}^{J-1} i\, X\eD_i + J\, X\eD_{(J)}.$$
Under $\PR\ejD_x$,
for any $x\in \cX_d$,
Itô's lemma then entails that the process $X\eD_{(J)}$
is solution to the following SDE,
for any $t\ge 0$:
\begin{equation}
\RMd X\eD_{(J)}(t) 
= V\eD_{(J)}(t) \RMd t 
+ \RMd\cN\eD_{(J)}(t)
- X\eD_{(J)}(t) \, \RMd W\sD(t),
\label{eq_def_XJD}
\end{equation}
with the definition of $W\sD$ from \eqref{eq_def_WsD2}
and the following definitions of  the process $V\eD_{(J)}$:
\begin{equation}
	V\eD_{(J)}
	= \alpha\cdot (M\eijD_1 - J)\cdot X\eD_{(J)}
	+\lambda\cdot X\eD_{J-1},
\label{eq_def_VJD}
\end{equation}
and of the martingale $\cN\eD_{(J)}$:
\begin{equation}
\RMd \cN\eD_{(J)}(t)
= \sum\nolimits_{j= J}^{d} \sqrt{X\eD_j(t)}\, \RMd W_j(t),
\quad \wtd \cN\eD_{(J)}(0) = 0.
\label{eq_def_NJD}
\end{equation}
Since the sequence $(W_i)_{i\in \II{0, d}}$ 
defines a mutually independent family of Brownian motions,
the following identity holds for any $t\ge 0$:
\begin{equation*}
	\RMd\LAg \cN\eD_{(J)}\RAg_t 
	= \sum\nolimits_{j= J}^{d} X\eD_j(t) \RMd t
= X\eD_{(J)}(t) \RMd t.
\end{equation*}
Thanks to \cite[Theorem 18.12]{Ka02},
we can define a Brownian motion $W\eD_{(J)}$
such that the following SDE hold:
$\RMd \cN\eD_{(J)}(t)
= \sqrt{X\eD_{(J)}(t)} \RMd W\eD_{(J)}(t).$
Recalling \eqref{eq_def_WsD2},
it yields the following alternative identity
for the Brownian motion $W\sD$:
\begin{equation}
		\RMd W\sD(t)= \sum\nolimits_{i=0}^{J-1}  \sqrt{X\eD_i(t)} \RMd W_i(t)
	+ \sqrt{X\eD_{(J)}(t)} \RMd W\eD_{(J)}(t).
\label{eq_def2_WsD}
\end{equation}
The correlation between  $W\eD_{(J)}$ and the $(W_i)_{i\le J-1}$
remains zero,
while they constitute a system of Brownian motions 
under the same filtration.
$W\eD_{(J)}$
is thus independent of $\sigma(W_i; i\le J-1)$,
so that the system of equations satisfied by $\pi_J(X)$ 
is equivalent for any $\PR\ejD_x$.

Thanks to Itô's lemma, the following identity holds for any $t\ge 0$
and any $i\in \II{J, d}$:
\begin{equation}
	d\bar{X}\eD_i(t)
	= \bar V\eD_i(t) \RMd t + \RMd\bar \cN\eD_i(t)
	+ \tfrac{1}{2}\,  \RMd \LAg \cN\eD_{(J)}, \cN\eF_i\RAg_t,
	\label{eqXbI}
\end{equation}
with $\cN\eD_{(J)}(t)$ and $\cN\eF_i(t)$
defined respectively in \eqref{eq_def_XJD} and \eqref{eq_def_NFi},
while the martingale component $\bar \cN\eD_i$
is expressed as follows
after simplifications:
\begin{equation*}
\begin{split}
\RMd\bar \cN\eD_i(t)
&= X\eD_{(J)}(t) \RMd\bar \cN\eF_i(t)
+ \bar X\eF_i(t) \RMd \cN\eD_{(J)}(t)
\\ &= \sqrt{\bar{X}\eD_i(t)} \RMd W_i(t)
- \bar{X}\eD_i(t) \RMd W\sD(t),
\quad \bar \cN\eD_i(0) = 0,
\end{split}
\end{equation*}
and the process $\bar V\eF_i(t)$ as follows,
thanks to \eqref{eq_def_VJD} and to \eqref{eq_def_VFi}
after simplifications:
\begin{equation*}
	\begin{split}
	\bar V\eF_i
&	= X\eD_{(J)} \cdot \bar V\eF_i
	+ \bar X\eF_i \cdot V\eD_{(J)}
		\\&= \alpha \cdot(M_1\eijD - J)\cdot \bar{X}\eD_i 
		+ \lambda  \cdot(\bar{X}\eD_{i-1} - \idc{i<d}\bar{X}\eD_i).
	\end{split}
\end{equation*}
$\cN\eF_i$ is stated in terms of the sequence $(W\eF_i)_{i\in \II{0, d}}$
of Brownian motions, 
which is independent of the sequence $(W_i)_{i\in \II{0, d}}$
in terms of which the martingale $\cN\eD_{(J)}$ is stated.
Therefore, $\RMd \LAg \cN\eD_{(J)}, \cN\eF_i\RAg \equiv 0$.
The system of SDEs \eqref{eqXbI} satisfied by $(\bar{X}\eD_i)_{i\in \II{0, d}}$
thus coincides with the system \eqref{Sd} satisfied by $(X\eD_i)_{i\in \II{0, d}}$.
By the uniqueness of the whole system,
cf Proposition \ref{prop_Trep}
in the appendix,
$X\eD$ coincide with $\bar{X}\eD$.
%
This ends the proof of Proposition \ref{M_Kcomp}.
\epf

\subsubsection{The connexion formula between $\PR\ejD$ and $\PR$}
\label{sec_Girs}	

Thanks to the Girsanov transform,
we will establish a relevant quantification 
for the transfer between the original law of the process $\PR$
and the law of $\PR\ejD$.  
For the upcoming Proposition~\ref{M_Girs}, 
we shall exploit a control on moments of order $k$.
We recall the definition of $T\eDk_{m}$
from \eqref{eq_def_TkM} for any $m>0$:
\begin{equation*}
T\eDk_{m}:= \inf\big\{s \ge 0\pv M\eD_{k}(s) \ge m\big\}.
\end{equation*}

\begin{prop}
\label{M_Girs}
Given any $t, \eps>0$, $k \ge 2$, there exists $C_M, C_G>0$ 
for which the following holds.
For any $m\ge 1$, with $m':= C_M\cdot m$,
for any $d\in \II{1, \infty}$, any $J\le d$, 
and  any $x\in \cX_{d}\cap \cX^{2 k}$ such that $M_{k}(x) \le m$,
there exists a coupling between $\PR\ejD$ and $\PR$
 such that the following upper-bound holds
 a.s. on the event $\{t < T\eDk_{m'}\}$:
 \begin{equation*}
 \quad \Bv \log \lp \frac{\RMd \PR\ejD_x}{\RMd \PR_x}\! 
 \left\vert \underset{[0, t]}{\quad}\right.\rp \Bv
 \le C_G\frac{m}{J^{k - 2}},
 \end{equation*}
while the event $\{t < T\eDk_{m'}\}$
happens with a probability close to 1 in the following sense:
\begin{equation*}
\PR_x\lp  T\eDk_{m'} \le t \rp
\le \eps.
\end{equation*}
\end{prop}
\begin{rem}
In this article, we exploit Proposition \ref{M_Girs}
only for $k = 3$. 
The proof is extended to any real number $k\ge2$
to explicit our motivation for choosing $k>2$
(see \eqref{ineq_t3m} below).
\end{rem}

\paragraph{The explicit connexion formula} 
\hfill\\
We express in this subsection the Girsanov transform 
that makes it possible to relate $\PR\ejD$ and $\PR$.
It is expressed in the upcoming Lemma \ref{M_dPkd}
in terms of the processes $R_1\eijD$
and $R_2\eijD$ with the following definitions: 
\begin{equation*}
	R_1\eijD:= \sum\nolimits_{i= J+1}^d (i-J)\cdot X\eD_i\mVg\quad
	R_2\eijD:= \sum\nolimits_{i= J+1}^d  (i-J)^2\cdot X\eD_i.
\end{equation*}
One can notice that they correspond to the expectation and variance
of the vector $(Y_i)_{i\in \Z_+}$
such that $Y_0 =  \sum_{j = 0}^{J} X\eD_j$
and for any $i\in \N$, $Y_i = X\eD_{J+i}$.

\begin{lem}
	\label{M_dPkd}
For any $d\in \II{1, \infty}$ and $J\in \II{1, d}$, there exists a coupling between $\PR\ejD$ and $\PR$
	such that:
	\begin{equation*}
\begin{split}
	\log \frac{\RMd \PR\ejD_x}{\RMd \PR_x}\! 
	 \left\vert \underset{[0, t]}{\quad}\right.
	= \alpha\cdot R_1\eijD(0) - \alpha\cdot R_1\eijD(t)
		+ \int_0^t G\eijD(s) \RMd s,
\end{split}
	\end{equation*}
with the following definition of the process $G\eijD$:
\begin{equation*}
\begin{split}
	G\eijD
&:= \alpha^2\,  (M\eD_1 - J)\, R_1\eijD
-  \alpha^2\, R_2\eijD
+ \alpha\, \lambda\, (X\eD_{(J)}
-\idc{d<\infty} X\eD_d)
\\&\quad - \frac{\alpha^2}{2} \, \lc R_2\eijD
- (R_1\eijD)^2\rc.
\end{split}
\end{equation*}
\end{lem}

\bpf: We define as follows the  martingale $\cL\eijD$,
starting at 0:
\begin{equation*}
\begin{split}
	\RMd	\cL\eijD(t) 
&:= - \alpha\sum_{i\ge J+1} (i-J) \sqrt{X\eD_i(t)}\, \RMd W_i(t) 
	+ \alpha \cdot R_1\eijD(t)\, \RMd W\sD(t)
\\&	=- \alpha\sum_{i\ge J+1} (i-J)\cdot \lc 
	\sqrt{X\eD_i(t)}\, \RMd W_i(t) -
	X\eD_i(t)\, \RMd W\sD(t)\rc.
\end{split}
\end{equation*}
By this choice,
we obtain the following identities,
for any $i\in \II{0, d}$:
\begin{equation}
\begin{split}
	&\RMd \LAg \cL\eijD, W_i\RAg_s
= \alpha \cdot
\lc R\eD_1(s) - (i-J)_+ \rc \cdot\sqrt{X\eD_i(s)}\, \RMd s,
\\&
\RMd \LAg \cL\eijD, W\sD\RAg_s
= 0.
\end{split}
\label{eq_LijD_id}
\end{equation}
Recalling the systems of SDEs \eqref{Sd}
and \eqref{Skd},
it entails that the Girsanov transform 
of the law $\PR$ with respect to the exponential martingale 
of $\cL\eijD$
generates $\PR\ejD$,
in the sense that
the following property holds for any $t>0$
and $x\in \cX_d$:
\begin{equation}
	\log \frac{\RMd \PR\ejD_x}{\RMd \PR_x}\! 
	\left\vert \underset{[0, t]}{\quad}\right.
	= \cL\eijD(t) 
	 -\frac{1}{2}\,   \RMd \LAg \cL\eijD\RAg_t.
	\label{eq_dec_dPdP}
\end{equation}
Thanks to \eqref{eq_LijD_id},
the quadratic variation $\RMd \LAg \cL\eijD\RAg_t$ 
satisfies the following identity:
\begin{equation}
\begin{split}
		\RMd \LAg \cL\eijD\RAg_t
	&= -\alpha\sum_{i\ge J+1}\, (i-J)\, \sqrt{X\eD_i(t)}
	\RMd \LAg \cL\eijD, W_i\RAg_t
\\	&=  \alpha^2\cdot \lc R_2\eijD(t)
	- R_1\eijD(t)^2\rc.
\end{split}
\label{eq_def_VarLijD}
\end{equation}
On the other hand,
we note the following identity
for any $s\ge 0$:
\begin{equation*}
\RMd R_1\eijD(s)
= V_1\eijD(s) \RMd s
- \frac{1}{\alpha} \RMd\cL\eijD(s),
\quad R_1\eijD(0) = \sum_{i\ge J+1} (i-J)\cdot x_i,
\end{equation*}
where the process $V_1\eijD$ is defined as follows:
\begin{equation*}
V_1\eijD
= \alpha\cdot\lc (M\eD_1 - J)\cdot R_1\eijD - R_2\eijD\rc \RMd s
+ \lambda\cdot (X\eD_{(J)}
-\idc{d<\infty} X\eD_d).
\end{equation*}
With this alternative expression for $\cL\eijD$,
recalling \eqref{eq_dec_dPdP} and
\eqref{eq_def_VarLijD},
we conclude the proof of Lemma \ref{M_dPkd}.
\hfill \epf

\paragraph{Concluding the proof of Proposition~\ref{M_Girs}}
\hfill\\
 The aim is now to get uniform upper-bound
on the expression given in Lemma \ref{M_dPkd}.
We consider any $k\ge 2$.
Firstly, the following inequalities
hold by definitions of  $R_1\eijD\ge 0$ and $M\eD_{k}$:
\begin{equation*}
R_1\eijD \le J^{-(k-1)} \Tsum{i\ge J+1} i^{k}\, X\eD_i
	\le J^{-(k-1)} M\eD_{k}.
	\label{M_R1}
\end{equation*}
Similarly for $R_2\eijD\ge 0$ and $X\eD_{(J)} - \idc{d<\infty} X\eD_d(s) \ge 0$:
\begin{equation*}
	R_2\eijD \le J^{-(k-2)} M\eD_{k},
	\qquad 	 X\eD_{(J)} - \idc{d<\infty} X\eD_d(s)\le J^{-k} M\eD_{k}.
\end{equation*}
Thanks to the Cauchy-Schwarz inequality,
$(R_1\eijD)^2 \le R_2\eijD$.
Thanks to Hölder's inequality,
both $M\eD_1 \le (M\eD_{k})^{1/k}$ and $M\eD_{k-1} \le (M\eD_{k})^{(k-1)/k}$
hold,
which finally yields the following inequalities
about $M\eD_1\cdot R_1\eijD \ge 0$:
\begin{equation*}
M\eD_1\cdot R_1\eijD 
	\le J^{-(k-2)} M\eD_1\cdot M\eD_{k-1}
	\le J^{-(k-2)} M\eD_{k},
	\label{M_MR1}
\end{equation*}
Thanks to Lemma \ref{M_dPkd}, 
we can thus choose a constant $C_1>0$
such that the following inequality
holds a.s. on the event $\{t< T\eDk_{m'}\}$
for any $m, m'\ge 1$ such that $m<m'$, any $d\in \II{1, \infty}$, any $J\in \II{1, d}$, and any $x\in \cX_d$ 
such that $M_k(x) \le m$:
\begin{equation}
	\Bv\log\Big( \frac{\RMd \PR\ejD_x}{\RMd \PR_x}\! 
	\left\vert \underset{[0, t]}{\quad}\right.\Big)\Bv
	\le C_1\, \frac{m'}{J^{k - 2}}.
	\label{M_dPdP}
\end{equation}
 
Thanks to Lemma \ref{M_moment_control},
we can choose $C_2>0$ 
such that
	the following inequality holds
for any $m, m'\ge 1$ such that $m< m'$, any $d\in \II{1, \infty}$ 
and any $x\in \cX_d$ such that $M_k(x)\le m$:
\begin{equation*}
\PR_x\lp  T\eDk_{m'}\le t \rp
\le \dfrac{C_2 m}{m'}.
\end{equation*}

Let $\eps>0$. We thus define $m' = C_M\cdot m$,
where $C_M:=  C_2/\eps$,
so that the above upper-bound is exactly $\eps$.
Thanks to \eqref{M_dPdP}
with $C_G = C_1\cdot C_2/\eps$, 
the following inequality holds a.s.
on the event  $\{t< T\eDk_{m'}\}$,
for any $m\ge 1$, any $d\in \II{1, \infty}$, any $J\in \II{1, d}$, and any $x\in \cX_d$ 
such that $M_k(x) \le m$:
\begin{equation*}
	\Bv\log \Big(\frac{\RMd \PR\ejD_x}{\RMd \PR_x}\! 
	\left\vert \underset{[0, t]}{\quad}\right.\Big)\Bv
	\le C_G\, \frac{m}{J^{k - 2}}.
\end{equation*}
This concludes the proof of Proposition \ref{M_Girs}.
\epf

\subsubsection{Upper-bound on the probability of moment increase
	for $X\eF$}
\label{M_sec_momCk}

Similarly as for Lemma~\ref{M_moment_control}, 
exploiting the decomposition in Proposition~\ref{M_Kcomp},
we define the third moment $M_3\eF$ of $X\eF$:
\begin{equation*}
M_3\eF:= \sum\nolimits_{i=J}^d i^3\, X\eF_i
\in [J^3, \infty),
\end{equation*}
and the corresponding hitting time $T\etF_{m}$ of the value $m>0$:
\begin{equation}
T\etF_{m}
:= \inf\Big\{s \ge 0\pv M_3\eF(s) \ge m\Big\}.
\label{M_tauM3k}
\end{equation}
The upper-bound is to be obtained up to the following hitting time $\tau_0\ejD$:
\begin{equation}
\tau_0\ejD:=\inf\Big\{t\ge 0 \pv X\eD_{(J)}(t) =0\Big\}.
\label{eq_def_tau0ejD}
\end{equation}

For clarity, we define
$\cF^{(J)} = \sigma\Big(W_i: i\le J-1\pv W\eD_{(J)}\Big)$.
Recall that the process $X\eF_i$
is driven by Brownian motions $(W\eF_i: i\ge J)$
that are independent of $\cF^{(J)}$.
The inclusion $\sigma(\pi_J(X)) \subset \cF^{(J)}$ is directly obtained
through the autonomous set of equation verified by $\pi_J(X)$.
The following control on $M_3\eF$ exploits the filtration $\cF^{(J)}_t:= \cF^{(J)}\vee \cF_t$.

\begin{lem}
	\label{M_moment_controlk}
	For any $t>0$, 
	there exists $C\ge 1$
	such that the following inequality holds 
	a.s. 
	for any $m, m'>0$, any $d\in \II{1, \infty}$ 
	and any $x\in \cX_d$ 
	such that $M_3\eF(x)\le m$:
	\begin{equation*}
\PR_x\eijD\lp  T\etF_{m'}\le t\wedge \tau_0\ejD
	\bv \cF^{(J)}\rp
	\le \dfrac{C m}{m'}.	
	\end{equation*}
\end{lem}

\bpf: Under $\PR\ejD$, 
we exploit the Itô formula 
to express 
$M_3\eF$ as the solution to the following SDE:
\begin{equation}\label{M_VMk}
\RMd M_3\eF(t) 
:= V\eF_3(t) \RMd t + \RMd \M\eF_3(t),
\end{equation}
where $V\eF_3$ is a bounded variation process
defined as:
\begin{equation*}
V\eF_3
:= \lambda\cdot \lc
\frac{X\eD_{J-1}}{X\eD_{(J)}}\cdot (J^3 - M_3\eF)
+ \sum_{\ell\ge J} (\ell + 1)^{3}\cdot X\eF_{\ell} 
- M\eF_3  + \idc{d<\infty} d^3\cdot X\eF_d\rc.
\label{M_Vd}
\end{equation*}
Note that whatever the values of $(X\eD_{J-1}/X\eD_{(J)})$,
with the rough estimate $(\ell + 1)^3 \le 8 \ell^3$ for $\ell \ge 1$,
the inequality $V\eF_3 \le 8\lambda\, M\eF_3$ holds. 
On the other hand,
the local martingale process $\M\eF_3$ is 
expressed as follows in terms of the martingales
$\cN\eF_i$ defined in \eqref{eq_def_NFi}:
\begin{equation}\label{M_Mk3}
\RMd \M\eF_3(t)
:= \sum_{i\ge J} i^3 \cdot\cN\eF_i(t),
\quad \M\eF_3(0)= 0.
\end{equation}
Relying on the same calculations as for $M_3$,
$\M\eF_3$ is a continuous local martingale starting from 0
for the filtration $\cF^{(J)}_t$
whose quadratic variation 
satisfies the following identity:
$$\RMd \langle\M\eF_3\rangle_t
= \dfrac{M\eF_6(t)- (M\eF_3(t))^2}{X\eD_{(J)}(t)},$$
where $M\eF_6(s) := \sum_{i\ge J} i^6\, X\eF_i(t)$.
The definition of the localization time can be adapted from \eqref{def_Tell}
as follows for any positive integer $\ell$:
\begin{equation}
	T_\ell:= \inf\{ s \ge 0\pv \LAg\M_3\eF\RAg_s \ge \ell\mVg 
M_3\eF(s) \ge \ell\}.
\label{def_Tell_2}
\end{equation}
Thanks to Proposition \ref{M_Kcomp}, the following inequalities 
	 hold for any $t$
such that $X\eD_{(J)}(t)>0$:
\begin{equation*}
M_3\eF(t)\le \frac{M_3(t)}{X\eD_{(J)}(t)},
\qquad
M_6\eF(t)\le \frac{M_6(t)}{X\eD_{(J)}(t)}.
\end{equation*}
Thanks to the proof of \cite[Theorem 3]{AP13}
(see also Propositon \ref{pr_Fin_mom} in the appendix), 
$M\eI_6$ is  locally upper-bounded, a.s. under $\PR\ejD_x$
for any $x\in \cX_\infty = \cX^6$.
It entails that both $M_3\eF$
	and $\LAg\M_3\eF\RAg$
	 are also a.s. upper-bounded
	 for any $n\ge 1$ on the time-interval $[0, t\wedge \tau_{\sfrac{1}{n}}\ejD]$,
	 where $\tau_{\sfrac{1}{n}}\ejD 
	 = \inf\{t\ge 0 \pv X\eD_{(J)}(t) = 1/n\}$. 
Taking the limit with $n$ tending to infinity, 
$\lim_\ell T_\ell \ge t\wedge \tau_{0}\ejD$ holds a.s.
The rest of the proof of Lemma \ref{M_moment_controlk}
can be taken mutatis mutandis from the one of Lemma \ref{M_moment_control}
(with $C = \exp[8\lambda t]\vee 1$).
\epf

\subsection{Mixing property and accessibility}
\label{sec_mixAcc}

Theorem~\ref{M_Mix.dInf}, stated and proved in Subsection~\ref{sec_MixdInf}, 
is the main result of the current Subsection~\ref{sec_mixAcc}.
It establishes the mixing estimate $(A2)$.

We consider three intermediate steps:
 first in Subsection~\ref{sec_Hdl},
 we justify that an interior subset of $\cX_d$
 with convenient properties can be accessed,
 cf Lemma \ref{lem_Hdl};
secondly in Subsection~\ref{sec_X0diff},
we show a mixing estimate on the two optimal subpopulation sizes $X\eD_0$ and $X\eD_1$,
cf Lemma \ref{lem_Xdiff};
thirdly in Subsection~\ref{sec_Xcut},
we prove the existence of a uniform lower-bound
on the probability of the event on which we will condition
to produce $\zeta\eD$ as a probability measure, 
cf Lemma \ref{sec_Xcut}.

We recall the definition of $\cD\eD_\ell$ from \eqref{def_DL}
for any integer $\ell$:
\begin{equation*}
\cD\eD_\ell:= \Big\{x\in \cX_d\pv  
x_0 \ge  \frac{1}{2\ell}\Big\},
\end{equation*}
and the fact that $T\eD_{\cD_\ell}$ denotes the exit time of $X\eD$ out of $\cD\eD_\ell$.

\begin{rem}
	It can be noted that for any $x\in \cD\eD_\ell$, $T\eD_{\cD_\ell}< \ext\eD$, 
	so that $T\eD_{\cD_\ell} = \ext\eD\wedge T\eD_{\cD_\ell}$.
\end{rem}

\begin{rem}
In the current Subsection~\ref{sec_mixAcc},
we apply the decomposition
introduced in Subsection~\ref{M_sec_Agg}
for $J = 2$.
The definitions \eqref{M_piK} for $\pi_J$ and \eqref{eq_def_XJ} for $X\eD_{(J)}$ will be used 
below in case $J= 2$.
The proofs given in the next Subsection~\ref{sec_AFdI}
will exploit a generalization of the argument for large values of $J$.
\end{rem}

\subsubsection{Access to an interior point}
\label{sec_Hdl}
Subsection \ref{sec_Hdl} is devoted to the proof of the upcoming Lemma~\ref{lem_Hdl},
in which we justify the accessibility of $\cH\eD_\ell$, 
which is an interior subset of $\cX_d$ with convenient properties.
These properties entail that any state $x\in \cH\eD_\ell$
will constitute a suitable initial condition in order to exploit Property~$(H)$. 

\begin{lem}
	\label{lem_Hdl}
	For any integer $\ell\ge 1$,
	there exist four real numbers $\tp, \mpr, y>0$, such that 
	the following inequality holds for any
	$d\in \II{2, \infty}$ and any
	$x\in \cD\eiD_\ell$:
	\begin{align*}
		\PR_x\Big(M\eD_3(\tp) \le m\mVg 
		X\eiD_0(\tp)\wedge X\eiD_1(\tp)\wedge X\eiD_{(2)}(\tp)\ge y
		\mVg \tp < T\eD_{\cD_{2\ell}}\Big)
		\ge \cp.
	\end{align*}
\end{lem}
This lemma leads to the introduction of the following subset of $\cX_d$
for any $d\in \II{2, \infty}$ and any $\ell\in \N$:
\begin{equation*}
\cH\eD_\ell = \cH\eD(m_\ell, y_\ell) :=\Big\{x'\in \cX_d\pv M_3(x')\le m_\ell\mVg 
x'_0\wedge x'_1\wedge \Tsum{k\ge 2} x'_k \ge y_\ell\},
\end{equation*}
where $m_\ell, y_\ell>0$ are the values of $m$ and $y$ 
associated to $\ell$ through Lemma~\ref{lem_Hdl},
so that the lower-bounded probability can be expressed as $\PR_x\big(X\eD(\tp)\in \cH\eD_\ell
\pv \tp < T\eD_{\cD_{2\ell}}\big)$.  
Without restriction, $y_\ell\le \sfrac{1\!}{4\ell}$ can be assumed.
\\

\bpf: Let $\ell\ge 1$.
We define $y_\wedge = \sfrac{1\!}{4\ell}$.
The process $X\eiD_0$, which is the initial component of the solution to \eqref{Sd}-\eqref{eq_def_Bi},
is lower-bounded under $\PR_x$ for any $d$ and $x\in \cD\eD_\ell$
by the solution $Y_0$ to the following SDE,
thanks to Lemma \ref{lem_comp}:
\begin{equation*}
	\RMd Y_0(s) = -\lambda\, \RMd s
	+ \sqrt{Y_0(s)\cdot (1-Y_0(s))}\, \RMd B_0(s)\mVg \quad Y_0(0) = 2 y_\wedge,
\end{equation*}
where $B_0$ is a Brownian motion.
 With $c_0:= \PR\big(\inf_{t\in [0, 1]} Y_0(t) > y_\wedge\bv Y_0(0) = 2 y_\wedge\big) >0$,
we thus deduce the following inequality for any $d\in \II{2, \infty}$
and any $x\in \cD\eD_\ell$:
\begin{equation}\label{M_0Svz}
	\PR_x\Big(1<T\eD_{\cD_{2\ell}}\Big) \ge c_0.
\end{equation}
Thanks to Lemmas~\ref{M_Lm1i} and \ref{M_Lmki}
(similarly as in Subsection~\ref{sec_eT_dInf}),
there exists $m_D>0$
such that the following inequality holds for any $d\in \II{2, \infty}$
and any $x\in \cD\eD_\ell$:
\begin{equation}\label{eq_min_tautDm}
	\PR_x\Big(\tau\etD_{m_D} \le \tdiv{3}\Big) \ge \sqrt{1-\frac{c_0}{2}}.
\end{equation}
Thanks to Lemma~\ref{M_moment_control},
there exists $m\ge m_D$
such that the following inequality holds for any $d\in \II{2, \infty}$
and any $x\in \cX_d$ such that $M_3(x)\le m_D$:
\begin{equation}\label{eq_min_TtDm}
g_m[x]	:=\PR_x\Big(1<T\etD_m\Big) \ge \sqrt{1-\frac{c_0}{2}}.
\end{equation}
Thanks to the strong Markov property at time $\tau\etD_m$,
recalling \eqref{M_0Svz}, \eqref{eq_min_tautDm} and \eqref{eq_min_TtDm},
the following inequalities hold for any $d\in \II{2, \infty}$
and any $x\in \cD\eD_\ell$:
\begin{equation}\label{eq_min_T0tDm}
\begin{split}
&		\PR_x\Big(\tau\etD_{m_D} \le \tdiv{3}\mVg
	1<T\etD_m\wedge T\eD_{\cD_{2\ell}}\Big) 
\\&\quad \ge \PR_x\Big(1<T\eD_{\cD_{2\ell}}\Big) - 
\Big[1-\bE_x\Big(g_m\big[X\eD(\tau\etD_{m_D})\big]\pv \tau\etD_{m_D} \le \tdiv{3}\Big) \Big]
\\&\quad	\ge \frac{c_0}{2}.
\end{split}
\end{equation}
Let $b_1 := (\lambda y_\wedge)\wedge (\tdiv{4})$.
Thanks to Lemma~\ref{lem_comp},
the process $X\eD_1$ is lower-bounded 
a.s. on the event $\{1<T\eD_{\cD_{2\ell}}\}$
under $\PR_x$
on the time-interval $[\tdiv{3}, 1]$
by the solution $Y_1$ to the following SDE:
\begin{equation}
	\RMd Y_1(s) = \varphi_1\, \RMd s
	- (\alpha + \lambda)\cdot Y_1(s)\, \RMd s
+ \sqrt{Y_1(s)\cdot (1-Y_1(s))}\, \RMd B_1(s)\mVg \quad Y_1(\tdiv{3}) = 0.
\label{eq_def_Y1}
\end{equation}
Thanks to Corollary~\ref{cor_bound_oneD},
 since $\varphi_1\in (0, \tdiv{4})$,
0 is a regular reflecting boundary for this process $Y_1$.
Therefore, there exists $y_1\in (0, \tdiv{2\lambda})$ such that 
$\PR(Y_1(\sfrac{2\!}{3})< 2 y_1)\le \sfrac{c_0\!}{4}$.
Recalling \eqref{eq_min_T0tDm}, 
the following inequality thus holds for any $d\in \II{2, \infty}$
and any $x\in \cD\eD_\ell$:
\begin{equation}\label{eq_min_X_1}
\PR_x\Big(\tau\etD_{m_D} \le \tdiv{3}\mVg
		X\eD_1(\sfrac{2\!}{3})\ge 2 y_1\mVg
		1<T\etD_m\wedge T\eD_{\cD_{2\ell}}\Big) 
\ge \frac{c_0}{4}.
\end{equation}
Thanks to Lemma~\ref{lem_comp},
the process $X\eD_1$ is lower-bounded 
a.s. on the event $\{1<T\eD_{\cD_{2\ell}}\}\cap \{X\eD_1(\sfrac{2\!}{3})\ge 2 y_1\}$
under $\PR_x$
on the time-interval $[\sfrac{2\!}{3}, 1]$
by the solution $Y_1^I$ to the following SDE,
where the difference with \eqref{eq_def_Y1} lies in the initial condition 
at time $\sfrac{2\!}{3}$:
\begin{equation*}
	\RMd Y_1^I(s) = \varphi_1\, \RMd s
	- (\alpha + \lambda)\cdot Y_1^I(s)\, \RMd s
	+ \sqrt{Y_1^I(s)\cdot (1-Y_1^I(s))}\, \RMd B_1(s)\mVg \quad Y_1^I(\sfrac{2\!}{3}) =  2 y_1.
\end{equation*}
We consider the two following stopping times $T^1_{y_1}$ and $T^{1, (d)}_{y_1}$:
\begin{equation*}
T^1_{y_1} := \inf\{s\ge \sfrac{2\!}{3}\pv Y_1^I(s) \le y_1\},
\qquad 
T^{1, (d)}_{y_1} := \inf\{s\ge \sfrac{2\!}{3}\pv X\eD_1(s) \le y_1\},
\end{equation*}
namely the hitting time of $y_1$ after time $\sfrac{2\!}{3}$
 by the processes respectively $Y_1^I$ and $X\eD_1$.
There exists $t_1\in (0, \tdiv{3})$ such that: 
\begin{equation*}
\PR\Big(T^1_{y_1} \le \sfrac{2\!}{3} + t_1\Big)\le \dfrac{c_0}{8}.
\end{equation*}
Let $t= \sfrac{2\!}{3} + t_1$.
The following inequality thus holds for any $d\in \II{2, \infty}$
and any $x\in \cD\eD_\ell$:
\begin{equation}\label{eq_min_XL1}
	\PR_x\Big(\tau\etD_{m_D} \le \tdiv{3}\mVg
	X\eD_1(\sfrac{2\!}{3})\ge 2 y_1\mVg
	t < T^{1, (d)}_{y_1} \wedge T\etD_m\wedge T\eD_{\cD_{2\ell}}\Big) 
	\ge \dfrac{c_0}{8}.
\end{equation}
Let $\varphi_2 := \lambda\, y_1$.
Thanks to Lemma~\ref{lem_comp},
the process $X\eD_2$ is lower-bounded 
a.s. on the event $\{X\eD_1(\sfrac{2\!}{3})\ge 2y_1\}\cap\{
t < T^{1, (d)}_{y_1}\}$
under $\PR_x$
on the time-interval $[\sfrac{2\!}{3}, t]$
by the solution $Y_2$ to the following SDE:
\begin{equation}
	\RMd Y_2(s) = \varphi_2\, \RMd s
	- (2 \alpha + \lambda)\cdot Y_2(s)\, \RMd s
	+ \sqrt{Y_2(s)\cdot (1-Y_2(s))}\, \RMd B_2(s)\mVg \quad Y_2(\sfrac{2\!}{3}) = 0.
	\label{eq_def_Y2}
\end{equation}
Thanks to Corrolary~\ref{cor_bound_oneD} since $\varphi_2\in (0, \tdiv{2})$,
0 is a regular reflecting boundary for this process $Y_2$.
Therefore, there exists $y_2\in (0, y_1)$ such that 
$\PR(Y_2(t)< y_2)\le c_0/16$.
Recalling \eqref{eq_min_XL1}, 
the following inequality thus holds for any $d\in \II{2, \infty}$
and any $x\in \cD\eD_\ell$:
\begin{equation}\label{eq_min_X_2}
\begin{split}
&\PR_x\Big(M\eD_3(\tp) \le m\mVg 
X\eiD_0(\tp)\wedge X\eiD_1(\tp)\wedge X\eiD_{(2)}(\tp)\ge y_2
\mVg \tp < T\eD_{\cD_{2\ell}}\Big)
\\&\ge \PR_x\Big(\tau\etD_{m_D} \le \tdiv{3}\mVg
X\eD_1(\sfrac{2\!}{3})\ge 2 y_1\mVg
X\eD_2(t)\ge y_2\mVg
t<T^{1, (d)}_{y_1} \wedge T\etD_m\wedge T\eD_{\cD_{2\ell}}\Big) 
\\&\ge \frac{c_0}{16}.
\end{split}
\end{equation}
This concludes the proof of Lemma \ref{lem_Hdl}.\epf

\subsubsection{Mixing estimate on the two optimal subpopulation sizes}
\label{sec_X0diff}

Subsection \ref{sec_X0diff} is devoted to the proof of the upcoming Lemma~\ref{lem_Xdiff},
in which we both establish a mixing property for the process $(X\eD_0, X\eD_1)$
and obtain a control on $X\eD_1$ and on $X\eD_{(2)}$.
As a reference initial condition, we consider $z\eD\in \cX_d$ to be such that $z\eD_k = 2^{-k-1}$ for any $k\in \II{0, d-1}$ 
(that is $k \in \Z_+$ for $d=\infty$)
and $z\eD_d = 2^{-d}$ (in the case $d<\infty$).
Note that $z\eD\in \cD\eD_1$.

We introduce as follows some subsets $\Y_2(y)$ of $\cX_2$
in terms of some parameter $y\in (0, y_\ell)$
that describes its gap to the boundary of $\cX_2$:
\begin{equation*}
\Y_2(y)
:= \Big\{z\in \cX_2\pv z_0\wedge z_1\wedge z_2 \ge y\Big\}.
\end{equation*}
There exist two connected open relatively compact sets $\mathfrak{K}_2^\wedge(y), \mathfrak{K}_2^\vee(y)$
with $C^\infty$-boundaries
with the following properties:
\begin{equation*}
	\Y_2(2y)
	\subset \mathfrak{K}_2^\wedge(y),
	\quad 
	\overline{\mathfrak{K}_2^\wedge(y)}
	\subset\mathfrak{K}_2^\vee(y)
	\subset \Y_2(y).
\end{equation*}
 The stopping time $T_y^{2, (d)}$ is defined for any $y\in (0, y_\ell)$ as follows:
\begin{equation}\label{eq_def_Ty2d}
T_y^{2, (d)} 
:= \inf\Big\{t\ge 0\pv \pi_2(X\eiD(\tp)) \notin \mathfrak{K}_2^\vee(y)\Big\},
\end{equation}
namely the exit time of $\pi_2(X\eiD)$ outside of $\mathfrak{K}_2^\vee(y)$. 
Note that for any $y>0$, there exists an integer $L$ such that $T_y^{2, (d)} < T\eD_{\cD_L}$.

\begin{lem}
	\label{lem_Xdiff}
	For any $\ell\ge 1$ and $y_C\in (0, 1)$,
	there exists $y\in (0, y_C]$ and four real numbers $\cp_D, \cp_Z, \mpr_U, \mpr_D>0$ 
such that
	the two following inequalities  hold for any
	$d\in \II{2, \infty}$ and any
	$x\in \cH\eiD_\ell$:
	\begin{multline*}
		\idc{z\in \Y_2(2y)}
		\PR_x\Big(\pi_2(X\eiD(2))\in \RMd z
		\pv 2 <  T_y^{2, (d)} \wedge T\etD_{\mpr_U}\Big)
	\\	\ge \cp_D\, \idc{z\in \Y_2(2y)} 
	\PR_{z\eD}\Big(\pi_2(X\eiD(1))\in \RMd z
	\pv 1 <  T_y^{2, (d)}\wedge T\etD_{\mpr_D}\Big),
\end{multline*}
and:
\begin{equation*}
\cZ\eD:= \PR_{z\eD}\Big(X\eiD(1)\in \cR_y
\mVg 1 <  T_y^{2, (d)}\wedge T\etD_{\mpr_D}\Big)
\ge c_Z,
\end{equation*}
where the subset $\cR_y$ of $\Y_2(2y)$ is defined as follows:
\begin{equation}\label{eq_def_Ry}
	\cR_y := \Big\{z\in \Y_2(2y)\pv z_0\in (1-5y, 1-4y)\mVg z_1\in (2y, 3y)\Big\}.
\end{equation}
\end{lem}

\bpf: Let $\ell\ge 1$, $m_H$ and $y_H$ be the constants associated to $\cH\eD_\ell$
and $y_C$ be given.
Let $y = y_C\wedge (\sfrac{y_H\!}{2})\wedge (\tdiv{7})$.
Recall from Proposition \ref{M_Kcomp} 
that the system $(X\eiD_0, X\eiD_1, X\eiD_{(2)}) = \pi_2(X\eiD)$ 
is autonomous under $\PR\ejdD_x$,
whatever $d\in \II{2, \infty}$ and $x\in \cH\eD_\ell$,
with a common infinitesimal generator $\cL^{(2)}$ on $\cX_2$.
As stated in Lemma~\ref{lem_Har}, the process $\pi_2(X\eiD)$ satisfies Property~$(H)$
on $\Y_2(y)$.
Since  $\pi_2(\cH\eD_\ell)\subset \Y_2(2y)\subset \mathfrak{K}_2^\wedge(y)$,
we deduce as in the proof of Proposition \ref{M_Mixd}
that there exists $\cp_D^1>0$
such that the following inequality holds
for any $d\in \II{2, \infty}$ and $x\in \cH\eD_\ell$:
\begin{multline}\label{eq_min_dens_Ay}
		\idc{z\in \Y_2(2y)}
\PR\ejdD_x\big(\pi_2(X\eiD(2))\in \RMd z
\pv 2 <  T_y^{2, (d)}\big)
\\\ge \cp^1_D\, \idc{z\in \Y_2(2y)} 
\PR\ejdD_{z\eD}\big(\pi_2(X\eiD(2))\in \RMd z
\pv 1 <  T_y^{2, (d)}\big).
\end{multline}

To prepare for the second inequality,
recall \eqref{eq_def_Ry}.
The set $\cR_y$ has a non empty interior, 
so that we can find a smooth function $f:\mathfrak{K}_2^\vee(y)\mapsto [0, 1]$ 
with support on $\cR_y$ that is non-zero.
We are led to consider 
the corresponding Cauchy problem on $\bR_+\times \mathfrak{K}_2^\vee(y)$
with the value at the boundary given by the function $u_{\partial \mathfrak{K}_2^\vee(y)}$
defined as $u_{\partial \mathfrak{K}_2^\vee(y)}(0, z) = f(z)$ for any $z\in \mathfrak{K}_2^\vee(y)$
and as $u_{\partial \mathfrak{K}_2^\vee(y)}(t, z) = 0$ for any $t\ge 0$ and $z\in \partial \mathfrak{K}_2^\vee(y)$.
Thanks to Property~$(H)$, 
there exists a unique positive strong solution $u$ to this Cauchy problem.
Note that by construction  $\pi_2(z\eD) = z^{(2)}$ is independent of $d$.
As a consequence:
\begin{equation}\label{eq_min_Ry}
\begin{split}
	\PR\ejdD_{z\eD}\big(\pi_2(X\eD(1))\in \cR_y
\pv 1 <  T_y^{2, (d)}\big)
&\ge \bE\ejdD_{z\eD}\Big(f\Big[\pi_2(X\eD(1))\Big]
\pv 1 <  T_y^{2, (d)}\Big)
\\&\quad = u(1, z^{(2)})
>0.
\end{split}
\end{equation}
So as to relate to the original dynamics prescribed by $\PR_x$,
we need upper-bounds of the third moment 
that are given independently of $\pi_2(X\eD)$,
by referring to Subsection~\ref{M_sec_momCk}.
Note that $M\eF_3(x)\le m_\ell$ holds for any $x \in \cH\eD_\ell$
and that $T_y^{2, (d)} < \tau_0\ejdD$ hold also a.s.
Thanks to Lemma \ref{M_moment_controlk},
there exists $m_U'>0$
such that the following inequality holds for any $d\in \II{2, \infty}$
and any $x \in \cH\eD_\ell$:
	\begin{equation}\label{eq_min_M3U}
	\PR_x\ejdD\lp  T\etF_{m_U'}\le 2\wedge T_y^{2, (d)}
	\bv \cF^{(2)}\rp
	\le \tdiv{2}.	
\end{equation}
Similarly, 
there exists $m_D'>0$
such that the following inequality holds for any $d\in \II{2, \infty}$:
\begin{equation}\label{eq_min_M3D}
	\PR_{z\eD}\ejdD\lp  T\etF_{m_D'}\le 1\wedge T_y^{2, (d)}
	\bv \cF^{(2)}\rp
	\le \tdiv{2}.	
\end{equation}
Thanks to \eqref{eq_min_dens_Ay} and to \eqref{eq_min_M3U},
the following inequality holds for any $d\in \II{2, \infty}$
and any $x \in \cH\eD_\ell$:
\begin{multline*}
	\idc{z\in \Y_2(2y)}
\PR\ejdD_x\Big(\pi_2(X\eiD(2))\in \RMd z
\pv 2 <  T_y^{2, (d)}\wedge T\etF_{m_U'}\Big)
\\
\ge 2\cp^1_D\, \idc{z\in \Y_2(2y)} 
\PR\ejdD_{z\eD}\Big(\pi_2(X\eiD(1))\in \RMd z
\pv 1 <  T_y^{2, (d)}\Big).
\end{multline*}
Thanks to Proposition \ref{M_Girs},
noting that $M_3(x)\le M\eF_3(x) + 1$ (it holds for any $d\ge 2$ and $x\in \cX_d$),
there exists $\cp_D>0$ 
such that the following inequality holds for any $d\in \II{2, \infty}$
and any $x \in \cH\eD_\ell$, with $m_U = m_U' +1$ and $m_D = m_D' + 1$:
\begin{multline*}
\idc{z\in \Y_2(2y)}
\PR_x\Big(\pi_2(X\eiD(2))\in \RMd z
\pv 2 <  T_y^{2, (d)}\wedge T\etD_{m_U}\Big) 
\\
\ge \cp_D\, \idc{z\in \Y_2(2y)} 
\PR_{z\eD}\Big(\pi_2(X\eiD(1))\in \RMd z
\pv 1 <  T_y^{2, (d)}\wedge T\etD_{m_D}\Big),
\end{multline*}
which is the first inequality in Lemma \ref{lem_Xdiff}.

Similarly, thanks to Proposition \ref{M_Girs},
recalling \eqref{eq_min_Ry} and \eqref{eq_min_M3D},
there exists $c_Z>0$
such that the following inequality holds for any $d\in \II{2, \infty}$:
\begin{equation*}
	\PR_{z\eD}\Big(\pi_2(X\eiD(1))\in \cR_y
	\mVg 1 <  T_y^{2, (d)}\wedge T\etD_{\mpr_D}\Big)
	\ge c_Z,
\end{equation*}
which concludes the proof of Lemma \ref{lem_Xdiff}.
\epf

\subsubsection{The event of a regular behavior
	happen  with a lower-bounded probability}
\label{sec_Xcut}

Subsection~\ref{sec_Xcut} is devoted to the proof of the upcoming Lemma~\ref{M_Xcut},
in which we justify 
for well-prepared initial conditions
that the dependency in the different components of $X\eD_{(2)}$
can be forgotten in an event of non-negligible probability.

\begin{lem}
	\label{M_Xcut}
For any $\ell\ge 1$,
	there exists $y_C>0$
	such that the following property holds for any $y\in (0, y_C)$
	and any $m>0$.
	There exists $\cp>0$ 
	such that  the following inequality holds 
	for  any $d\in \II{2, \infty}$ 
 and any $x\in \cX_d$ such that 
	$\pi_2(x)\in \cR_y$ and $M_3(x)\le m$:	
	\begin{equation*}
		\PR_x\Big(\tau_0\ejdD < 1
		< T\eD_{\cD_3}\Big)
		\ge c.
	\end{equation*}
\end{lem}
\bpf:
	We wish to control the extinction of a uniform upper-bound 
of  $X\eD_{(2)}$, 
that is to be solution to the following SDE:
\begin{equation}
	\RMd Z(t):= \frac{\RMd t}{4}
	+ \sqrt{Z(t)\cdot (1-Z(t))}\, \RMd B_{(2)}(t)
	\mVg \quad Z(0) = y,
	\label{eq_J2_zDef1}
\end{equation}
where $W$ is a standard Brownian motion.
Note that 0 is an accessible boundary for this process $Z$
(it is actually regular refecting),
thanks to Corollary \ref{cor_bound_oneD}.

To ensure that the flux of population 
into $X\eD_{(2)}$
due to mutations remains lower than $\tdiv{4}$
in the time-interval $[0, \tpF]$ for some $\tpF>0$,
we wish to impose that  $X\eD_{1}$ 
remains upper-bounded by $\tdiv{4\lambda}$
on this time-interval.
For practical reasons, the considered upper-bound is actually slightly adjusted:
\begin{equation}
	y_M := \frac{1}{4\lambda} \wedge \frac{1}{2}.
	\label{eq_J2_yMDef}
\end{equation}
Similarly as for $X\eD_{(2)}$,
we thus consider an upper-bound for the process $X\eD_{1}$,
as the solution to the following SDE:
\begin{equation}
	\RMd Y(t):= (\lambda+\alpha)\, \RMd t
	+ \sqrt{Y(t)\cdot (1-Y(t))}\, \RMd B_1(t)
	\mVg \quad Y(0) = y_M/4,
	\label{eq_J2_YDef}
\end{equation}
where $B_1$ is a standard Brownian motion
(recall that $M_1\ejdD\le 2$ under $\PR\ejdD_x$).
We consider the following stopping time 
$T^Y_{y_M}$ 
(as a function of $y_M$):
\begin{equation}\label{eq_def_TYyM}
	T^Y_{y_M}
	:= \inf\Big\{t\ge 0\pv Y(t)\ge y_M\Big\},
\end{equation}
namely the hitting time of $y_M$ by the process $Y$.
We consider also the martingale process $\cN_1(t)$,
defined as follows for any $t>0$:
\begin{equation}
	\cN_1(t) 
	:= \int_{0}^t \Big[\sqrt{Y(t)\cdot (1-Y(t))} \wedge \sqrt{y_M\cdot (1-y_M)}\Big] \RMd B_1(t).
	\label{eq_def_J2_NmU}
\end{equation}
Its quadratic variation 
is upper-bounded as follows
at time $\tpF$:
\begin{equation*}
	\LAg \cN_1\RAg_{\tpF} 
	\le \tpF\cdot y_M\cdot (1-y_M).
\end{equation*} 
Doob's inequality thus entails the following upper-bound:
\begin{equation}\label{eq_J2_supNP}
	\PR\Big(\sup_{t\le \tpF} |\cN_1(t)| \ge \dfrac{y_M}{2}\Big)
	\le   \frac{16\tpF}{y_M}.
\end{equation}
Thanks to \eqref{eq_J2_YDef} and \eqref{eq_def_J2_NmU}
since $y_M \le \tdiv{2}$, 
the following identity holds a.s. on the event $\Lbr t \le T^Y_{y_M}\Rbr$, 
for any $t\in [0, \tpF]$:
\begin{equation*}
	Y(t) = \frac{y_M}{4} + (\lambda+\alpha)\cdot t+ \cN_1(t).
\end{equation*}
On the event $\{T^Y_{y_M} \le \tpF\}$,
the evaluation of this identity at time $T^Y_{y_M}$ entails 
the following inequality a.s.:
\begin{equation}
	\sup_{t\le \tpF} |\cN_1(t)| \ge \frac{3y_M}{4} - (\lambda+\alpha)\cdot \tpF.
	\label{ineq_J2_YN}
\end{equation}
\medskip

On the other hand, 
thanks to 
Property $(H)$ (see Subsection \ref{sec_Harnack})
and to Proposition~\ref{M_Kcomp}
since the event $\{\tF < T\eD_{\cD_2}\}$
only depends on the process $X\eD_0$,
there exists $c_0>0$ such that 
the following inequality holds
for any $d\in \II{2, \infty}$ and
for any $x\in \cX_d$ such that $x_0 \in [\tdiv{2}, 1]$:
\begin{equation}	\label{eq_tauT}
	\PR\ejdD_x \Big(\tF < T\eD_{\cD_2}\Big) \ge c_0.
\end{equation}
We then choose $\tpF>0$
as follows:
\begin{equation}\label{eq_def_J2_tF}
	\tpF := \Big(\frac{c_0\cdot y_M}{2^5}\Big)
	\wedge \Big(\frac{y_M}{4(\lambda+\alpha)}\Big)
	\wedge \tF,
\end{equation} 
which ensures thanks to \eqref{eq_J2_supNP} and to \eqref{ineq_J2_YN} 
the following inequality:
\begin{equation}
	\PR\Big(T^Y_{y_M} \le \tpF\Big) 
	\le \frac{c_0}{2}.
	\label{ineq_J2_TYyM}
\end{equation}
Recalling that 0 is an accessible boundary for the process $Z$
defined in \eqref{eq_J2_zDef1},
we then choose $z\in (0, y_M/4)$ sufficiently small 
for the following inequality to hold:
\begin{equation}
	\PR\Big(\tpF \le \ext^Z\Big) \le \frac{c_0}{4}, 
	\label{eq_J2_zDef}
\end{equation}
where $\ext^Z:= \inf\{t\ge 0\pv Z_t = 0\}.$

The definition of $y_C>0$ is as follows:
\begin{equation*}
y_C = \frac{z}{3}.
\end{equation*}
For any $y\in (0, y_C)$, $z\in \cR_y$ thus implies that
$z_0> 1-5y_C> 1-(\sfrac{5\!}{12})\cdot y_M>\tdiv{2}$
and $z_1\wedge z_2 < 3y_C = z< y_M/4$.
Thanks to Lemma \ref{lem_comp},
for any $y\in (0, y_C)$,
any $d\in \II{2, \infty}$ and any $x\in \cX_d$ such that $\pi_2(x)\in \cR_y$,
$X\eD_1$ under the law $\PR\ejdD_x$
 is upper-bounded by the process $Y$,
 which entails with \eqref{ineq_J2_TYyM} the following inequality:
 \begin{equation}\label{eq_T1DyM}
 \PR\ejdD_x\Big(T^{1, (d)}_{y_M} \le \tpF\Big) 
 \le \frac{c_0}{2},
 \end{equation}
where the stopping time $T^{1, (d)}_{y_M}$ is defined as follows:
\begin{equation*}
	T^{1, (d)}_{y_M}
	:= \inf\Big\{t\ge 0\pv X\eD_1\ge y_M\Big\},
\end{equation*}
namely the hitting time of $y_M$ by the process $X\eD_1$.

Thanks similarly to Lemma \ref{lem_comp},
$X\eD_{(2)}$ under the law $\PR\ejdD_x$
is upper-bounded by the process $Z$.
Recalling \eqref{eq_J2_zDef}  and \eqref{eq_T1DyM},
the following inequalities thus hold
for any $y\in (0, y_C)$,
any $d\in \II{2, \infty}$ and any $x\in \cX_d$ such that $\pi_2(x)\in \cR_y$:
\begin{equation}\label{eq_tau02D}
	\begin{split}
			\PR\ejdD_x\Big(\tau_0^{(2:d)} < \tpF\Big)
		&\ge \PR\ejdD_x\Big(\tau_0^{(2:d)} < \tpF< T^{1, (d)}_{y_M}\Big) 
	\\&\ge 1- \frac{3c_0}{4}.
	\end{split}
\end{equation}
Recalling \eqref{eq_tauT}, \eqref{eq_T1DyM} and \eqref{eq_tau02D},
the following inequality thus holds for any $y\in (0, y_C)$,
any $d\in \II{2, \infty}$ and any $x\in \cX_d$ such that $\pi_2(x)\in \cR_y$,
\begin{equation}
\PR\ejdD_x\Big(\tau_0^{(2:d)} < \tpF\wedge T\eD_{\cD_2}\Big) 
\ge \frac{c_0}{4}.
\end{equation}
Since the event $\{\tau_0^{(2:d)} < \tpF\wedge T\eD_{\cD_2}\}$ 
belongs to the sigma-field $\cF^{(2)}_{\tpF}$,
we can follow the same reasoning regarding the third moment
as in Subsection \ref{sec_X0diff}
and find a constant $\cp_1>0$
such that the following inequality holds
for any $y\in (0, y_C)$,
any $d\in \II{2, \infty}$ and any $x\in \cX_d$ such that both $\pi_2(x)\in \cR_y$
and $M_3(x)\le m$:
\begin{equation}
	\PR_x\Big(\tau_0^{(2:d)} < \tpF\wedge T\eD_{\cD_2}\Big) 
	\ge \cp_1.
\end{equation}

On the other hand, 
there exists $\cp_2>0$ such that the following inequality holds for any $d\in \II{2, \infty}$ and $x\in \cD\eD_2$:
\begin{equation*}
	\PR_x\Big(\tF< T\eD_{\cD_3}\Big) 
\ge \cp_2.
\end{equation*}
It is exactly \eqref{eq_min_X0} from Subsection~\ref{M_SHPos},
which extends to the case $d=\infty$.
Thanks to the strong Markov property at time $\tau_0^{(2:d)}$,
it concludes the proof of Lemma \ref{M_Xcut}
with $\cp = \cp_1\cdot \cp_2>0$.
\epf 

\subsubsection{Main mixing estimate}
\label{sec_MixdInf}

After the proofs of the three Lemmas~\ref{lem_Hdl}, \ref{lem_Xdiff} and \ref{M_Xcut}, 
we are in conditions to prove the mixing estimate,
as stated in the next theorem
that is the main result of this Subsection~\ref{sec_mixAcc}.
It exploits the following definition of $\zeta\eiD$.
In this formula, 
the two real numbers $y$ and $\mpr_D$ are associated by Lemma \ref{M_Xcut}
with the (arbitrary) choice  $\tp:=1$.
\begin{equation}
	\label{M_zeta}
	\zeta\eiD(\RMd x):= \int_{\cX_d}
	\PR_{z}\lp X\eD(1)\in \RMd x
	\bv \tau_0\ejdD < 1
	< T\eD_{\cD_3}\rp\, \nu\eD({\rm d} z)\mVg
\end{equation}
where  the measure $\nu\eD$ is defined as follows:
\begin{equation*}
\nu\eD({\rm d} z)
:= \dfrac{\idc{z\in \cR_{y}} }{\cZ\eD}\PR_{z\eD}\Big(\pi_2(X\eiD(2))\in \RMd z
\pv 1 <  T_y^{2, (d)}\wedge T\etD_{\mpr_D}\Big),
\end{equation*}
while $z\in \cR_y$ is seen as an element of $\cX_d$ 
by defining $z_k = 0$ for any $k\in \II{3, d}$.

Note that $\nu\eD$ is related to the measure that appears in the lower-bound in Lemma~\ref{lem_Xdiff},
with a restriction to the set $\cR_{y}\subset \Y_2(2y)$
followed by its renormalisation.
The second inequality in Lemma~\ref{lem_Xdiff}
ensures that the normalizing term $\cZ\eD$
is lower-bounded away from 0 uniformly in $d\in \II{2, \infty}$.

\begin{rem} 
As exploited in the final Subsection \ref{M_sec_TECV},
	the constraint $1< T\eD_{\cD_3}$
 ensures that $\zeta\eiD$ is supported on $\cD\eD_3$.
\end{rem}
\begin{theo}
	\label{M_Mix.dInf}
	For any $\ell\ge 1$, 
	there exist an integer $L>\ell$ and two real number $\tp, \cp>0$
	such that 
 the following inequality
holds  for any $d\in \II{2, \infty}$ and 
any	$x\in \cD\eiD_\ell$,
where $\zeta\eiD$ is defined in \eqref{M_zeta}:
	\begin{equation*}
	\PR_x\Big(X\eD(\tp)\in \RMd x'
		\pv \tp < T\eiD_{\cD_L}\Big)
		\ge c\, \zeta\eiD(\RMd x').
	\end{equation*}
\end{theo}

\bpf: Let $\ell\ge 1$.
Thanks to Lemma \ref{lem_Hdl}, 
the set $\cH\eD_\ell$, the two real numbers $\cp_H, \tp_H>0$
and the integer $L_H\ge 1$ are such that
 the following 
inquality holds  
for any
$d\in \II{2, \infty}$ and any
$x\in \cD\eiD_\ell$:
\begin{equation}\label{eq_Hdl}
	\PR_x\Big(X\eiD(\tp_H)\in \cH\eD_\ell
	\pv \tp_H < T\eD_{\cD_{L_H}}\Big)
	\ge \cp_H.
\end{equation}
 We then define $y_C>0$ according to Lemma~\ref{M_Xcut}.
 Thanks to Lemma~\ref{lem_Xdiff},
there exists $y\in (0, y_C]$ and four real numbers $\cp_D, \cp_Z, \mpr_U, \mpr_D>0$ 
such that
the two following inequalities  hold for any
$d\in \II{2, \infty}$ and any
$x\in \cH\eiD_\ell$:
\begin{multline}\label{eq_min_Xdiff}
	\idc{z\in \Y_2(2y)}
	\PR_x\Big(\pi_2(X\eiD(2))\in \RMd z
	\pv 2 <  T_y^{2, (d)} \wedge T\etD_{\mpr_U}\Big)
	\\	\ge \cp_D\, \idc{z\in \Y_2(2y)} \PR_{z\eD}\Big(\pi_2(X\eiD(1))\in \RMd z
	\pv 1 <  T_y^{2, (d)}\wedge T\etD_{\mpr_D}\Big),
\end{multline}
and:
\begin{equation}\label{eq_Zmin}
	\cZ\eD:= \PR_{z\eD}\Big(\pi_2(X\eiD(1))\in \cR_y
	\mVg 1 <  T_y^{2, (d)}\wedge T\etD_{\mpr_D}\Big)
	\ge \cp_Z.
\end{equation}
Note that there exists an integer $L\ge 3\vee L_H$ 
such that $T_y^{2, (d)} < T\eD_{\cD_{L}}$
(recall \eqref{eq_def_Ty2d}).
Thanks to Lemma~\ref{M_Xcut}
with this value of $y$ and $m_C = m_D\vee m_U$,
there exists $\cp_C>0$ 
such that  the following inequality holds 
for  any $d\in \II{2, \infty}$ 
and any $x\in \cX_d$ such that 
$\pi_2(x)\in \cR_y$ and $M_3(x)\le m_C$:	
\begin{equation}\label{eq_def_cC}
	\PR_x\Big(\tau_0\ejdD < 1
	< T\eD_{\cD_3}\Big)
	\ge \cp_C.
\end{equation}
Thanks to Proposition \ref{M_Kcomp}
since $X\eD_j(\tau_0\ejdD)= 0$ for any $j\ge 2$,
the following identity holds for any $d$ and any $x, x'\in \cX_d$
such that $\pi_2(x) = \pi_2(x')$:
\begin{multline*}
\PR\ejdD_x\Big(X\eD(\tau_0\ejdD)\in \RMd z\mVg
\tau_0\ejdD \in \RMd t\pv \tau_0\ejdD < 1\wedge T\eD_{\cD_3}\Big)
\\= \PR\ejdD_{x'}\Big(X\eD(\tau_0\ejdD)\in \RMd z\mVg
\tau_0\ejdD \in \RMd t \pv \tau_0\ejdD < 1\wedge T\eD_{\cD_3}\Big).
\end{multline*}
As in the proof given in Subsection \ref{sec_X0diff},
this entails that there exists $\cp_G>0$
such that the following inequality 
holds for any $d$ and any $x, x'\in \cX_d$
 such that $\pi_2(x) = \pi_2(x')$, $M_3(x)\le m_U$ and $M_3(x)\le m_D$:
\begin{multline*}
	\PR_x\Big(X\eD(\tau_0\ejdD)\in \RMd z\mVg
	\tau_0\ejdD \in \RMd t
	\pv \tau_0\ejdD < 1\wedge T\eD_{\cD_3}\Big)
	\\ \ge \cp_G \PR_{x'}\Big(X\eD(\tau_0\ejdD)\in \RMd z\mVg
	\tau_0\ejdD \in \RMd t 
	\pv \tau_0\ejdD < 1\wedge T\eD_{\cD_3}\Big).
\end{multline*}
Thanks to the strong Markov property at time $\tau_0\ejdD$,
the following inequality 
thus holds for any $d$ and any $x, x'\in \cX_d$
such that $\pi_2(x) = \pi_2(x')$, $M_3(x)\le m_U$ and $M_3(x)\le m_D$:
\begin{equation*}
	\PR_x\Big(X\eD(1)\in \RMd z\pv \tau_0\ejdD < 1< T\eD_{\cD_3}\Big)
	\ge \cp_G \PR_{x'}\Big(X\eD(1)\in \RMd z \pv \tau_0\ejdD < 1< T\eD_{\cD_3}\Big).
\end{equation*}

Thanks to the Markov property at time $1$,
recalling \eqref{M_zeta}, \eqref{eq_min_Xdiff}, \eqref{eq_Zmin} and \eqref{eq_def_cC},
we deduce
the following inequality for any
$d\in \II{2, \infty}$ and any
$x\in \cH\eiD_\ell$:
	\begin{equation*}
	\PR_x\Big(X\eD(3)\in \RMd x'
	\pv 3 < T\eiD_{\cD_{L_D}}\Big)
	\ge (\cp_D\cdot \cp_Z\cdot \cp_C\cdot \cp_G)\cdot \zeta\eiD(\RMd x').
\end{equation*}
Thanks to  the Markov property at time $t_H$, recalling \eqref{eq_Hdl} and that $L\ge L_H$,
with  $\cp = c_H\cdot \cp_D\cdot \cp_Z\cdot \cp_C\cdot \cp_G>0$
and $\tp = t_H +3$,
the following inequality thus holds 
for any $d\in \II{2, \infty}$ and 
any	$x\in \cD\eiD_\ell$:
\begin{equation*}
	\PR_x\Big(X\eD(\tp)\in \RMd x'
	\pv \tp < T\eiD_{\cD_L}\Big)
	\ge c\cdot \zeta\eiD(\RMd x').
\end{equation*}
The proof of Theorem \ref{M_Mix.dInf}
is complete.\epf

\begin{rem}
	In the above proof of Theorem \ref{M_Mix.dInf},
 Lemma~\ref{M_moment_controlk} was exploited to deduce
	upper-bounds of the third moments. 
	Instead of third moments, we could have concluded to a similar result
	by considering second moments instead.
	However, the proof of Theorem~\ref{thm_AF.dInf} below
	strongly exploits a control of a moment with exponent strictly greater than 2.
\end{rem}

\subsection{Almost perfect harvest}
\label{sec_AFdI}

The crucial result of this subsection is the upcoming Theorem~\ref{thm_AF.dInf},
whose conclusion is very close to the property $(A3_F)$ 
of ``almost perfect harvest",
see Subsection~\ref{M_sec_SeCo}.
\begin{theo}
	\label{thm_AF.dInf}
	Given any $\rho,  \mpr, \eta, y> 0$ 
	and any $\epsilon\in (0,\, 1)$,
	there exists  $\cp_V >0$ and an integer $J\ge 1$
	such that the following property holds 
	for any  $d\in \II{J, \infty}$ and any $x \in \cX_d$.
	There exists two stopping times $U_H\eD$ and $V\eD$ 
	such that the following three conditions hold
	for any $x_\zeta \in \cD_3\eD$:
	\begin{equation*}
		\begin{split}
			&  \Lbr\ext\eD \wedge \tF < U_H\eD \Rbr
			= \Lbr U_H\eD = \infty\Rbr\mVg
			\qquad  
			\PR_{x} (U_H\eD = \infty, \,  \tF< \ext\eD) 
			\le \epsilon\, e^{-\rho}
			\\&
			\PR_{x} \lp X\eD(U_H\eD) \in \RMd y
			\pv U_H\eD < \ext\eD \rp
			\le c_V\,\PR_{x_\zeta} \lp
			X\eD(V\eD)  \in \RMd y
			\pv V\eD < \ext\eD\rp.
		\end{split}
	\end{equation*} 
In addition, the probability space $\Omega$ 
and the filtration $\cF_t$ according to which $U_H\eD$ and $V\eD$ are stopping times
can be chosen to be the canonical representation 
of the process $X\eD$,
see Remark~\ref{rem_Omega_path}.
\end{theo}
\begin{rem}\label{rem_D3}
	The definition of $\zeta\eD$ in \eqref{M_zeta} 
	makes it supported on $\cD_3$.
	To emphasize that this property is sufficient for our purposes,
	we consider for the upper-bound in the last inequality
 Dirac initial conditions 
	of the form $x_\zeta\in \cD_3$.
\end{rem}
The proof of Theorem~\ref{thm_AF.dInf} is split into three parts:
we start in Subsection~\ref{sec_param} with the choice of the parameters and the statement of the corresponding properties,
then introduce the definition of $U_H\eD$ and $V\eD$ with their intrinsic properties in  Subsection~\ref{sec_def_UH}
before we conclude with the comparison of densities at time $U_H\eD$ versus $V\eD$
in Subsection~\ref{sec_Comp_Dens}.

\subsubsection{Choice of the parameters}
\label{sec_param}

The choice of $t_F = 1$ is made for simplicity.
In the first time-interval of length $\tdiv{3}$, 
we justify the access to a suitable set $E\eD(m, y, \eta)$ 
for $\eta$ well-chosen as a function of $\epsilon$, according to \eqref{eq_etaD},
then $m, 1/y$ sufficiently large according to Theorem~\ref{eT_dInf}.
In the next time-interval of length $t_H \le \tdiv{3}$, 
we couple the first $J$ coordinates 
between two processes with different initial conditions.
The duration $t_H$ is to be fixed below (in \eqref{def_THprop}),
sufficiently small to restrict on trajectories away from the boundaries.
Then, we impose that $X\eD_{(J)}$ gets extinct in the next time-interval 
of length $\tpF$,
with a similar approach as in Subsection~\ref{sec_Xcut}
for Lemma~\ref{M_Xcut}. 
Remark that we rely for the two last time-intervals
on the notations and results of Subsection~\ref{M_sec_Agg}.
Notably,  we  recall that the process $X\eD$ is solution under
$\PR\ejD$ to \eqref{Skd},
that $X\eD_{(J)}$ is defined in \eqref{eq_def_XJ} 
and that the process $M_1\ejD$ is upper-bounded by $J$
by virtue of its definition in \eqref{eq_def_MJD}. 

	We wish to control the extinction of a uniform upper-bound $(Z(t))_{t\ge 0}$
of the process $\big(X\eD_{(J)}(\tau_E\eD+t_H+t)\big)_{t\ge 0}$, 
that is to be solution to the following SDE,
analogous to \eqref{eq_J2_zDef1},
with $t\in [0, \tpF]$:
\begin{equation}
	\RMd Z(t):= \frac{\RMd t}{4}
	+ \sqrt{Z(t)\cdot (1-Z(t))}\, \RMd B(t)
	\mVg \quad Z(0) = z,
	\label{eq_zDef1}
\end{equation}
where $B$ is a standard Brownian motion.
Note that 0 is an accessible boundary for this process $Z$,
thanks to Corollary \ref{cor_bound_oneD}.

\paragraph{Upper-bound on the incomming flux of population:}\hfill\\
To ensure that the flux of population 
into $X\eD_{(J)}$
due to mutations remains lower than $\tdiv{4}$
in the time-interval of length $\tpF$ that follows $\tau_E\eD+t_H$,
we wish to impose that the process $\big(X\eD_{J-1}(\tau_E\eD+t_H +t)\big)_{t\in [0, \tpF]}$ 
remains upper-bounded by $\tdiv{4\lambda}$.
As in \eqref{eq_J2_yMDef}, we define:
\begin{equation}
	y_M := \frac{1}{4\lambda} \wedge \frac{1}{2}.
	\label{eq_yMDef}
\end{equation}
Similarly as in \eqref{eq_J2_YDef},
we thus consider an upper-bound $(Y(t))_{t\ge 0}$ 
for the process \\$\big(X\eD_{J-1}(\tau_E\eD+t_H +t)\big)_{t\in [0, \tpF]}$,
as the solution to the following SDE
with $t\in [0, \tpF]$:
\begin{equation}
	\RMd Y(t):= (\lambda+\alpha)\; \RMd t
	+ \sqrt{Y(t)\cdot (1-Y(t))}\, \RMd B_{(-1)}(t)
	\mVg \quad Y(0) = y_M/4,
	\label{eq_YDef}
\end{equation}
where $B_{(-1)}$ is a standard Brownian motion
(recall that $M_1\ejD\le J$
and that $X\eD_{J-1}$ is solution to \eqref{Skd} under $\PR\ejD_{x}$).
Let $T^Y_{y_M}$ denote the hitting time of $y_M$ by this process $Y$
(posterior to $t_H$).
The definition of the martingale $\cN_{(-1)}(t)$ 
is then slightly adapted from \eqref{eq_def_J2_NmU}
with a start at time $t_H$ with value~$0$,
so that the two following inequalities hold,
generally for the first one and a.s.
on the event $\{T^Y_{y_M} \le \tpF\}$ for the second one:
\begin{align}\label{eq_supNP}
	\PR\Big(\sup_{t\le \tpF} |\cN_{(-1)}(t)| \ge y_M/2\Big)
	\le   \frac{16 \tpF}{y_M},
\\	\sup_{t\le \tpF} |\cN_{(-1)}(t)| \ge \frac{3y_M}{4} - (\lambda+\alpha)\cdot \tpF.
	\label{ineq_YN}
\end{align}

\paragraph{Lower-bound on the survival with an initial condition in $\cD_3$}
\hfill

Thanks to Lemma \ref{lem_comp},
for any initial condition $x_\zeta\in \cD_3$,
$X\eiD_0$ 
is lower-bounded by the solution $Y_0$ to the following SDE:
\begin{equation*}
	\RMd Y_0(s) = -\lambda\, \RMd s
	+ \sqrt{Y_0(s)\cdot (1-Y_0(s))}\, \RMd B_0(s)\mVg \quad Y_0(0) = \tdiv{6}. 
\end{equation*}
Thus, denoting $c_\zeta:= (\tdiv{2})\cdot\PR_{\tdiv{6}}(\inf_{t\in [0, \tdiv{3}]} Y_0(t) >0) >0$,
we deduce that the following inequality holds for any $d\in \II{2, \infty}$
and any $x_\zeta\in \cD\eD_3$:
\begin{equation}\label{M_Svz}
	\PR_{x_\zeta}(\tdiv{3}<\ext\eiD) \ge 2 c_\zeta.
\end{equation}

\paragraph{Choice of $\tpF$ and $z$:} \hfill\\
Let $\epsilon>0$,
that we assume without loss of generality to be smaller than $c_\zeta$. 
We then choose $\tpF>0$
as follows:
\begin{equation}\label{eq_def_tF}
\tpF := \Big(\frac{\epsilon\cdot y_M}{16}\Big)
\wedge \Big(\frac{y_M}{4(\lambda+\alpha)}\Big)
\wedge \Big(\frac{1}{3}\Big),
\end{equation} 
which ensures thanks to \eqref{eq_supNP} and \eqref{ineq_YN} 
the following upper-bound in probability:
\begin{equation}
\PR\Big(T^Y_{y_M} \le \tpF\bv Y(0) = y_M/4\Big) 
\le \epsilon.
\label{ineq_TYyM}
\end{equation}
Recalling that 0 is an accessible boundary for the process $Z$
defined in \eqref{eq_zDef1},
we then choose $z\in (0, y_M/4)$ sufficiently small 
for the following inequality to hold:
\begin{equation}
	\PR\Big(\tpF \le \ext^Z\bv Z(0) = z\Big) \le \epsilon, 
	\label{eq_zDef}
\end{equation}
where $\ext^Z:= \inf\{t\ge 0 \pv Z_t = 0\}.$

\paragraph{Choice of $\eta$, $\mpr$ and $y$:}\hfill\\
Now, with the constants $C_G, C_M$ associated by Proposition \ref{M_Girs}
with the choices of $k = 3$, $\eps>0$ and $t=\tF$,
we can choose $\eta>0$ as follows:
\begin{equation}
	\eta = \Big( \frac{z}{C_M}\Big) \wedge \Big( \frac{\log(2)}{C_G}\Big),
	\label{eq_etaD}
\end{equation} 
where we recall that $z$ is an implicit function of $\epsilon$.
Thanks to Theorem \ref{eT_dInf},
given $\rho >0$,
we can choose 
the two real numbers $\mpr, y>0$ such that
the following inequality holds 
for any  $d\in \II{2, \infty}$ and  any $x\in \cX_d$:
\begin{equation}\label{eq_def_Ed}
	\PR_x\big( \tdiv{3} < \tau\eD_E\wedge \ext\eD\big)
	\le \eps,
\end{equation}
where we recall the definition of $E\eD$ given in \eqref{M_cDE}:
\begin{equation*}
	E\eD= E\eD(\mpr, y, \eta):= \big\{x\in \cX_d\pv 
	M_3(x) \le \mpr\mVg 
	\frl{j\le \Lfl \sfrac{\mpr\!}{\eta}\Rfl+1} x_j\ge y\big\}.
\end{equation*}
Recalling \eqref{M_Svz} and that $c_\zeta\ge \epsilon$, 
\eqref{eq_def_Ed}
has the following implication for 
any $d\in \II{2, \infty}$ and any initial condition $x_\zeta \in \cD_3$:
\begin{equation}\label{min_xzeta}
\PR_{x_\zeta}(\tau\eD_E \le \tdiv{3}<\ext\eiD) \ge c_\zeta.
\end{equation}

With $m_M = C_M\cdot \mpr$,
recalling the definition of $C_M$ in relation to Proposition \ref{M_Girs}, 
we deduce that the following inequalities hold
for any  $d\in \II{2, \infty}$ and 
any $x\in E\eD$:
\begin{equation}\label{ineq_t3m}
	\Bv \log\Big(\frac{\RMd \PR\ejD_x}{\RMd \PR_x}\! 
	\left\vert \underset{[0, 1]}{\quad}\right. 
	\Big)\Bv
	\le \frac{C_G\cdot \mpr}{J} 
\le \log(2),\qquad 
	\PR_x\Big(T\etD_{m_M} \le \tF\Big)
	\le \eps,
\end{equation}
with the following definition of the integer $J$:
\begin{equation}
	J:= \Lfl \sfrac{\mpr\!}{\eta}\Rfl+2.
	\label{eq_def_J}
\end{equation}
Since $\eta \le \sfrac{z}{C_M}$,
recalling the definition of $T\etD_m$ from \eqref{eq_def_TkM},
the following inequalities hold a.s. on the event $\{\tF < T\etD_{m_M}\}$,
for any $t_H \le \tdiv{3}$ and any $d\in \II{J, \infty}$:
\begin{equation}\label{ineq_XR0}
	X\eD_{(J-1)}(t_H) \le \frac{C_M\cdot \mpr}{(J-1)^3}\le z.
\end{equation}
Recalling that $z< y_M/4$, 
we deduce that $X\eD_{J-1}(t_H) \le y_M/4$ and that $X\eD_{(J)}(t_H)\le z$,
as intended for $Y$  (see \eqref{eq_YDef})
and for $Z$ (see \eqref{eq_zDef1})
to be appropriate upper-bounds.
\medskip

\paragraph{Choice of $t_H$:}\hfill\\
Note that
$\pi_J(x) 
\in \Y_J(y)$
holds for any $d\in \II{J, \infty}$ and any $x\in E\eD$, 
with the following definition of $\Y_J(y)$:
\begin{align*}
	\Y_J(y) 
	:= \textstyle \Lbr x\in \cX_J
	\pv \lp\bigwedge_{\{i\in \II{0, J}\}} x_i\rp  > y\Rbr,
\end{align*}
and the projection $\pi_J$ defined in \eqref{M_piK}.
On $\Y_J(y)$,  the diffusion term in the system  $(S^{(J)})$ of SDEs, 
i.e. the system $\eqref{Sd}$ with $d$ replaced by $J$,
is uniformly elliptic. In practice,
we need a bit more space for Property~$(H)$ to hold
(see Lemma~\ref{lem_Har}),
so that we consider the following exit time $T\eJD_{y'}$ generally for any $y'\in (0, y)$:
\begin{equation}
	T\eJD_{y'}:= \inf\Big\{t\ge 0\pv 
	\pi_J(X\eD(t))\notin \Y_J(y')\Big\}
	< \ext\eD.
	\label{def_Ty}
\end{equation}
The probability of such an escape is required to be very small,
uniformly in $d\in \II{J, \infty}$ and in $x\in E\eD$, 
as stated in the upcoming Lemma~\ref{lem_tH},
where  $\PR^{(J)}$ denotes the law of the system given by $(S^{(J)})$:
\begin{lem}
	\label{lem_tH}
The following supremum tends to 0
	as $t_H$ tends to 0: 
	\begin{align*}
		\sup\Lbr\PR^{(J)}_{\pi_J(x)}\Big(T\eJD_{y/2} \le t_H\Big)\bv
		d \in \II{J, \infty},\; x\in E\eD\Rbr.
	\end{align*}
\end{lem}
\bpf: Since the system $(S^{(J)})$ is uniformly elliptic on 
	some connected open relatively compact subset $\mathfrak{K}_J$ of $\Y_J(y/2)$
	with $\cC^\infty$-boundary
	that contains $\Y_J(y)$,
and recalling Proposition \ref{M_Kcomp},
Lemma~\ref{lem_tH} is deduced thanks to Lemma~\ref{lem_Har}
(also thanks e.g. to \cite[Proposition V.2.5]{Bass}).
\epf
\\

Thanks to Proposition \ref{M_Kcomp},
we can thus choose $t_H \le \tdiv{3}$ sufficiently small such that
the following inequality holds 
for any $d \in \II{J, \infty}$, and any $x\in E\eD$:
\begin{equation}
	\PR\ejD_x\Big(T\eJD_{y/2} \le t_H\Big) 
	\le \eps.
	\label{def_THprop}
\end{equation}

\subsubsection{Definition of $U_H\eD$ with a control of exceptional events}
\label{sec_def_UH}
For the definition of $U_H\eD$,
the following extinction time $\widehat\tau_0\ejD$
is considered for the process $X\eD_{(J)}$  after time $\tau_E\eD + t_H$:
$$\widehat\tau_0\ejD:=\inf\Big\{t\ge \tau_E\eD + t_H \pv X\eD_{(J)}(t) =0\Big\}.$$

In view of Theorem \ref{thm_AF.dInf},
we define
$U_H\eD:= \widehat\tau_0\ejD$
on the following event: 
\begin{equation}
\Big\{\tau_E\eD < \tdiv{3}\Big\}
\cap \Big\{\tau_E\eD + t_H < \wht T\eJD_{y/2}\Big\}
	\cap 
	\Big\{\widehat\tau_0\ejD < (\tau_E\eD + t_H+\tpF)\wedge  \wht T\etD_{m_M}\wedge \ext\eD\Big\}
	\label{eq_def_UH}
\end{equation}
and otherwise $U_H\eD:= \infty$,
where the definitions of $\wht T\eJD_{y/2}$ and $\wht T\etD_{m_M}$
are  as follows, 
adjusted from those of $T\eJD_{y/2}$ and $T\etD_{m_M}$
with a specific time-shift:
\begin{equation*}
\begin{split}
	\wht T\eJD_{y/2}
	&:= \inf\Big\{t\ge \tau_E\eD\pv 
\pi_J(X\eD(t))\notin \Y_J(y/2)\Big\}, 
\\
\wht T\etD_{m_M}
&:= \inf\Big\{t\ge \tau_E\eD + t_H\pv 
M_3\eD(t)\ge m_M\Big\}.
\end{split}
\end{equation*}
Since $t_H, \tpF\le \sfrac{1\!}{3}$, we deduce that $\big\{\tF\wedge \ext\eD < U_H\eD\big\} = \big\{U_H\eD= \infty\big\}$.
On the other hand, the stopping time $V\eD$ is defined as follows:
\begin{equation}
V\eD 	:=\inf\Big\{t\ge \tau_E\eD + 2 t_H \pv X\eD_{(J)}(t) =0\Big\}.
\label{eq_def_V}
\end{equation}
Remark that $U_H\eD$ and $V\eD$ are regularly expressed in terms of the process $X\eD$,
so that they can be expressed as stopping times 
for a path space representation of $\Omega$ and $(\cF_t)$,
namely the canonical representation of $X\eD$,
as stated in Theorem~\ref{thm_AF.dInf}.

Thanks to the strong Markov property at time 
$\tau_E\eD$,
the following inequality holds for any $d\in \II{J, \infty}$
and $x\in \cX_d$:
\begin{equation}\label{eq_new_gH}
	\PR_x\Big(U_H\eD = \infty,  \tF < \ext\eD\Big)
\le \PR_x \Big( \tdiv{3} <  \tau_E\eD \wedge \ext\eD \Big)
+ \bE_x\Big[g\big(X\eD(\tau_E\eD)\big)\pv \tau_E\eD<\ext\eD\Big],
\end{equation}
where the function $g$
is expressed as follows 
 for any  $x\in E\eD$:
\begin{multline}
g(x) 
= \PR_x\Big(T\etD_{m_M}\le \tF \Big)
+ \PR_x\Big(T\eJD_{y/2} \le t_H\mVg \tF < T\etD_{m_M} \Big)
\\+ \PR_x\Big(t_H + \tpF < \wtd \tau_0\ejD\mVg \tF < T\etD_{m_M} \Big),
\end{multline}
where $\wtd\tau_0\ejD:=\inf\Big\{t\ge t_H \pv X\eD_{(J)}(t) =0\Big\}.$
Recalling \eqref{ineq_t3m},
the function $g$ is upper-bounded as follows in terms of the probability measure $\PR\ejD_x$
for any $x\in E\eD$:
\begin{multline}\label{eq_def_gH}
	g(x) 
	\le \PR_x\Big(T\etD_{m_M}\le \tF \Big)
	+ 2\PR\ejD_x\Big(T\eJD_{y/2} \le t_H \Big) 
	+ 2\PR\ejD_x\Big(t_H + \tpF < \wtd \tau_0\ejD \Big).
\end{multline}
Thanks to the Markov Property at time $t_H$
and to Lemma \ref{lem_comp}
with \eqref{ineq_XR0},
the following upper-bound for the last term
is expressed in terms of the processes $Y$ (see \eqref{eq_YDef})
and $Z$ (see \eqref{eq_zDef1})
for any $x\in E\eD$:
\begin{equation*}
\PR\ejD_x\Big(t_H + \tpF < \wtd \tau_0\ejD \Big)
\le \PR\Big(T^Y_{y_M} \le \tpF\bv Y(0) = y_M/4\Big)
+ \PR\Big(\tpF \le \ext^Z\bv Z(0) = z\Big).
\end{equation*}
Recalling \eqref{ineq_TYyM}
and \eqref{eq_zDef},
this term is thus upper-bounded by $2\eps$.
Injecting this upper-bound in \eqref{eq_def_gH}
and the similar ones deduced from \eqref{ineq_t3m} and \eqref{def_THprop},
we deduce that $g(x)\le 7\epsilon$
for any $x\in E\eD$.
Injecting this upper-bound 
and the one from \eqref{eq_def_Ed}
into \eqref{eq_new_gH},
the inequality $\PR_x\Big(U_H\eD = \infty,  \tF < \ext\eD\Big)\le 8\eps$
holds for any $d\in \II{J, \infty}$
and $x\in \cX_d$.
To arrive at the upper-bound 
that $\PR_x(U_H\eD = \infty,  \tF < \ext\eD)\le \epsilon'\cdot e^{-\rho}$,
as stated in Theorem \ref{thm_AF.dInf} for some $\epsilon'>0$,
we just have to choose $\epsilon = \epsilon'\cdot e^{-\rho}/8$.
Since $\epsilon$ is freely chosen, so is $\epsilon'$.

\subsubsection{Comparison of densities}
\label{sec_Comp_Dens}

The core argument for this Subsection~\ref{sec_Comp_Dens} is Harnack's inequality 
(from Subsection \ref{sec_Harnack})
applied on a reduction of the system 
to a finite dimensional projection. 
We need however to handle carefully both the distance from the boundary 
of the finite dimensional projection,
since we need the diffusion to be elliptic,
and the control of the third moment,
to convert $\PR$ into $\PR\ejD$ and vice-versa.

For any $d\in \II{J, \infty}$,
let us first assume that $x$ and $x_\zeta$ both belong to  $E\eD$,
so that $\tau_E\eD= 0$ for both initial conditions $x$ and $x_\zeta$.
By virtue of
the definition of $U_H\eD$ (see \eqref{eq_def_UH}), 
then thanks to \eqref{ineq_t3m}:
\begin{equation}
	\begin{split}
		\PR_x\Big(X\eD(U_H\eD) \in \RMd x'\mVg U_H\eD < \infty\Big)
		&\le 2\, \PR\ejD_x\Big(X\eD(\widehat\tau_0\ejD) \in \RMd x'\mVg U_H\eD <\infty\Big).
	\end{split}
	\label{eq_dens_UHG}
\end{equation}
Recall the definition of $\tau_0\ejD$ from \eqref{eq_def_tau0ejD}.
Thanks to the Markov property at time $t_H$, 
the following inequality holds:
\begin{equation}
	\PR\ejD_x\Big(X\eD(\widehat\tau_0\ejD) \in \RMd x'\mVg U_H\eD <\infty\Big)
	\le \int_{\cX_d} \nu^1_x(\RMd x_H) \nu^2_{x_H}(\RMd x'),
	\label{eq_kern}
\end{equation}
where:
\begin{equation*}
	\nu^1_x(\RMd x_H) 
	:= \PR\ejD_x\Big(X\eD(t_H) \in \RMd x_H 
	\pv t_H < T\eJD_{y/2}\wedge T\etJ_{m_M}\Big),
\end{equation*}
so that $\nu^1_x(\RMd x_H)$ describes the transition of the process $X\eD$
during the time-interval $[0, t_H]$
on the event $\{t_H < T\eJD_{y/2}\wedge T\etJ_{m_M}\}$,
while:
\begin{equation}
	\nu^2_{x_H}(\RMd x')
	:= \PR\ejD_{x_H}\Big( X\eD(\tau_0\ejD) \in \RMd x'
	\pv \tau_0\ejD < \tpF \wedge T\etJ_{m_M}\wedge \ext\eD\Big),
	\label{eq_def_nu2}
\end{equation}
so that  $\nu^2_{x_H}(\RMd x')$ describes the transition of $X\eD$
during the time-interval $[0, \tau_0\ejD]$
on the event $\{\tau_0\ejD < \tpF \wedge T\etJ_{m_M}\wedge \ext\eD\}$.
Note  on this event that $X\eD(\tau_0\ejD)_i = 0$ holds for any $i\ge J$,
so that $X\eD(\tau_0\ejD)$ is directly expressed in terms of $\pi_J(X\eD(\tau_0\ejD))$.
The measure
$\nu^2_{x_H}$ only depends on $\pi_{J}(x_H)$
thanks to Proposition~\ref{M_Kcomp},
so that the formula $\bar \nu^2_{\pi_{J}(x_H)} = \nu^2_{x_H}$
produces a well-defined quantity.
Therefore:
\begin{equation}
	\int_{\cX_d} \nu^1_x(\RMd x_H) \nu^2_{x_H}(\RMd x')
	= \int_{\cX_{J}} \bar\nu^1_x(\RMd x'_H) \bar\nu^2_{x'_H}(\RMd x'),
	\label{eq_proj}
\end{equation}
where $\bar\nu^1_x$ is the image of $\nu^1_x$ by the projection
$\pi_{J}: \cX_d\mapsto \cX_{J}$, i.e.:
\begin{equation*}
	\bar\nu^1_x(\RMd x'_H) 
	:= \PR\ejD_x\Big(\pi_{J}(X\eD(t_H)) \in \RMd x'_H 
	\pv t_H < T\eJD_{y/2}\wedge T\etJ_{m_M}\Big).
\end{equation*}
\medskip

Note that $T\eJD_{y/2}\wedge T\etJ_{m_M}$
corresponds to the exit time of $\pi_J(X\eD)$
out of some domain $\mathfrak{H}_J(y/2, m_M)$.
There exist two connected open relatively compact sets $\mathfrak{K}_J^\wedge, \mathfrak{K}_J^\vee$
with $\cC^\infty$-boundaries such that the following inclusions hold:
\begin{eqnarray*}
	\mathfrak{H}_J(y/2, m_M)
	\subset \mathfrak{K}_J^\wedge,
	\quad 
\overline{\mathfrak{K}_J^\wedge}	\subset\mathfrak{K}_J^\vee
	\subset \mathfrak{H}_J(y/4, 2 m_M).
\end{eqnarray*}
Considering Property $(H)$
(see Subsection \ref{sec_Harnack})
with Dirichlet boundary conditions on $\mathfrak{K}_J^\vee$
(that is with $u_{\partial \mathfrak{K}_J^\vee}(z, t) \equiv 0$
for any $z\in \partial \mathfrak{K}_J^\vee$ and $t\in [0, t_H]$)
amounts to the following inequality, which holds with a fixed constant $C_H>0$
for any $d\in \II{J, \infty}$ and any $x_1, x_2 \in \mathfrak{K}_J^\wedge $:
\begin{multline*}
	\idc{x'\in \mathfrak{K}_J^\wedge}
	\PR\ejD_{x_1}\Big(\pi_J(X\eD(t_H)) \in \RMd x'\mVg t_H < T_{\mathfrak{K}_J^\vee}\Big)
	\\ \textstyle
	\le C_H\, \idc{x'\in \mathfrak{K}_J^\wedge}
	\PR\ejD_{x_2}\Big(\pi_J(X\eD(2t_H)) \in \RMd x'\mVg 2t_H < T_{\mathfrak{K}_J^\vee}\Big).
\end{multline*}
It implies in particular that 
$\bar\nu^1_x\le C_H\cdot \check\nu^1_{x_\zeta}$,
where the time $t_H$ and the constants $y/2$ and $m_M$ are slightly relaxed:
\begin{equation*}
	\check\nu^1_{x_\zeta}(\RMd x'_H) 
	:= \PR\ejD_{x_\zeta}\Big(\pi_J(X\eD(2 t_H)) \in \RMd x'_H 
	\pv 2 t_H < T\eJD_{y/4}\wedge T\etJ_{2m_M}\Big).
\end{equation*}
Thanks in addition to \eqref{eq_dens_UHG}, \eqref{eq_kern}, \eqref{eq_proj},
we obtain as an intermediate step the following inequality,
which holds for any $d\in \II{J, \infty}$ and any $x, x_\zeta\in E\eD$:
\begin{equation}
	\PR_x\Big(X\eD(U_H\eD) \in \RMd x'\mVg U_H\eD < \infty\Big)
	\le 2 C_H \int_{\cX_J} \check\nu^1_{x_\zeta}(\RMd x'_H) \bar\nu^2_{x'_H}(\RMd x').
	\label{eq_aH}
\end{equation} 
\medskip

To go back to an upper-bound in terms of the original process,
we want again to exploit  Proposition \ref{M_Girs}.
So we need to again ensure upper-bounds on the third moments 
for the last components for which we lost the information.
We exploit the representation given in Proposition \ref{M_Kcomp},
firstly on the time-interval $[0, 2t_H]$.
Since $x_\zeta \in E\eD$, we have
$M_3\eF(0)\le \mpr/y$.
Thanks to Lemma \ref{M_moment_controlk},
we can thus define $m_H>0$ such that
the following inequality holds 
a.s. on the event $\{2 t_H < T\eJD_{y/4}\}$,
for any $d\in \II{J, \infty}$ and any $x_\zeta \in E\eD$:
\begin{equation*}
	\PR_{x_\zeta}\ejD\lp  T\etF_{m_H}\le 2 t_H 
	\bv \cF^{(J)}\rp
	\le \tdiv{2}.
\end{equation*}
With the fact that $\ext\eD> T\eJD_{y/4}$,
it implies that $\check\nu^1_{x_\zeta}\le 2 \tilde\nu^1_{x_\zeta}$,
where:
\begin{equation}
	\tilde\nu^1_{x_\zeta}(\RMd x_H)
	:= \PR\ejD_{x_\zeta}\Big(\pi_J\big(X\eD(2 t_H)\big) \in \RMd x'_H 
	\pv 2 t_H < \ext\eD\wedge T\etJ_{2m_M}\wedge T\etF_{m_H}\Big).
	\label{eq_def_tnu1}
\end{equation}
Recalling that the formula $\bar \nu^2_{\pi_J(x_H)} = \nu^2_{x_H}$
produces a well-defined quantity,
it implies that $\int_{\cX_J} \tilde\nu^1_{x_\zeta}(\RMd x'_H) \bar\nu^2_{x'_H}(\RMd x')
= \int_{\cX_d} \hat\nu^1_{x_\zeta}(\RMd x_H) \nu^2_{x_H}(\RMd x')$, where:
\begin{equation}
	\hat\nu^1_{x_\zeta}(\RMd x_H)
	:= \PR\ejD_{x_\zeta}\Big(X\eD(2 t_H) \in \RMd x_H 
	\pv 2 t_H < \ext\eD\wedge T\etJ_{2m_M}\wedge T\etF_{m_H}\Big).
	\label{M_m3tH}
\end{equation}
Note that this measure $\hat\nu^1_{x_\zeta}$
is supported on the set $\big\{x_H\in \cX_d\pv M_3\eF(x_H)\le m_H\big\}$.
Thanks again to Lemma \ref{M_moment_controlk},
we can thus define $m_F\ge m_H$ such that
the following inequality holds 
a.s.
for any $d\in \II{J, \infty}$ and any $x_H \in \cX_d$ 
such that $M_3\eF(x_H)\le m_H$:
\begin{equation*}
	\PR_{x_H}\ejD\lp  T\etF_{m_F}< \tpF\wedge \tau_0\ejD 
	\bv \cF^{(J)}\rp
	\le \tdiv{2}.
\end{equation*}
It entails that $\nu^2_{x_H}\le 2 \tilde\nu^2_{x_H}$,
where:
\begin{equation}
	\tilde\nu^2_{x_H}(\RMd x')
	:= \PR\ejD_{x_H}\Big(X\eD(\tau_0\ejD) \in \RMd x' 
	\pv \tau_0\ejD <\tpF \wedge \ext\eD\wedge T\etJ_{m_M}\wedge T\etF_{m_F}\Big).
	\label{eq_def_tnu2}
\end{equation}
Since $M\eD_3(x) \le M\ejD_3(x) + M\eF_3(x)$,
and thanks to the Markov property at time $t_H$,
recalling \eqref{eq_aH}, \eqref{eq_def_tnu1}, \eqref{M_m3tH} and \eqref{eq_def_tnu2},
the following inequality holds
with $V\eD$ as defined in \eqref{eq_def_V}
and $\hat m:= [2m_M+m_H]\vee [m_M+m_F]$,
for any $x, x_\zeta \in E\eD$:
\begin{equation*}
	\PR_x\Big(X\eD(U_H\eD) \in \RMd x'\mVg U_H\eD < \infty\Big)
	\le 8 C_H \cdot 
	\PR_{x_\zeta}\ejD\Big(X\eD(V\eD) \in \RMd x'
	\pv V\eD < \ext\eD\wedge T\etD_{\hat m}\Big).
\end{equation*}
We can then relate to the original law $\PR_{x_\zeta}$
thanks to Proposition~\ref{M_Girs}
with an additional factor $C_G>0$
such that  the following inequality holds for any $x, x_\zeta \in E\eD$:
\begin{equation}\label{ineq_UV}
	\PR_x\Big(X\eD(U_H\eD) \in \RMd x'\mVg U_H\eD < \infty\Big)
	\le 8 C_H\cdot C_G\cdot 
	\PR_{x_\zeta}\Big(X\eD(V\eD) \in \RMd x'
	\pv V\eD < \ext\eD\Big).
\end{equation}
To conclude, we consider general initial conditions $\bar x\in \cX_d$ and $\bar x_\zeta \in \cD_3$
and define:
\begin{equation*}
\nu^E_x(\RMd x') := \PR_x\Big(X\eD(\tau_E\eD)\in \RMd x'\pv \tau_E\eD <\tdiv{3}\wedge \ext\eD\Big).
\end{equation*}
We deduce  from \eqref{min_xzeta} that $\nu^E_{\bar x_\zeta}(E\eD) \ge c_\zeta$
while $\nu^E_{\bar x}(E\eD) \le 1$.
Thanks (twice) to the strong Markov property at time $\tau_E\eD$
and by virtue of the definitions of $U_H\eD$ and $V\eD$
in \eqref{eq_def_UH} and \eqref{eq_def_V},
\eqref{ineq_UV} is extended as follows 
for any $d\in \II{J, \infty}$
to any $\bar x\in \cX_d$ and $\bar x_\zeta \in \cD_3$:
\begin{equation*}
\begin{split}
	&\PR_{\bar x}\Big(X\eD(U_H\eD) \in \RMd x'\mVg U_H\eD < \infty\Big)
	= \int_{E\eD} \nu^E_{\bar x}(\RMd x) \PR_{x}\Big(X\eD(U_H\eD) \in \RMd x'\mVg U_H\eD < \infty\Big)
	\\&\quad \le 8 C_H\cdot C_G\cdot \frac{\nu^E_{\bar x}(E\eD)}{\nu^E_{\bar x_\zeta}(E\eD)}
	\int_{E\eD} \nu^E_{\bar x_\zeta}(dx_\zeta) 
	\PR_{x_\zeta}\Big(X\eD(V\eD) \in \RMd x'
	\pv V\eD < \ext\eD\Big)
	\\&\quad \le c_V\cdot  \PR_{\bar x_\zeta}\Big(X\eD(V\eD) \in \RMd x'
	\pv V\eD < \ext\eD\Big),
\end{split}
\end{equation*}
where $c_V = 8 C_H\,C_G/c_\zeta>0$.
This concludes the proof of Theorem \ref{thm_AF.dInf}.
\epf

\subsection{Concluding the proof of Theorem \ref{M_ECVdInf}}
\label{M_sec_TECV}

We plan to exploit Theorem \ref{M_th:ECV}
and ensure Assumption $\mathbf{(AF)}$
(recall Subsection~\ref{M_sec_SeCo}).
Firstly, the sets $\cD\eD_\ell$  satisfy Assumption $(A0)$
(recall \eqref{def_DL}).
Thanks to Theorem \ref{M_Mix.dInf}, 
Assumption $(A1)$ 
holds true for 	the reference measure $\zeta\eD$ defined by \eqref{M_zeta}.
Thanks to  $(A1)$ and \cite[Lemma 3.0.2]{AV_QSD}, 
$\rho_S$ is upper-bounded 
by some real number $\wtd{\rho}_S$
that only depends on the constants involved in $(A1)$,
and can thus be chosen independent of $d$.
In order to satisfy  $\rho>\rho_S$,
we set $\rho:= 2 \wtd{\rho}_S$.
For this choice of $\rho$,
we deduce Assumption~$(A2)$ 
as a consequence of Theorem~\ref{eT_dInf}
(recall also Remark~\ref{rem_eTdInf}),
where the complementary to the transitory domain can be taken as 
the proposed set $E\eD(m, y, 1)$ for some $m, y>0$
independent of $d\in \II{1, \infty}$,
or simply $\cD\eD_L$ with $L\ge \tdiv{y}$
(we choose $\eta = 1$ for simplicity).


When $d=\infty$
we infer Assumption $(A3_F)$ thanks to Theorem~\ref{thm_AF.dInf}
(regardless of the transitory domain that is chosen),
for this choice of $\rho$ and $\zeta\eI$,
noting that $\zeta\eI$ is supported on $\cD\eI_3$.
Assumption $\mathbf{(AF)}$ thus holds true for $d=\infty$.
Thanks to   Theorem~\ref{M_th:ECV},
we then  conclude the proof 
that the semi-group displays QSC
and that the Q-process exists on $\cX_\infty$.

In addition, recall that $\rho\eI_0= -\log[\PR_{\nu\eI}(1<\ext\eI)]$ 
(see Remark \ref{rem_nuChar}).
Since $\cup_\ell \cD\eI_\ell = \cX_\infty$
(recall \eqref{def_DL}),
there exists $\ell$ such that $\nu\eI(\cD\eI_\ell) >0$.
With a slight adaptation of the proof of Lemma \ref{M_LX0i},
we deduce that $\PR\eI_x(\ext\eI \le 1)$
is uniformly bounded away from 0 for any $x\in \cD\eI_\ell$.
This concludes that $\rho\eI_0>0$.

As noted in \cite[Subsection 6.1]{AV_disc},
the statement of Assumption $(A3_F)$ is actually required for a single value of $\eps$
defined in terms of the parameters involved in Assumptions $(A1)$ and $(A2)$.
Since these parameters can be chosen uniformly in $d$,
the corresponding value of $\eps$ can also be chosen uniformly
so that the required properties in Assumption $\mathbf{(AF)}$ holds
uniformly for any $d$ larger than a given threshold
(see the definition of $J$ in \eqref{eq_def_J}).
In addition and for the same reason,
all parameter 
involved in the convergences
can be chosen 
independently of $d$ sufficiently large. 
There are indeed intricate yet explicit relations
between all these parameters introduced in \cite{AV_disc}
and the corresponding assumptions.
This concludes the proof of Theorem~\ref{M_ECVdInf}.
\epf
\hfill\\

No datasets were generated or analysed during the current study.

\appendix
\section*{Appendix}
\section{Two representations for the same process}
\label{app_TwoRep}
In this Section \ref{app_TwoRep}, 
we justify more precisely (see  Proposition \ref{prop_Trep}) than we did in Remark~\ref{rem_martRep}
that the process we consider is exactly the same as defined
notably in \cite{AP13},
though the representation has been adapted for our purposes.

Let us recall from \eqref{Sd}
for any $d\in \II{2, \infty}$
the definition of our focal process $X\eD$ as follows
for any  $i\in \II{0, d}$:
\begin{equation}
		\RMd X\eD_i(t) = b\eD_i(X\eD(t))\, \RMd t
		+ \RMd \cN\eD_i(t),
		\label{eq_def_XdM}
\end{equation}
where the vectorial function $b$ encodes the drift term
 as follows for any $x\in \cX_d$ and $i\in \II{0, d}$:
\begin{equation*}
b\eD_i(x)	
:= \alpha\cdot \Big(\Tsum{j\in \II{0, d}} j\cdot x_j - i\Big)\cdot x_i
+ \lambda\cdot \Big(x_{i-1}-\idc{i<d} x_{i}\Big),
\end{equation*}
while the martingale process $\cN\eD_i$, which starts at 0 for $t=0$, is expressed as follows:
	\begin{equation}
\begin{split}
		\RMd \cN\eD_i(t)
	&:= \sum_{j = 0}^{d}
	\Big(\delta_{ij} - X\eD_i(t)\Big)\cdot \sqrt{X\eD_j(t)} \, \RMd W_j(t)
	\\&	=\sqrt{X\eD_i(t)} \, \RMd W_i(t)
	- X\eD_i(t) \, \RMd W\sD(t)
\end{split}
	\label{eq_def_MDi}
\end{equation}
in terms of 
a family $(W_i)_{i\in \Z_+}$  of mutually independent Brownian motions
and an aggregated martingale $W\sD$ with the following definition:
\begin{equation*}
\RMd	W\sD(t):=  \textstyle \sum_{j=0}^{d} \sqrt{X\eD_j(t)} \RMd W_j(t),\quad 
W\sD(0)=0.
\end{equation*}

On the other hand, we next consider the more classical way that Muller's ratchet diffusion $\wht X$ is defined,
notably in \cite{AP13}:
\begin{equation*}
	\RMd \wht X\eD_i(t) = b\eD_i\big(\wht X\eD(t)\big)\, \RMd t
	+ \RMd \wht \cN\eD_i(t),
\end{equation*}
where  the martingale term $\wht \cN\eD_i$ is now expressed as follows:
	\begin{equation*}
	\RMd \widehat \cN\eD_i(t) = \Tsum{j\neq i} \sqrt{\widehat X\eD_i(t)\cdot \widehat X\eD_j(t)} \, \RMd W_{i, j}(t),
	\quad \widehat \cN\eD_i(0) = 0.
\end{equation*}
in terms of 
a family $(W^{i, j})_{i<j}$ of independent Brownian motions,
extended to any $i\neq j$
by the symmetry property $W_{i, j}(t) = W_{j, i}(t)$.

For the infinite-dimensional case, we recall the definition of $\cX^\eta$ in \ref{def_cXeta}
as the set of probabilities on $\N$ with finite $\eta$-th moment
($\eta =6$ being considered for simplicity in our proofs).
The upcoming Proposition~\ref{prop_Trep} summarises the result of this Section~\ref{app_TwoRep}
and notably implies Proposition~\ref{prop_ex_uniq}.
\begin{prop}
	\label{prop_Trep}
The processes $(X\eD_i)$ and $(\widehat X\eD_i)$ share the same law.
Existence and uniqueness in law holds,
for processes on $\cX_d$ for any $d\in \N$, 
and, for the case $d=\infty$, on $\cX^\eta$ for any $\eta>2$.
For any $d$, 
the associated infinitesimal generator takes the following nondivergence form, for any $u\in C^2(\cD\eD_\ell)$ and $x \in \cD\eD_\ell$:
\begin{equation*}
	\cL^{(d)} u(x)
	= \frac{1}{2}\sum_{i, j = 1}^d \sigma^{(d)}_{i, j}(x) \partial^2_{i, j} u(x) 
	+ \sum_{i= 1}^d b^{(d)}_{i}(x) \partial_{i} u(x),
\end{equation*}
where the state-dependent diffusion matrix $\sigma^{(d)}$ is expressed as follows:
\begin{equation*}
	\sigma^{(d)}_{i, j}(x) :=\left\{
	\begin{split}
&- x_i\cdot x_j\qquad &\text{ if } i\neq j,
\\&  x_i\cdot (1-x_i)\qquad &\text{ if } i= j,
	\end{split}
\right.
\end{equation*}
and the state-dependent drift term $b^{(d)}$ is expressed as follows:
\begin{equation*}
b^{(d)}_{i}(x)
:= \alpha\cdot\Big(M\eD_1(x) -i\Big)\cdot x_i
+ \lambda\cdot\Big(x_{i-1} - \idc{i<d}\cdot x_i\Big). 
\end{equation*}
\end{prop}
\begin{rem}
Recall that our representation of the martingale term 
corresponds to the idea that internal demographic fluctuations $W^i$
specific to each subpopulation level $i$
get compensated thanks to the process $W$.
On the other hand, 
the classical representation is associated 
to the principles of Moran models,
where demographic fluctuations are generated 
by choosing one parent at random whose offspring replaces 
another inidividual also chosen at random. 
By distinguishing these contributions in terms
of the two involved types $(i, j)$
(be it the parental type or the killed type),
we can derive the contribution
$\sqrt{\widehat X\eD_i(t)\cdot \widehat X\eD_j(t)} \, \RMd W_{i, j}(t)$
to the fluctuations.
It explains why $W_{i, j} = - W_{j, i}$ and 
how the previous term can be interpreted as a random flux, 
which does not need to be compensated.
\end{rem}

\bpf: According to \cite[Theorem 2.1]{Shi87},
the existence and the uniqueness in law of solutions hold
for the system \ref{eq_def_XdM} of SDEs
(for both finite and infinite $d$)
provided that the drift term $b_i\eD$ satisfies the following 
three conditions: 
\begin{itemize}
	\item $b_{i}(x)\ge 0$
	uniformly in $x\in \cX_d$,
	\item $\sum_{i\in \II{0, d}}b_{i}(x)=0$,
	\item there exists a matrix $(q_{ij})_{i,j\in \II{0, d}}$
	such that $q_{ij}\ge 0$
	and $\sup_{k\in \II{0, d}}\sum_{i=0}^d
	q_{ik}<+\infty$
	for every $i$ and $j$ of $\II{0, d}$,
	and such that the following inequality holds 
	for any $x$ and $x’$ of $\cX_d$ $(i\in \II{0, d})$: 
$$|b\eD_{i}(x)-b\eD_{i}(x’)| \le \Tsum{j\in S}q_{ij}\cdot|x_{j}-x_{j}’|.$$
\end{itemize}
The first two properties are satisfied for any $d\in \II{2, \infty}$.
For finite $d\in \N$,
the following estimate on the drift term holds 
for any $i\in \II{0, d}$ and $x,x'\in \cX_d$:
\begin{equation*}
|b_i(x) - b_i(x')|
\le \textstyle (2\alpha \cdot d + \lambda)\cdot \sum_{j=0}^d |x_j - x_j'|.
\end{equation*}
The corresponding choice of $q_{i j} \equiv 2\alpha \cdot d + \lambda$
ensures for any $d\in \N$ 
the existence and uniqueness in law for the system \ref{eq_def_XdM} of SDEs.
Thanks to \cite[Theorem 3]{AP13}
(which relies on the above result of $\cite{Shi87}$),
 existence and uniqueness in law
also hold 
on the set $\cX^\eta$ for any $\eta>2$.
The associated infinitesimal generator
is expressed in \cite[Theorem 1.1]{Shi87}
as stated in Proposition~\ref{prop_Trep}.

For completeness and to check that the laws of $\cN\eD_i(t)$
and $\widehat \cN\eD_i(t)$ actually coincide,
we provide next the computations of the quadratic variations 
of these martingale terms.
We first consider the quadratic variation of each component $i$
 for the process $\cN\eD_i$:
\begin{equation*}
	\begin{split}
		\RMd \LAg \cN\eD_i\RAg_t 
		&= X\eD_i(t) \,\RMd \LAg W_i\RAg_t + (X\eD_i(t))^2\, \RMd \LAg W\sD\RAg_t
		- 2 \sqrt{X\eD_i(t)}\cdot X\eD_i(t)\, \RMd \LAg W_i, W\sD\RAg_t
		\\&=  X\eD_i(t)\cdot (1-X\eD_i(t))\, \RMd t,
	\end{split}
\end{equation*}
then for the process $\widehat \cN\eD_i$:
\begin{equation*}
	\begin{split}
		\RMd \LAg \widehat \cN\eD_i\RAg_t		
		&= \Tsum{j\neq i} \widehat X\eD_i(t)\cdot \widehat X\eD_j(t) \, \RMd t
		\\&= \widehat X\eD_i(t)\cdot (1-\widehat X\eD_i(t))\, \RMd t.
	\end{split}
\end{equation*}

We secondly consider the cross-variations for $i\neq j$,
which yields for the process $M\eD$:
\begin{equation*}
	\begin{split}
		&\RMd \LAg \cN\eD_i, \cN\eD_j\RAg_t 
		\\&= \sqrt{X\eD_i(t)\cdot X\eD_j(t)} \, \RMd \LAg W_i, W_j\RAg_t
		-  \sqrt{X\eD_i(t)}\cdot X\eD_j(t) \, \RMd \LAg W_i, W\sD\RAg_t
		\\&
		\hcm{1}- X\eD_i(t)\cdot \sqrt{X\eD_j(t)} \, \RMd \LAg W\sD, W_j\RAg_t
		+X\eD_i(t)\cdot X\eD_j(t) \, \RMd \LAg W\sD\RAg_t
		\\&= \lc 0 - \sqrt{X\eD_i(t)}\cdot X\eD_j(t) \cdot\sqrt{X\eD_i(t)} 
		- X\eD_i(t)\cdot \sqrt{X\eD_j(t)} \cdot\sqrt{X\eD_j(t)} 
		+ X\eD_i(t)\cdot X\eD_j(t)\rc \, \RMd t,
		\\&= - X\eD_i(t)\cdot X\eD_j(t)\,  \RMd t,
	\end{split}
\end{equation*}
then yields for the process $\wht M\eD$:
\begin{equation*}
	\begin{split}
		\RMd \LAg \widehat \cN\eD_i, \widehat \cN\eD_j\RAg_t
		&= - \Tsum{k\neq i}\Tsum{\ell \neq j} 
		\sqrt{\widehat X\eD_i(t)\cdot \widehat X\eD_k(t)\cdot \widehat X\eD_j(t)\cdot \widehat X\eD_\ell(t)}\,
		\RMd \LAg W_{i, k}, W_{\ell, j}\RAg_t
		\\&= -\widehat X\eD_i(t)\cdot \widehat X\eD_j(t)\, \RMd t,
	\end{split}
\end{equation*}
since $\RMd \LAg W_{i, k}, W_{\ell, j}\RAg_t = \idc{i = \ell, j=k} \RMd t$.

Since the quadratic variations are expressed in the same ways 
for the solutions $X\eD$ and $\widehat X\eD$ to the systems respectively \ref{eq_def_XdM} and \ref{eq_def_MDi},
we conclude thanks to \cite[Theorem 1.1]{Shi87}
that the laws of these two processes coincide.
This concludes the proof of Proposition~\ref{prop_Trep}.

\section{Justification for the boundedness of the moments}
\label{sec_Fin_mom}

In this Section \ref{sec_Fin_mom}, 
we derive from the proof of \cite[Theorem 3]{AP13}
the local boundedness of the moment processes, as stated in this last Proposition~\ref{pr_Fin_mom}.
\begin{prop}
\label{pr_Fin_mom}
The process $M\eI_k$ is a.s. locally upper-bounded
for any $k \ge 1$, any $\alpha, \lambda>0$ 
and any $x\in \cX_\infty \cap \cX^{k}$
 under both $\PR_x$ and $\PR\ejD_x$,
 whatever $J\in \N$.
\end{prop}
\bpf: Without loss of generality, we assume that $k\ge 3$
and consider any $x\in \cX^k$, any $\lambda >0$
and any $t>0$.
The conclusion of the proof of \cite[Proposition 2.2]{AP13}
can be restated as the following upper-bound,
which holds in the case where $\alpha = 0$
(as indicated with the notation $\PR^{[0]}$):
\begin{equation}
\PR^{[0]}_x\Big(\Tsup{s\le t} M\eI_k(s) < \infty\Big) = 1.
\label{eq_locB_Mk}
\end{equation}
We consider any $\alpha>0$
and denote for clarity $\PR^{[\alpha]}$ instead of $\PR\eI$ 
for the law of the process $X\eI$ solution to \eqref{Sd}.
Thanks to \cite[Subsection 2.1, proof of Theorem 3]{AP13},
$\PR^{[\alpha]}$ is related to $\PR^{[0]}$ 
with a Girsanov transform,
namely $\PR^{[\alpha]}_x\vert_{\cF_t} = Z_\alpha(t) \cdot \PR^{[0]}_x\vert_{\cF_t}$
holds for any $x\in \cX_\infty$ and any $t>0$,
which involves the following martingale process $(Z_\alpha(s))_{s\ge 0}$:
\begin{equation*}
\begin{split}
	Z_\alpha(s) 
&:= \exp\lp -\alpha \M\eI_1(s) - (\alpha^2/2)\cdot \int_0^s M\eI_2(s) \RMd s\rp,
\\&\quad 
\le \exp\lp - \alpha\cdot[M\eI_1(x) + \lambda s]\rp.
\end{split}
\end{equation*}
In the above expression, $\M\eI_1(s)$ is the martingale component of $M\eI_1(s)$,
as stated in Lemma \ref{lem_MkSM}.
Since $Z_\alpha(t)$ is uniformly  upper-bounded,
\eqref{eq_locB_Mk} entails that
$\Tsup{s\le t} M\eI_k(s) < \infty$ holds also a.s. under $\PR^{[\alpha]}_x$,
provided $x\in \cX^k$.

We next consider $J\in \N$ and the law $\PR\ejaD_x
= \PR^{[\alpha], (J:\infty)}_x$
as stated in Section \ref{M_sec_Agg}.
Thanks to Proposition \ref{M_Kcomp},
and with the same arguments as above,
$\PR\ejaD$ is related to $\PR^{[0]}$ 
with another Girsanov transform,
namely $\PR\ejaD_x\vert_{\cF_t} = Z^J_\alpha(t) \cdot \PR^{[0]}_x\vert_{\cF_t}$
holds for any $x\in \cX_\infty$ and any $t>0$,
which involves the following martingale process $(Z^J_\alpha(s))_{s\ge 0}$:
\begin{equation*}
	\begin{split}
		Z^J_\alpha(s) 
		&:= \exp\lp -\alpha \M^{(J:\infty)}_1(s) - (\alpha^2/2)\cdot \int_0^s M^{(J:\infty)}_2(s) \RMd s\rp,
		\\&\quad 
		\le \exp\lp - \alpha\cdot\big[(M^{(J:\infty)}_1(x)\wedge J) + \lambda s\big]\rp.
	\end{split}
\end{equation*}
Similarly, $\M^{(J:\infty)}_1(s)$ is the martingale component of $M^{(J:\infty)}_1(s)$.
With the same arguments as for the proof of Lemma \ref{lem_MkSM},
its quadratic variation involves the process $M^{(J:\infty)}_2$,
which is the second moment saturated at value $J$, i.e.:
\begin{equation*}
\begin{split}
	M^{(J:\infty)}_2(t) 
&:= \sum_{i\le J-1} i^2\cdot X\eI_i(t) + J^2\cdot X\eI_{(J)}(t)
\\&= \sum_{i\ge 0} \big(i\wedge J\big)^2\cdot X\eI_i(t).
\end{split}
\end{equation*}
Since $Z^J_\alpha(t)$ is uniformly  upper-bounded,
\eqref{eq_locB_Mk} entails that
$\Tsup{s\le t} M\eI_k(s) < \infty$ holds a.s. under $\PR\ejaD_x$,
provided $x\in \cX^k$.
This concludes the proof of Propostion \ref{pr_Fin_mom}.\epf

%
%

%
%
%

%
\end{document}